\RequirePackage{fix-cm}
\documentclass[smallextended]{svjour3}       
\smartqed  
\usepackage{graphicx}
\usepackage{epstopdf, epsfig}
\usepackage{amsmath}
\usepackage{mathtools}
\usepackage{subfigure}
\usepackage[export]{adjustbox}
\usepackage{booktabs,array}
\usepackage{tikz}
\usepackage{multirow}
\usepackage{pstool}
\usepackage{psfrag}
\usepackage{pictex}
\usepackage{pifont}
\usepackage{xcolor}
\usepackage[shortlabels]{enumitem}
\usepackage{algorithm}
\usepackage{algorithmic}
\usetikzlibrary{shapes}
\usetikzlibrary{shapes.misc}
\usepackage{float}
\usepackage{amssymb}
\usepackage{bm}

\newcommand\Rey{\mbox{\textit{Re}}}
\newcommand\Real{\mbox{Re}}
\newcommand\Imag{\mbox{Im}}
\newcommand{\bluedot}{\tikz{\draw[blue,fill=blue] (0,0) circle (.4ex)}}

\DeclarePairedDelimiter\norm{\lVert}{\rVert}

\begin{document}

\title{Resolvent-based approach for $H_2$-optimal estimation and control: an application to the cylinder flow
}

\titlerunning{Resolvent-based $H_2$-optimal estimation and control}        

\author{Bo Jin \and Simon J.  Illingworth \and Richard D. Sandberg 
}


\institute{Bo Jin \at
              Department of Mechanical Engineering, University of Melbourne, VIC, 3010, Australia \\
              \email{bjin1@student.unimelb.edu.au}           
           \and
           Simon J.  Illingworth \at
             Department of Mechanical Engineering, University of Melbourne, VIC, 3010, Australia  
            \and
            Richard D. Sandberg \at
            Department of Mechanical Engineering, University of Melbourne, VIC, 3010, Australia  
}

\date{Received: date / Accepted: date}

\maketitle

\begin{abstract}
We consider estimation and control of the cylinder wake at low Reynolds numbers. A particular focus is on the development of efficient numerical algorithms to design optimal linear feedback controllers when there are many inputs (disturbances applied everywhere) and many outputs (perturbations measured everywhere).  We propose a resolvent-based iterative algorithm to perform i) optimal estimation of the flow using a limited number of sensors; and ii) optimal control of the flow when the entire flow is known but only a limited number of actuators are available for control. The method uses resolvent analysis to take advantage of the low-rank characteristics of the cylinder wake and solutions are obtained without any model-order reduction. Optimal feedback controllers are also obtained by combining the solutions of the estimation and control problems.  We show that the performance of the estimators and controllers converges to the true global optima, indicating that the important physical mechanisms for estimation and control are of low rank.
\keywords{vortex dynamics \and resolvent analysis \and flow control}
\end{abstract}

\section{Introduction}\label{intro}
Flow control is either passive or active, and can be of tremendous benefit in a number of applications \cite{gad2007flow}. A robust and energetically efficient way to control fluid flows is active closed-loop control or feedback control, which comprises actuators operating according to real-time information from the flow field obtained by sensors. In the last few decades, efforts have been made to improve the efficiency and effectiveness of feedback control schemes and linear control theory (e.g.~optimal control) has received considerable attention, providing a rigorous and systematic approach for designing control laws \cite{choi2008control}. However, direct implementation of linear control theory is computationally challenging for most fluid-flow problems due to the high dimensionality of the discretized Navier-Stokes equations. 
\subsection{Reduced-order modelling}
One common approach for solving flow control problems is to form reduced-order models (ROMs) which capture the essential mechanisms in fluid flows. Traditional control design tools can then be directly applied to design low-dimensional controllers. For fluid flows dominated by vortices, e.g.~bluff body flows, a simple way of obtaining ROMs is to use vortex models \cite{li2003feedback,protas2004linear,henning2007feedback}. Another alternative approach to build ROMs is based on the proper orthogonal decomposition (POD) \cite{lumley1970stochastic}, the standard method of balanced truncation \cite{moore1981principal} or the eigensystem realization algorithm (ERA) \cite{juang1985eigensystem}. In particular POD has shown its effectiveness for flow control of the cylinder wake \cite{gillies1998low,singh2001optimal} and of channel flow \cite{ilak2008modeling}. The work of Rowley \cite{rowley2005model} further extended balanced truncation to high-dimensional systems where approximate balanced truncation can be achieved using a snapshot-based algorithm named balanced proper orthogonal decomposition (BPOD). The ERA is a system identification method that is equivalent to BPOD and allows direct application to flow systems using only DNS or experimental data \cite{belson2013feedback,illingworth2016model,flinois2016feedback,yao2017feedback,yao2017model} without requiring simulations of the adjoint system. These methods were originally designed for stable systems of large dimension and were subsequently extended to unstable systems either by a state-projection method \cite{ahuja2010feedback} or by limiting the sampling time \cite{flinois2016feedback}. It was also shown that a frequency-domain version of BPOD was directly amenable to unstable systems \cite{dergham2011model}. For a deeper insight into these methods, a detailed comparison is presented by Ma et al. \cite{ma2011reduced} using the example of the flow past an inclined flat plat.

A fundamental goal of reduced-order modelling is to seek a low-dimensional description of a fluid flow. A direct application is the work of Dahan et al. \cite{dahan2012feedback} where frequency responses were captured from a series of numerical simulations under harmonic actuation over a broad range of frequencies. Recently, model decomposition based on resolvent analysis has shown good potential for control purposes. The method has been applied to efficiently characterise the input-output behaviour of a broad range of flows, such as cavity flow \cite{gomez2016reduced}, flat-plate boundary layer flow \cite{sipp2013characterization}, pipe flow \cite{mckeon2010critical} and the flow past a cylinder \cite{symon2018non}. The method allows direct formulation of the flow system in the frequency domain. Therefore no special treatment is required for unstable systems and there are no restrictions on the frequency range over which the linear dynamics can be captured. 

\subsection{Full-dimensional control}
Feedback control for two- or three-dimensional fluid flows will typically require some form of model reduction \cite{bagheri2009input}.
In the design procedure for full-dimensional control, although no direct model reduction is employed, the number of inputs or outputs (i.e.~terminals) in the control problem is generally limited to a small number. This is true for most control problems, even multi-input–multi-output (MIMO) systems, where the number of terminals $m$ is still far less than the dimensionality of the flow system $N$. This improves the feasibility of a number of algorithms for designing full-dimensional controllers directly from the original system. 

Methods that bypass open-loop model reduction for high-dimensional control problems are summarised in Bewley et al. \cite{bewley2016methods}, who propose several Riccati-less control design methods. In particular, the Minimum Control Energy (MCE) algorithm developed by Bewley et al. \cite{bewley2007minimal} uses eigenanalysis for full-dimensional optimal control problems under the limitation of an infinite control penalty. In later research, Pralits and Luchini \cite{pralits2010riccati} circumvented this limitation using the \textit{adjoint of the direct-adjoint} (ADA) algorithm which extends adjoint-based optimisation by considering the adjoint of the original problem instead. The dimension of the corresponding optimisation problem can therefore be converted from size $N$ to size $m$, where $m\ll N$ denotes the number of outputs. The exact solution of the Riccati equation associated with the original high-dimensional control problem can then be approximated without directly solving the Riccati equation itself. The method has been applied to optimal full-dimensional control design for a boundary layer \cite{semeraro2013riccati}, in which the resulting full-dimensional controllers served as excellent benchmarks for the performance evaluation of the ROM-based controllers. 

As an alternative, solving Riccati equations to obtain linear quadratic regulators (LQR) is a more direct approach to the design of full-dimensional controllers. A classic method for large-scale Riccati equations is the Chandrasekhar algorithm employed by Kailath \cite{kailath1973some} where the original problem is replaced by a series of partial differential equations. The method is feasible provided the number of inputs or outputs is much smaller than the dimension of the original system. Efforts have been devoted to improving the efficiency of the algorithm by incorporating iterative methods such as the Newton-Kleinman method \cite{kleinman1968iterative,banks1991numerical} and the Krylov subspace projection method \cite{benner2004solving}. Recently, an efficient solution method formulated by B{\"a}nsch et al. \cite{bansch2015riccati} extended the original Newton-ADI method for problems governed by ordinary differential equations (ODEs) to those governed by Differential Algebraic Equations (DAEs). The work of Benner et al. \cite{benner2016inexact} further proposed an improved inexact low-rank Newton-ADI method which incorporated many improvements \cite{benner1998exact,feitzinger2009inexact,benner2013efficient,benner2013reformulated}. An extension of the method was then presented by Benner et al. \cite{benner2020efficient} specifically for the Hessenberg index-2 DAE case arising for fluid flow problems, which shows a speed-up of approximately 100 times when compared to the algorithm used in the previous work of B{\"a}nsch et al. \cite{bansch2015riccati}. These efficient algorithms for solving large-scale, sparse Riccati equations have been incorporated into the M.E.S.S. library \cite{SaaKB19-mmess-2.0} which allows one to design full-dimensional optimal estimators and controllers for high-dimensional control problems arising from fluid flows.

\subsection{Objectives of the present work}
As mentioned above, many previous studies concerning dimension reduction and control design suffer from two major limitations: i) optimal control based on ROMs usually considers only the open-loop behaviour of the flow; ii) direct design of full-dimensional controllers is limited to cases with a small number of inputs and outputs. To address these issues, we use a resolvent-based algorithm for full-dimensional control design, which achieves globally optimal performance for various control setups in high-dimensional flow systems. For each feedback control setup, the approach finds optimal solutions to two independent problems. The first problem is the optimal estimation problem in which the objective is to estimate the flow from limited sensor measurements. The second problem is the full-state information control problem in which the entire flow field is known but only limited actuators are available for control. A major contribution of the approach is its handling of the optimal estimation and control problems i) when disturbances are applied everywhere in the domain (so-called full-state inputs) and ii) when the flow is measured everywhere (so-called full-state output). The method exploits an iterative scheme and the convergence to the global optima can be shown from the convergence of the closed-loop energy spectra.

In this study we employ the $H_2$ optimal control tools established by Doyle et al. \cite{doyle1988state} to solve the optimal feedback control problem for the two-dimensional flow past a cylinder. A complete introduction to the method can be found in Skogestad and Postlethwaite \cite{skogestad2007multivariable} and a direct application to a one-dimensional flow system can be found in Chen and Rowley \cite{chen2011h}. The article is organised as follows. In \S\ref{sec:math_formu} we formulate and discretize the time-invariant model governing linear flow perturbations. The optimal estimation and control problems arising from the flow system are presented in \S\ref{sec:intro_lqg}. A framework for solving the corresponding high-dimensional Riccati equations is then introduced in \S\ref{sec:method_intro}. Here, we present an iterative algorithm to obtain the global optimal estimator and controller, which together lead to an optimal feedback controller. In \S\ref{sec:results}, we test our method for the feedback stabilisation of the two-dimensional cylinder flow in which the feedback setup aims to attenuate linear perturbations around the base flow. Conclusions are drawn in \S\ref{sec:conclusion}

\section{Mathematical formulation}\label{sec:math_formu}
\subsection{Governing equations}
The feedback control of small-amplitude perturbations is analysed for the two-dimensional (2D) cylinder flow governed by the incompressible Navier-Stokes equations:
\begin{equation}\label{equ:nsequation}
\begin{gathered}
    \partial_t\textbf{\textit{u}}=-\textbf{\textit{u}}\cdot\nabla\textbf{\textit{u}}-\nabla \textit{p}+\Rey^{-1} \nabla^2\textbf{\textit{u}}\ ,\\
    \nabla\cdot\textbf{\textit{u}}=0\ ,
\end{gathered}
\end{equation}
where $\textbf{\textit{u}}(\textbf{\textit{x}},t)$ and $\textit{p}(\textbf{\textit{x}},t)$ are the velocity and pressure field respectively. The Reynolds number $\Rey={\textit{U}_{\infty}D}/{\nu}$ is defined using the free-stream velocity $\textit{U}_{\infty}$, the cylinder diameter $D$ and the kinematic viscosity $\nu$ to make all variables dimensionless. The nonlinear equations are then linearised about the laminar base flow $(\textbf{\textit{U}},\textit{P})$, providing linear governing equations for small-amplitude perturbations:
\begin{equation}\label{equ:pertur}
\begin{gathered}
    \partial_t\textbf{\textit{u}}'=- \textbf{\textit{U}}\cdot\nabla\textbf{\textit{u}}'-\textbf{\textit{u}}'\cdot\nabla\textbf{\textit{U}}-\nabla\textit{p}'+\Rey^{-1}\nabla^2\textbf{\textit{u}}'+\textbf{\textit{f}}'(\textbf{\textit{x}},t)\ ,\\ \nabla\cdot\textbf{\textit{u}}'=0\ .
\end{gathered}
\end{equation}
The nonlinear term $\textbf{\textit{u}}'\cdot\nabla\textbf{\textit{u}}'$ is neglected due to our assumption of small perturbations and $\textbf{\textit{f}}'(\textbf{\textit{x}},t)$ models any external forcing, such as stochastic disturbances or actuation. The governing equations \eqref{equ:pertur} can also be written compactly as
\begin{equation}\label{equ:compact_form}
    \dfrac{\partial}{\partial t}
    \underbrace{
    \begin{bmatrix}
        \text{I}&\text{0}\\
        \text{0}&\text{0}
    \end{bmatrix}}_{\mathcal{E}}
    \begin{bmatrix}
        \textbf{\textit{u}}'\\
        \textit{p}'
    \end{bmatrix}
    =
    \underbrace{
    \begin{bmatrix}
    \mathcal{L}&-\nabla()\\
    -\nabla\cdot()&\text{0}
    \end{bmatrix}}_{\mathcal{A}}
    \begin{bmatrix}
        \textbf{\textit{u}}'\\
        \textit{p}'
    \end{bmatrix}+
    \underbrace{
    \begin{bmatrix}
        \text{I}\\
        \text{0}
    \end{bmatrix}}_{\mathcal{P}}\textbf{\textit{f}}'(\textbf{\textit{x}},t)\ ,
\end{equation}
where the linear operator $\mathcal{L} = {-\textbf{\textit{U}}\cdot\nabla()} - {()\cdot\nabla\textbf{\textit{U}}} + {\Rey^{-1}\nabla^2()}$ and $\mathcal{A}$ denotes the linearized Navier–Stokes operator around the base flow. The prolongation operator $\mathcal{P}$ maps a velocity vector $\textbf{\textit{u}}'$ to a velocity-zero-pressure vector $[\textbf{\textit{u}}',\ \text{0}]^T$. 

\subsection{The discretised time-invariant model}
\begin{figure}
  \centerline{\includegraphics[width=0.75\textwidth]{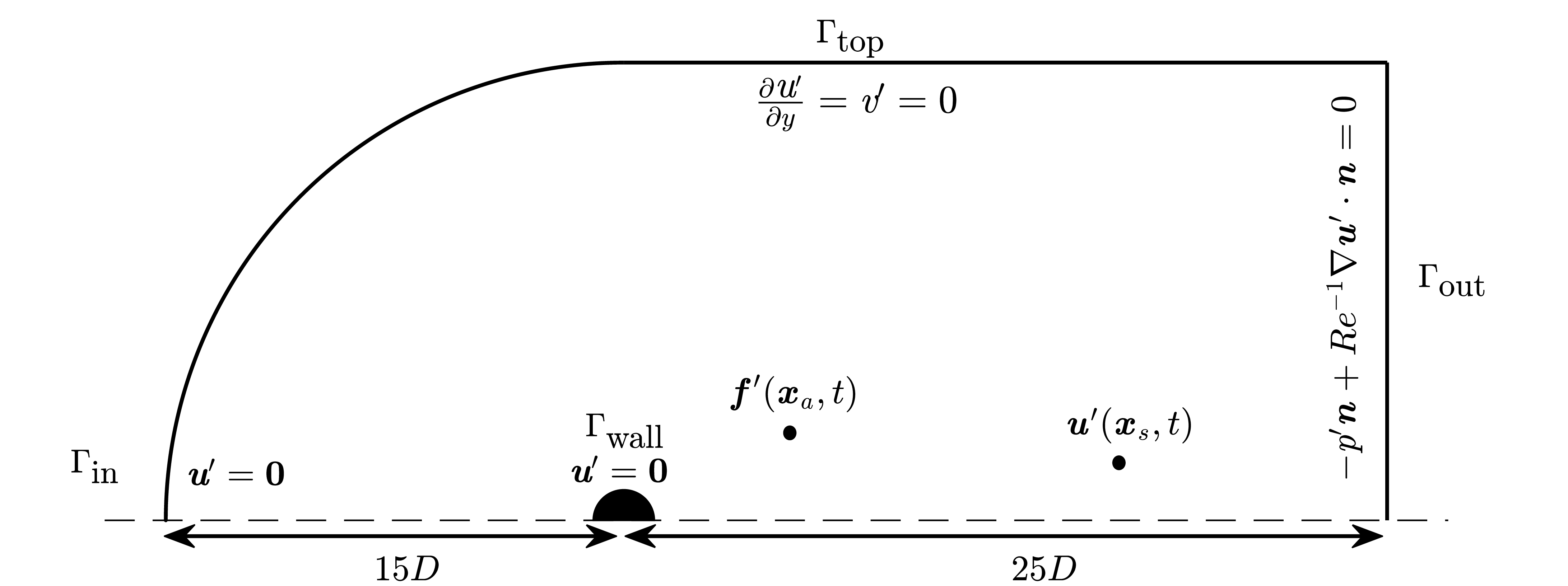}}
  \caption{Computational domain and boundary conditions for the cylinder flow (only the top half segment of the entire domain is shown).}
\label{fig:compudomain}
\end{figure}

The linearised Navier-Stokes equations \eqref{equ:compact_form} are discretized over the 2D computational domain shown in figure \ref{fig:compudomain} using Taylor-Hood mixed finite elements on the computing platform FEniCS \cite{logg2012automated}. We use a structured mesh with mesh points clustered smoothly near the cylinder and in the wake to appropriately resolve the details of the flow. The mesh consists of $2.7\times 10^4$ triangles with a minimum wall-normal spacing of $0.01$ at the cylinder surface. The compound state vector $\textbf{\textit{w}}=[\textbf{\textit{u}}',\ \textit{p}']^T\in \mathbb{R}^{N}$ thus corresponds to approximately $N\approx 1.2\times10^5$ degrees of freedom. 

Having discretised the continuous equations \eqref{equ:compact_form}, we can now express them in the form of a linear time-invariant state-space model including any inputs and outputs:
\begin{gather}\label{equ:state_space_model}
    \begin{aligned}
        \textbf{E}\Dot{\textbf{\textit{w}}}&=\textbf{A}\textbf{\textit{w}}+\textbf{B}\textbf{\textit{q}}\\
        \textbf{\textit{y}}&=\textbf{C}\textbf{\textit{w}}\ .
    \end{aligned}
\end{gather}
The system matrix $\textbf{A}$ denotes the discretized linear Navier-Stokes operator and the global mass matrix $\textbf{E}=\textbf{P}\textbf{M}\textbf{P}^T$ is defined using the prolongation matrix $\textbf{P}$ and the finite element mass matrix $\textbf{M}$ associated with the velocity field. Note that the superscript $(\cdot)^T$ designates the transpose and we use $(\cdot)^H$ for the conjugate transpose. The spatial-temporal discretization of any external forcing $\textbf{\textit{f}}'$ (e.g.~an actuator at the position $\textbf{\textit{x}}_a$) is characterised by the matrix $\textbf{B}$ and the time signal $\textbf{\textit{q}}$. The measurement $\textbf{\textit{y}}$ is an output of interest (e.g.~velocity $\textbf{\textit{u}}'$ at the position $\textbf{\textit{x}}_s$) and is characterised by the sensing matrix $\textbf{C}$. A feedback connection between the sensor output and the actuator input is then established to attenuate flow perturbations.

The boundary conditions used for the linear model \eqref{equ:state_space_model} are summarised in figure \ref{fig:compudomain}, and have been used in many previous studies \cite{leontini2006wake,jin2019feedback}. The base flow has the same boundary conditions as those depicted in the figure except at the inlet $\rmGamma_{\text{in}}$ where a uniform freestream velocity (i.e.~$\textbf{\textit{U}}=[\textit{U}_{\infty},0]$) is enforced. Note that the laminar base flows and discretised perturbation systems have been validated by comparing them with the stability analysis results \cite{Barkley06}. A sparse direct LU solver (MUMPS \cite{amestoy2001fully}) and iterative Arnoldi methods (ARPACK \cite{lehoucq1998arpack}) are used for all linear problems encountered in the study.

\section{Optimal estimator and controller design}\label{sec:intro_lqg}
\newcolumntype{a}{>{\centering\arraybackslash}m{0.47\textwidth}}
\newcolumntype{b}{>{\centering\arraybackslash}m{0.94\textwidth}}
\aboverulesep=0ex
\belowrulesep=0ex
\begin{figure*}
\begin{tabular}{@{}|a|a|@{}}
\toprule
\multicolumn{2}{|b|}{
    \hspace{0mm}
	\centerline{
	\includegraphics[width=0.94\textwidth]{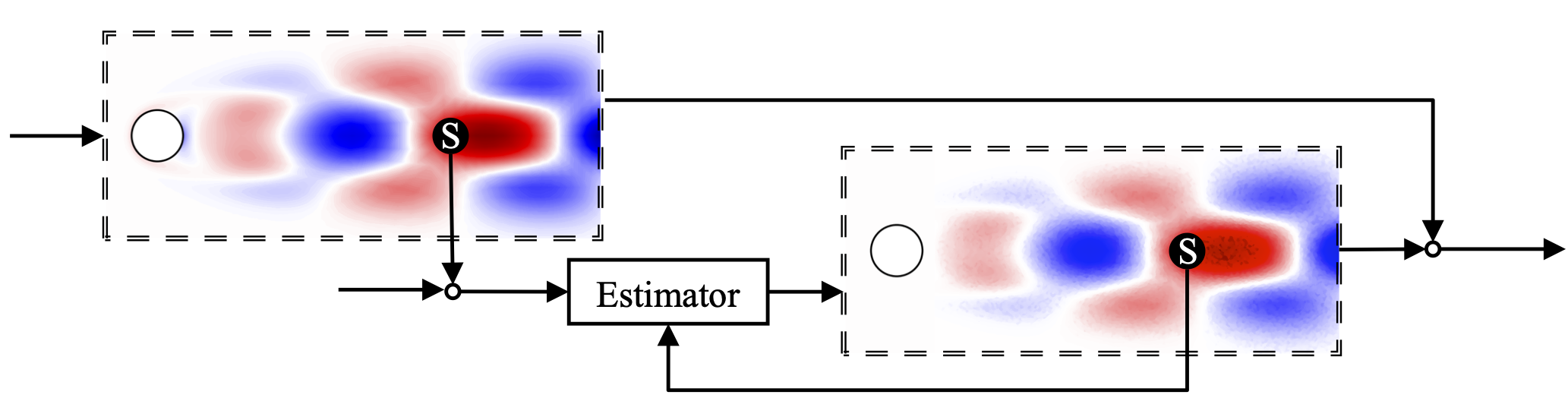} 
	\llap{\parbox[b]{4.5in}{(a) OE \\\rule{0ex}{1.2in}}}
    \llap{\parbox[b]{3.1in}{\textit{\textbf{y}}\\\rule{0ex}{0.24in}}}
    \llap{\parbox[b]{2.775in}{$\textit{\textbf{y}}_e$\\\rule{0ex}{0.0in}}}
    \llap{\parbox[b]{2.105in}{$\textit{\textbf{w}}_e$\\\rule{0ex}{0.6in}}}
    \llap{\parbox[b]{4.225in}{$\textit{\textbf{w}}$\\\rule{0ex}{0.95in}}}
    \llap{\parbox[b]{0.5in}{$\textit{\textbf{e}}$\\\rule{0ex}{0.485in}}}
    \llap{\parbox[b]{4.55in}{$\textit{\textbf{d}}$\\\rule{0ex}{0.8in}}}
    \llap{\parbox[b]{3.65in}{$\textit{\textbf{n}}$\\\rule{0ex}{0.265in}}}
	}
} \\[0pt] \midrule
velocity sensor at $\textbf{\textit{x}}_s$ & no actuator \\ \midrule
\multicolumn{2}{b}{}\\\midrule
\multicolumn{2}{|b|}{
    \vspace{2mm}
	\centerline{
	\hspace{7.5mm}
	\includegraphics[width=0.94\textwidth]{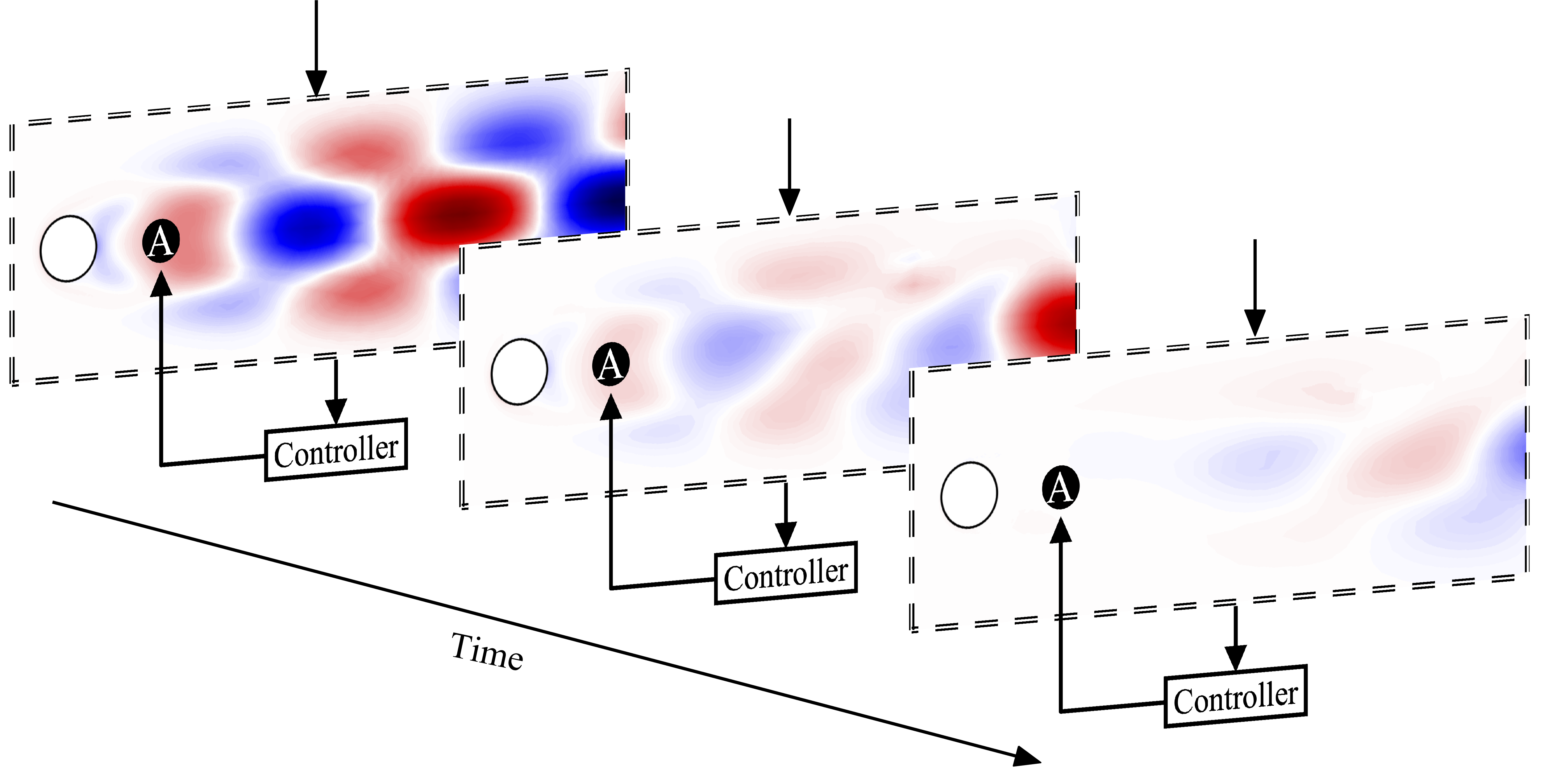} 
	\llap{\parbox[b]{4.625in}{(b) FIC\\\rule{0ex}{2.1in}}}
    \llap{\parbox[b]{4.375in}{$\textit{\textbf{w}}(t_0)$\\\rule{0ex}{1.7in}}}
    \llap{\parbox[b]{4.15in}{\textit{\textbf{q}}\\\rule{0ex}{1.0in}}}
    \llap{\parbox[b]{3.775in}{\textit{\textbf{d}}\\\rule{0ex}{2.05in}}}
    \llap{\parbox[b]{3.225in}{$\textit{\textbf{w}}(t_1)$\\\rule{0ex}{1.34in}}}
    \llap{\parbox[b]{3.00in}{\textit{\textbf{q}}\\\rule{0ex}{0.65in}}}
    \llap{\parbox[b]{2.565in}{\textit{\textbf{d}}\\\rule{0ex}{1.7in}}}
    \llap{\parbox[b]{2.05in}{$\textit{\textbf{w}}(t_2)$\\\rule{0ex}{1.0in}}}
    \llap{\parbox[b]{1.85in}{\textit{\textbf{q}}\\\rule{0ex}{0.3in}}}
    \llap{\parbox[b]{1.36in}{\textit{\textbf{d}}\\\rule{0ex}{1.36in}}}
	}
} \\[-10pt] \midrule
no sensor & body force actuator at $\textbf{\textit{x}}_a$ \\ \midrule
\multicolumn{2}{b}{}\\\midrule
\multicolumn{2}{|b|}{
    \vspace{2mm}
	\centerline{
	\hspace{12.5mm}
    \includegraphics[width=0.94\textwidth]{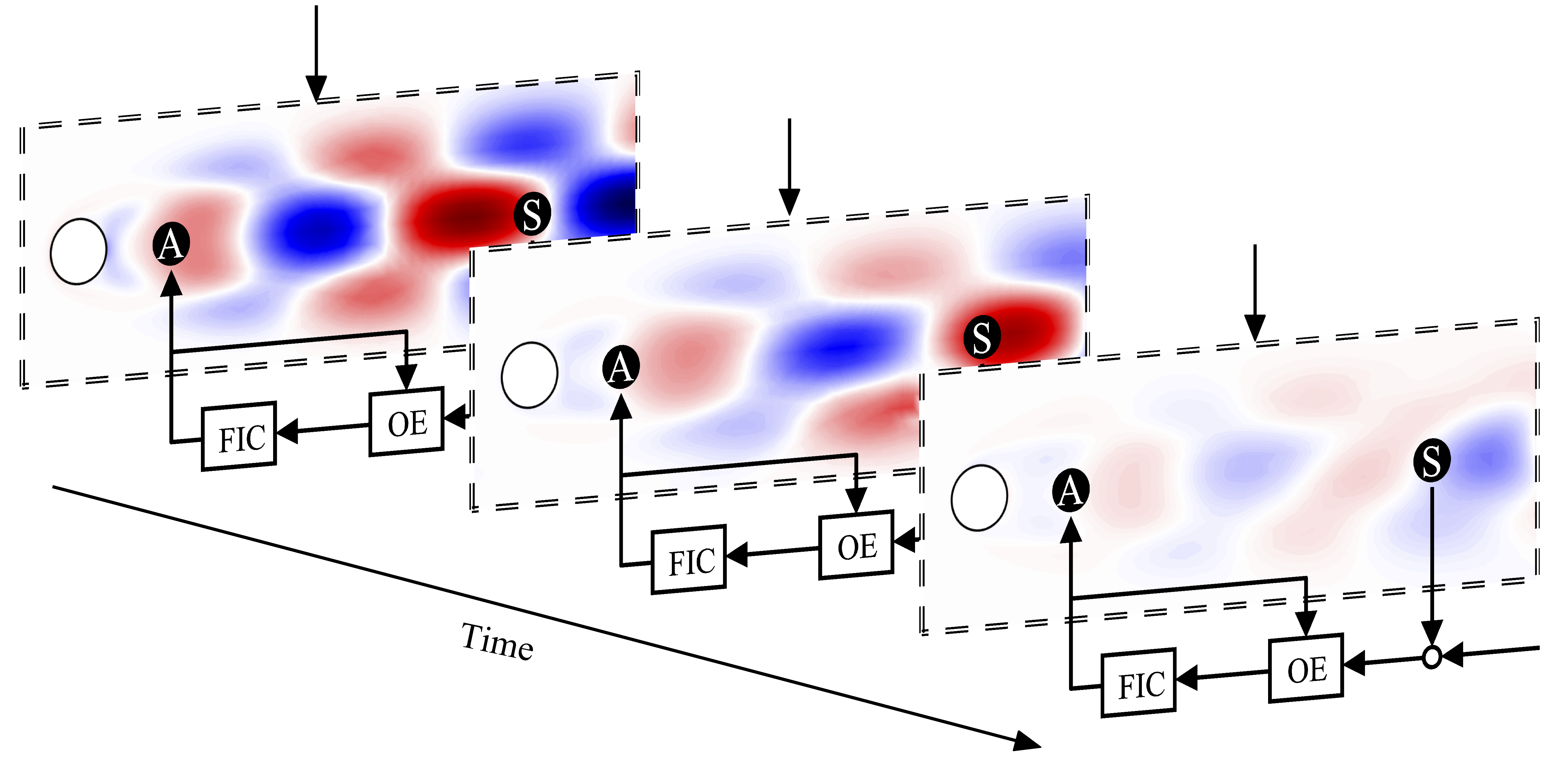}
    \llap{\parbox[b]{4.625in}{(c) IOC\\\rule{0ex}{2.1in}}}
    \llap{\parbox[b]{0.275in}{\textit{\textbf{n}}\\\rule{0ex}{0.275in}}}
    \llap{\parbox[b]{0.6in}{\textit{\textbf{y}}\\\rule{0ex}{0.23in}}}
    \llap{\parbox[b]{1.64in}{\textit{\textbf{q}}\\\rule{0ex}{0.32in}}}
    \llap{\parbox[b]{1.175in}{$\textit{\textbf{w}}_e$\\\rule{0ex}{0.2in}}}
    \llap{\parbox[b]{1.2in}{\textit{\textbf{d}}\\\rule{0ex}{1.36in}}}
    
    \llap{\parbox[b]{3.0in}{\textit{\textbf{q}}\\\rule{0ex}{0.65in}}}
    \llap{\parbox[b]{2.575in}{$\textit{\textbf{w}}_e$\\\rule{0ex}{0.55in}}}
    \llap{\parbox[b]{2.625in}{\textit{\textbf{d}}\\\rule{0ex}{1.7in}}}
    
    \llap{\parbox[b]{4.4in}{\textit{\textbf{q}}\\\rule{0ex}{1.0in}}}
    \llap{\parbox[b]{3.95in}{$\textit{\textbf{w}}_e$\\\rule{0ex}{0.9in}}}
    \llap{\parbox[b]{4.075in}{\textit{\textbf{d}}\\\rule{0ex}{2.025in}}}
    
    \llap{\parbox[b]{4.775in}{$\textit{\textbf{w}}(t_0)$\\\rule{0ex}{1.7in}}}
    \llap{\parbox[b]{3.55in}{$\textit{\textbf{w}}(t_1)$\\\rule{0ex}{1.34in}}}
    \llap{\parbox[b]{2.3in}{$\textit{\textbf{w}}(t_2)$\\\rule{0ex}{1.0in}}}
    }
} \\[-15pt] \midrule
velocity sensor at $\textbf{\textit{x}}_s$ & body force actuator at $\textbf{\textit{x}}_a$ \\ \bottomrule
\end{tabular}
\caption[Setups of the OE, FIC and IOC problems]{(a) Optimal estimation (OE) of the whole flow field using only limited sensor measurements. (b) Full-state information control (FIC) using only limited actuators for control when the entire flow field is known. (c) Input-Output control (IOC) with only limited sensor measurements and limited actuation.} \label{fig:controlsetup}
\end{figure*}

The aim of this section is to introduce the optimal estimation and control problems concisely, as these will be used in \S\ref{sec:method_intro}. We describe the optimal estimation and control problems in \S\ref{sec:intro_lqg_1} and their common governing equations are formulated in \S\ref{sec:intro_matrices}. The numerical scheme to solve the corresponding optimisation problems is introduced in \S\ref{sec:intro_solver}.

\subsection{The estimation and control problems}\label{sec:intro_lqg_1}
Our objective is to find a control signal from sensor measurements such that the flow perturbation is minimised under the excitation of external disturbances and in the presence of sensor noise. The control that will be optimal in minimisation flow perturbations arises, in one part, from the optimal filtering of any noise that corrupts our signal and, in another part, from the optimal control when the entire flow state is assumed to be known. This corresponds to two separate problems: i) the optimal estimation (OE) problem and ii) the full-information control (FIC) problem, which can be solved independently based on the Separation Theorem \cite{georgiou2013separation}. The optimal feedback control (i.e.~input-output control) law can then be constructed once these two problems are solved. Figure \ref{fig:controlsetup} shows a summary of the problems investigated.


\subsubsection{The optimal estimation (OE) problem}\label{sec:intro_oe}
In order to employ feedback, one needs to know the flow state. However, only the sensor measurement is available in any practical control setup. In the optimal estimation (OE) problem we therefore aim to estimate or observe the whole flow field from limited sensor measurements. More specifically, the flow estimate $\textbf{\textit{w}}_e$ is generated by applying the estimator gain $\textbf{K}_f$ to the difference between the sensor measurements from the actual flow $\textbf{\textit{y}}(t)$ and the estimate $\textbf{\textit{y}}_e(t)$, as shown in figure \ref{fig:controlsetup}(a). A major challenge for the OE problem is to form the estimator such that the differences between the actual flow field and the estimate are small in the presence of stochastic disturbances $\textbf{\textit{d}}$ and sensor noise $\textbf{\textit{n}}$. This issue is addressed by solving the optimisation problem that minimises the mean (time-averaged) kinetic energy of the estimation error $\textbf{\textit{J}}_{\text{OE}}$:
\begin{equation}\label{equ:lqe_joe}
    \textbf{\textit{J}}_{\text{OE}}=\lim_{t\rightarrow\infty}\dfrac{1}{T}\int_0^T\textbf{\textit{e}}^T\textbf{Q}\textbf{\textit{e}}\ dt\ ,
\end{equation}
where $\textbf{\textit{e}}(t)=\textbf{\textit{w}}(t)-\textbf{\textit{w}}_e(t)$ denotes the estimation error and the matrix $\textbf{Q}=\textbf{P}\textbf{M}\textbf{P}^T$ represents spatial integration for computing the total kinetic energy. 

\subsubsection{The full-state information control (FIC)  problem}\label{sec:intro_fic}
The second step in the design of the optimal feedback control involves the solution of an optimal state-feedback control problem in which the entire flow field is assumed to be perfectly measured or estimated (i.e.~$\textit{\textbf{y}}=\textit{\textbf{w}}$). The task of the FIC problem is therefore to attenuate flow perturbations using only limited actuators in the presence of stochastic disturbances $\textbf{\textit{d}}$, as shown in figure \ref{fig:controlsetup}(b). We aim to find the optimal full-state controller gain $\textbf{K}_r$ that minimises flow perturbations while also maintaining a sensible control effort:
\begin{align}\label{equ:lqr_jfic}
    \textbf{\textit{J}}_{\text{FIC}}=\lim_{t\rightarrow\infty}\dfrac{1}{T}\int_0^T\left(\textbf{\textit{w}}^T\textbf{Q}\textbf{\textit{w}}+\textbf{\textit{q}}^T\textbf{R}\textbf{\textit{q}}\right) dt\ .
\end{align}
The cost function $\textbf{\textit{J}}_{\text{FIC}}$ is composed of two contributions: i) the kinetic energy of perturbations and ii) the control cost scaled by a matrix $\textbf{R}$. The matrix $\textbf{R}$ is a weight matrix that controls the balance between minimising flow perturbations and minimising the control cost. 

\subsubsection{The input-output control (IOC) problem}\label{sec:intro_ioc}
Once the above two problems are solved, the optimal control signal $\textbf{\textit{q}}$ can be computed by applying the optimal full-state controller to the estimated state, as shown by figure \ref{fig:controlsetup}(c). As suggested by the separation principle, the feedback control law constructed from an optimal estimator and an optimal full-state controller is itself optimal. In this case, the optimal estimator provides the best prediction of the entire flow field from limited sensors contaminated by measurement noise. Based on this prediction, the optimal control signal is then generated for limited actuators and gives the best attenuation of flow perturbations while maintaining a sensible control cost, i.e.~minimising the cost function $\textbf{\textit{J}}_{\text{IOC}}$ defined in the same form as equation \eqref{equ:lqr_jfic} for a specific choice of $\textbf{R}$. 

\subsection{The systems for estimation and control}\label{sec:intro_matrices}

\begin{table}
\small\addtolength{\tabcolsep}{-0.175pt}
\begin{tabular}{|c|c|ccc|cc|cc|}
\hline
\multirow{2}{*}{P} &
  State &
  \multicolumn{3}{c|}{Input} &
  \multicolumn{2}{c|}{Output $\textbf{\textit{y}}$} &
  \multicolumn{2}{c|}{Performance $\textbf{\textit{z}}$} \\ \cline{2-9} 
& \textbf{\textit{x}} & \textbf{\textit{q}} & $\textbf{B}_q$  & $\textbf{B}_d$ & $\textbf{C}_y$& \multicolumn{1}{c|}{$\textbf{D}_n$}& $\textbf{C}_z$& $\textbf{D}_q$\\ \hline
\multirow{4}{*}{OE} & 
\multirow{4}{*}{\textbf{\textit{e}}} & 
\multirow{2}{*}{-\textbf{\textit{y}}} & 
\multirow{2}{*}{$\textbf{K}_f$} & 
\multirow{2}{*}{$\textbf{W}^{1/2}$} & 
\multirow{2}{*}{$\textbf{C}$} & 
\multicolumn{1}{c|}{\multirow{2}{*}{$\textbf{V}^{1/2}$}} & 
\multirow{2}{*}{$\textbf{Q}^{1/2}$} & 
\multirow{2}{*}{$\textbf{0}$} \\
& & & & & & \multicolumn{1}{c|}{} & & \\ \cline{3-9} & &
\multicolumn{5}{c}{\multirow{2}{*}{$\textbf{E}\textbf{X}\textbf{A}^T+\textbf{A}\textbf{X}\textbf{E}-\textbf{E}\textbf{X}\textbf{C}^T\textbf{V}^{-1}\textbf{C}\textbf{X}\textbf{E}+\textbf{W}=0$}} & 
\multicolumn{2}{c|}{\multirow{2}{*}{$\textbf{K}_f=\textbf{E}\textbf{X}\textbf{C}^T\textbf{V}^{-1}$}} \\ & &
\multicolumn{5}{c}{} & 
\multicolumn{2}{c|}{} \\ \hline
\multirow{5}{*}{FIC} & 
\multirow{5}{*}{\textbf{\textit{w}}} & 
\multirow{3}{*}{-\textbf{\textit{y}}} & 
\multirow{3}{*}{$\textbf{B}$}& \multirow{3}{*}{$\textbf{W}^{1/2}$}   & 
\multirow{3}{*}{$\textbf{K}_r$} & \multicolumn{1}{c|}{\multirow{3}{*}{$\textbf{0}$}} & 
\multirow{3}{*}{\begin{tabular}[c]{@{}c@{}}
        $\begin{bmatrix}
        \textbf{Q}^{1/2}\\ \textbf{0}
        \end{bmatrix}$
    \end{tabular}} & 
\multirow{3}{*}{\begin{tabular}[c]{@{}c@{}}
        $\begin{bmatrix}
        \textbf{0} \\ \textbf{R}^{1/2}
        \end{bmatrix}$
    \end{tabular}} \\   
& & & & & & \multicolumn{1}{c|}{} & & \\
& & & & & & \multicolumn{1}{c|}{} & & \\ \cline{3-9} & &
\multicolumn{5}{c}{\multirow{2}{*}{$\textbf{E}^T\textbf{Y}\textbf{A}+\textbf{A}^T\textbf{Y}\textbf{E}-\textbf{E}^T\textbf{Y}\textbf{B}\textbf{R}^{-1}\textbf{B}^T\textbf{Y}\textbf{E}+\textbf{Q}=0$}} & 
\multicolumn{2}{c|}{\multirow{2}{*}{$\textbf{K}_r=\textbf{R}^{-1}\textbf{B}^T\textbf{Y}\textbf{E}$}} \\ & &
\multicolumn{5}{c}{} & 
\multicolumn{2}{c|}{} \\ \hline
\multirow{3}{*}{IOC} &
  \multirow{3}{*}{\begin{tabular}[c]{@{}c@{}}
    $\begin{bmatrix}
    \textbf{\textit{w}} \\ \textbf{\textit{e}}
    \end{bmatrix}$
  \end{tabular}} &
  \multirow{3}{*}{-\textbf{\textit{y}}} &
  \multirow{3}{*}{\begin{tabular}[c]{@{}c@{}}
    $\begin{bmatrix}
    \textbf{B}&\textbf{0}\\ \textbf{0}&\textbf{K}_f
    \end{bmatrix}$
  \end{tabular}} &
  \multirow{3}{*}{\begin{tabular}[c]{@{}c@{}}
    $\begin{bmatrix}
    \textbf{W}^{1/2} \\ \textbf{W}^{1/2}
    \end{bmatrix}$
  \end{tabular}} &
  \multirow{3}{*}{\begin{tabular}[c]{@{}c@{}}
    $\begin{bmatrix}
    \textbf{K}_r&-\textbf{K}_r\\ \textbf{0}&\textbf{C}
    \end{bmatrix}$
  \end{tabular}} &
  \multirow{3}{*}{\begin{tabular}[c]{@{}c@{}}
    $\begin{bmatrix}
    \textbf{0}\\ \textbf{V}^{1/2}
    \end{bmatrix}$
  \end{tabular}} &
  \multirow{3}{*}{\begin{tabular}[c]{@{}c@{}}
    $\begin{bmatrix}
    \textbf{Q}^{1/2}&\textbf{0}\\ \textbf{0}&\textbf{0}
    \end{bmatrix}$
  \end{tabular}} &
  \multirow{3}{*}{\begin{tabular}[c]{@{}c@{}}
    $\begin{bmatrix}
    \textbf{0}&\textbf{0}\\ \textbf{0}&\textbf{R}^{1/2}
    \end{bmatrix}$
  \end{tabular}} \\
 & & & & & & & & \\
 & & & & & & & & \\ \hline
\end{tabular}
\caption[A summary of system states, matrices and  Riccati equations for the OE, FIC and IOC problems.]{A summary of system states, input and output matrices, and Riccati equations for the OE, FIC and IOC problems.}
\label{tab:sysmats_summary}
\end{table}

\begin{figure}
  \centerline{\includegraphics[width=0.6\textwidth]{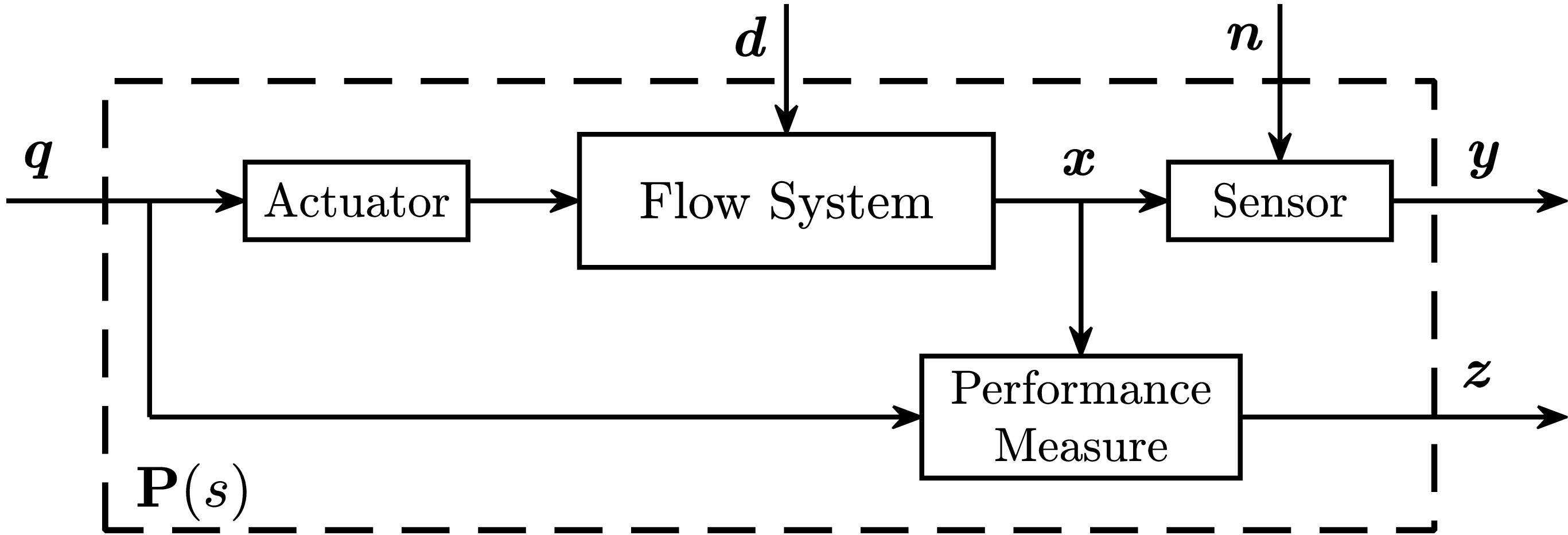}}
  \caption[Block diagram of the closed-loop state-space model.]{Block diagram of the state-space model \eqref{equ:ssmodel_summary}, which is denoted as a transfer function $\textbf{P}(s)$ with the Laplace variable $s$.}
\label{fig:blockdiagram_closed}
\end{figure}

We now consider the governing equations for the OE, FIC and IOC problems, which can be expressed as a linear time-invariant state-space model $\textbf{P}(s)$ shown by the block diagram in figure \ref{fig:blockdiagram_closed}:
\begin{gather}\label{equ:ssmodel_summary}
    \begin{aligned}
        \textbf{E}\Dot{\textbf{\textit{x}}}&=\textbf{A}\textbf{\textit{x}}+\textbf{B}_q\textbf{\textit{q}}+\textbf{\textit{B}}_d\textbf{\textit{d}}\\
        \textbf{\textit{y}}&=\textbf{C}_y\textbf{\textit{x}}+\textbf{D}_n\textbf{\textit{n}}\\
        \textbf{\textit{z}}&=\textbf{C}_z\textbf{\textit{x}}+\textbf{D}_q\textbf{\textit{q}}\ ,
    \end{aligned}
\end{gather}
where the system state $\textbf{\textit{x}}$ represents the internal state of interest for each of the three problems. Table \ref{tab:sysmats_summary} summarises the system states, inputs, outputs and the performance-measuring matrices for the OE, FIC and IOC problems. In particular, the optimal estimator gain $\textbf{K}_f$ and the optimal full-state controller gain $\textbf{K}_r$ are formed from the unique solutions of the generalised algebraic Riccati equations associated with the OE and FIC problems, respectively.

Both the disturbances $\textbf{\textit{d}}$ and sensor noise $\textbf{\textit{n}}$ are modelled as uncorrelated zero-mean Gaussian white noise. We inject stochastic disturbances over the entire velocity field (i.e.~as the external forcing), and $\textbf{\textit{d}}$ is an $N_u$-by-1 vector whose statistical properties are given by $\textbf{W}^{1/2}=\textbf{P}(\textbf{M}^{1/2})^T$ after the spatial discretization using a finite-element method \cite{croci2018efficient}. Here, $N_u$ denotes the size of the velocity state and the matrix square root $\textbf{M}^{1/2}$ computed from the Cholesky decomposition satisfies $(\textbf{M}^{1/2})^T\textbf{M}^{1/2}=\textbf{M}$. We also include a contribution from sensor noise in the measurement $\textbf{\textit{y}}$ and the diagonal matrix is set to $\textbf{V}^{1/2}=\alpha\textbf{I}$. 

The performance of either the estimator or the controller can be expressed by the mean kinetic energy of the output $\textbf{\textit{z}}$:
\begin{equation}\label{equ:cost_j}
    \textbf{\textit{J}}=\lim_{t\rightarrow\infty}\dfrac{1}{T}\int_0^T\textbf{\textit{z}}^T\textbf{\textit{z}}\ dt\ .
\end{equation}
Note that we choose $\textbf{Q}^{1/2}=\textbf{M}^{1/2}\textbf{P}^T$ to measure the kinetic energy of the whole perturbation field and $\textbf{R}^{1/2}=\alpha\textbf{I}$ to equally penalise all control signals.

\subsection{Numerical Riccati solver}\label{sec:intro_solver}
The state-space models arising from the discretization of two or three dimensional fluid flows are generally of high dimension. Although it is common to perform numerical simulations for high-dimensional systems (i.e.~$N>10^5$), traditional control design tools (e.g. Riccati solvers) typically become computationally intractable for $N>10^3$. The challenge of solving the large sparse Riccati equations associated with the OE and FIC problems has been partially overcome by a sparse Riccati solver in the M.E.S.S. library \cite{benner2019mess,SaaKB19-mmess-2.0} using an extended low-rank method. However, a major limitation of the method is that the number of inputs and outputs (so-called terminals) $m$ should be much smaller than the dimension of the control problem, i.e.~$m \ll\textit{N}$. For problems with either many inputs (e.g.~disturbances applied everywhere) or many outputs (e.g.~perturbations measured everywhere), no efficient numerical tools are available to directly handle large-scale systems and some special `terminal reduction' techniques need to be considered to properly construct low-rank matrices $\textbf{W}^{1/2}$ and $\textbf{Q}^{1/2}$ for the OE and FIC problems, respectively. In the next section, we will present an iterative algorithm that combines resolvent analysis and proper orthogonal decomposition (POD) for full-dimensional optimal estimator and controller design.

\section{Resolvent-based performance optimisation}\label{sec:method_intro}
To overcome the challenges associated with many inputs and outputs, we now introduce a `terminal reduction' technique that allows the design of optimal full-dimensional estimators and controllers in the presence of disturbances everywhere (i.e~input $\textit{\textbf{d}}$) and when measurement are applied everywhere (i.e~output $\textit{\textbf{z}}$). The method is developed based on resolvent analysis for mode selection (\S\ref{sec:method_intro_resol}) and the proper orthogonal decomposition (POD) in frequency space for efficient construction of low-rank orthonormal bases (\S\ref{sec:method_intro_pod}). The covariance and weight matrices appearing in the Riccati equations are then replaced by these low-rank linear bases. In \S\ref{sec:method_intro_algorithm}, we employ the method iteratively to find globally optimal estimators and controllers.

\subsection{Linear optimal forcing and response modes}\label{sec:method_intro_resol}
\subsubsection{The $H_2$ norm}\label{sec:method_intro_resol_2norm}
The $H_2$ norm of the transfer function $\textbf{Z}(s)$ is a common way to quantify the performance of estimators or full-state information controllers, which can be considered as the summation of the squared $H_2$ norms of two subsystems:
\begin{equation}\label{equ:h2norm_resovent}
    \textbf{\textit{J}}=\norm{\textbf{Z}(s)}^2_2=\norm{\textbf{Z}_1(s)}^2_2+\norm{\textbf{Z}_2(s)}^2_2\ ,
\end{equation}
where $\textbf{\textit{J}}$ denotes the cost function defined for the OE problem (see \eqref{equ:lqe_joe}) or the FIC problem (see \eqref{equ:lqr_jfic}). The subsystem $\textbf{Z}_2(s)$ represents the transfer function that governs either the effect of sensor noise in the OE problem or the control cost in the FIC problem (see equations (\ref{equ:app.OE_measure_zs}, \ref{equ:app.FIC_measure_zs}) in Appendix \ref{sec:app.a0}). Note that the dimensionality of $\textbf{Z}_2(s)$ is determined by the number of sensors or actuators, which is generally small and thus computing its norm is feasible and well handled by the Riccati solver. However,  $\textbf{Z}_1(s)$ represents the high-dimensional transfer function from the disturbances $\textbf{\textit{d}}$ to either the estimation error $\textbf{\textit{e}}$ in the OE problem or the perturbations $\textbf{\textit{w}}$ in the FIC problem. The corresponding $H_2$ norm squared is:
\begin{gather}\label{equ:h2norm_resovent_z1}
    \begin{aligned}
        \norm{\textbf{Z}_1(s)}^2_2&=\dfrac{1}{2\pi}\int_{-\infty}^{\infty}\textrm{Trace}\{\textbf{Z}_1^H(j\omega)\textbf{Z}_1(j\omega)\}\ d\omega\\
         &=\dfrac{1}{2\pi}\int_{-\infty}^{\infty}\sum_{i=1}^{N_u}\sigma^2_i(j\omega)\ d\omega\ ,
    \end{aligned}
\end{gather}
where $\sigma_i(j\omega)$ are the singular values of the transfer function $\textbf{Z}_1(s)$ at frequency $\omega$ arranged in descending order. The singular values of a transfer function can be considered as energy gains between a series of inputs and the corresponding outputs. We thus aim to minimise the integrated energy gain \eqref{equ:h2norm_resovent_z1} for inputs and outputs over all frequencies and all possible directions. 

However, it is not feasible to consider all inputs or outputs while designing estimators or controllers for a high-dimensional flow system (e.g.~$N>10^5$). Instead, we consider an alternative cost function $\gamma^2{(k,\omega_n)}$ which only considers a limited number of inputs and outputs within a specified frequency range:
\begin{equation}\label{equ:h2norm_resovent1}
     \gamma^2(k,\omega_n)=\dfrac{1}{2\pi}\int_{-\omega_n}^{\omega_n}\sum_{i=1}^{k}\sigma^2_i(j\omega)\ d\omega+\norm{\textbf{Z}_2(s)}^2_2\ ,
\end{equation}
in which only the first $k$ orthogonal inputs and outputs across a limited frequency range $\omega\in[-\omega_n,\ \omega_n]$ are considered. This cost function is constructed based on two facts: i) for fluid flows, only a limited number of dominant physical mechanisms occur within a finite frequency range, e.g.~the instability of the linearised cylinder flow occurs around $\omega_c\approx 0.8$; ii) these physical mechanisms can be approximated by a small number of orthogonal inputs and outputs that have large energy gains $\sigma_i^2$, which are also the most significant for estimation or control. Instead of minimising all energy gains over all frequencies and all possible directions as described by the original cost function \eqref{equ:h2norm_resovent}, it is more feasible to use the cost function \eqref{equ:h2norm_resovent1} that considers a significantly smaller number of inputs and outputs. Note that by choosing a sufficient number of orthogonal inputs and outputs over a sufficient frequency range, the performance of the estimator or controller should converge to the true global optimum. 

\subsubsection{Resolvent analysis}\label{sec:method_intro_resol_analysis}
\begin{figure}
    \centerline{\includegraphics[width=0.95\textwidth]{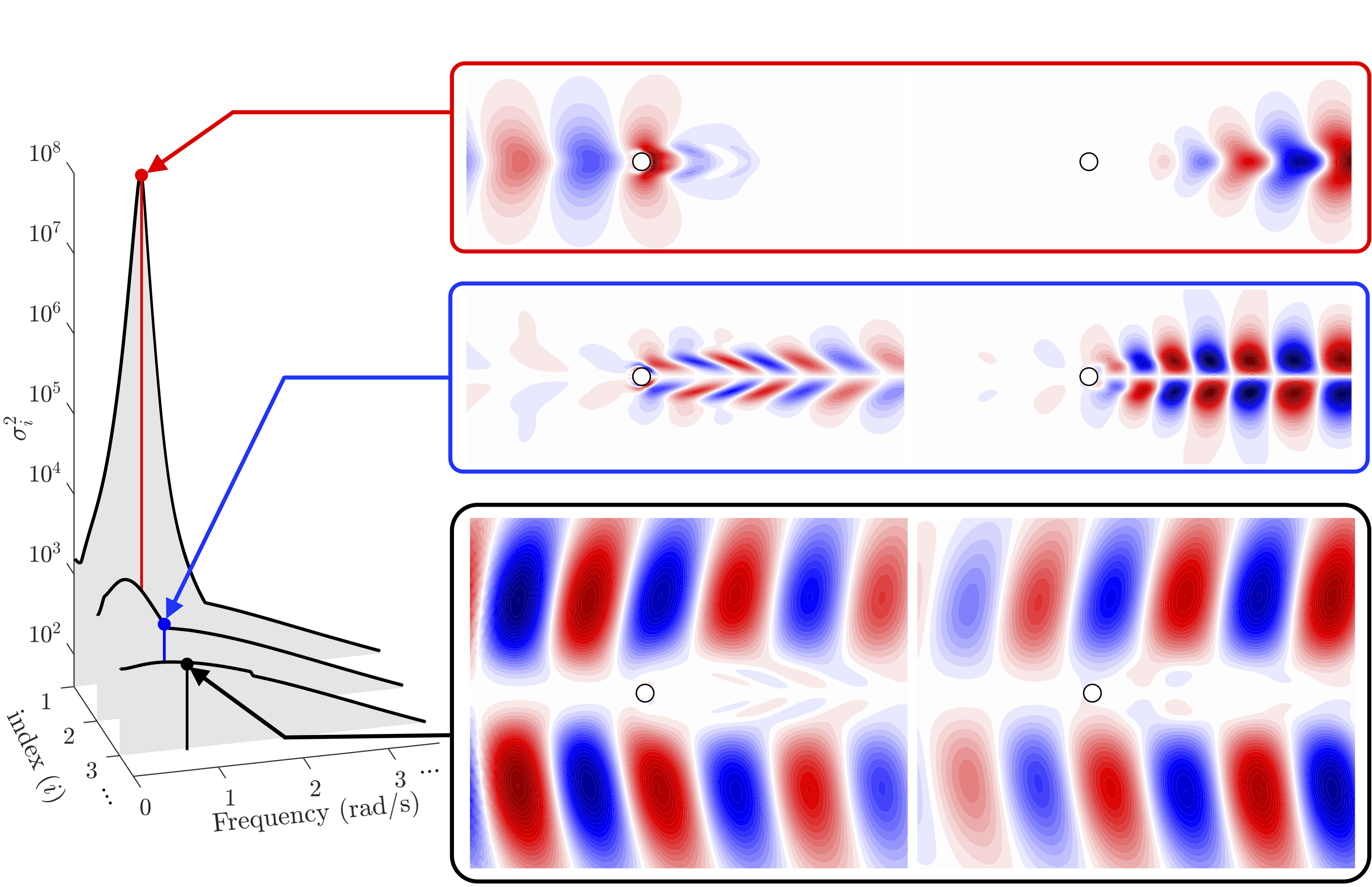}
    \llap{\parbox[b]{2.65in}{Input mode $\hat{\textbf{\textit{f}}}_i$ \\\rule{0ex}{2.7in}}}
    \llap{\parbox[b]{1.2in}{Output mode $\hat{\textbf{\textit{u}}}_i$\\\rule{0ex}{2.7in}}}
  }
  \caption[Resolvent analysis of the cylinder flow at $\Rey=90$.]{Resolvent analysis of the cylinder flow at $\Rey=90$. The first three energy spectra ($\sigma_i^2$) and the corresponding optimal input and output modes (transverse component, real part) at the unstable frequency $\omega\approx 0.76$ are shown for comparison.}
\label{fig:resolvent_diagram}
\end{figure}

We now introduce the framework of resolvent analysis which forms the optimal forcing that maximises the energy gain of the response. This is achieved using the generalised singular value decomposition of the transfer function $\textbf{Z}_1(s)$:
\begin{equation}\label{equ:sec4.1_svd}
    \textbf{Z}_1(s)=\textbf{Q}^{1/2}(s\textbf{E}-\textbf{A})^{-1}\textbf{W}^{1/2}=\left(\textbf{Q}^{1/2}\hat{\textbf{U}}\right)\rmSigma\left(\hat{\textbf{F}}\hspace{0mm}^H\textbf{P}^T\textbf{W}^{1/2}\right)\ .
\end{equation}
Here, we are considering the open-loop system (i.e.~before any estimators or controllers are designed). The Laplace variable $s$ is evaluated at a particular frequency $j\omega$ and $\rmSigma$ is a diagonal matrix containing the singular values $\sigma_i$ in descending order. The input and output matrices $\hat{\textbf{F}}$ and $\hat{\textbf{U}}$ each form an orthonormal basis with respect to the matrices $\textbf{W}$ and $\textbf{Q}$. 

Figure \ref{fig:resolvent_diagram} shows the first three singular values for the cylinder flow at $\Rey=90$, and compares the corresponding input and output modes at the unstable frequency $\omega\approx 0.76$. The first two output modes correspond to vortex shedding. In particular, the leading output mode $\hat{\textit{\textbf{u}}}_1$ corresponds to the amplification of upstream inputs while the second output mode corresponds to a local amplification of disturbances. There exists a large separation between the first ($\sigma_1^2$) and second ($\sigma_2^2$) energy gains. The third output mode, however, consists of waves that are concentrated in the free-stream. For estimation and control, we might expect that the physical mechanisms represented by the first two resolvent modes (e.g.~vortex shedding) are the most important. It is therefore acceptable to neglect free-stream modes and other less energetic modes while designing optimal estimators and controllers.

In the OE problem, the white noise disturbances $\textbf{\textit{d}}$ can be thought of as random combinations of the input basis in time, as given by $\textbf{W}^{1/2}=\textbf{P}\textbf{M}\hat{\textbf{F}}$. It is worthwhile to note that there are usually large separations between singular values $\sigma_i$ such that only a limited number of forcing modes give rise to energetic responses that are important for estimation. This allows us to consider only a truncated input matrix $\hat{\textbf{F}}_k$ (with only the first $k$ columns) for optimal estimator design. As for the FIC problem with perturbation measurements everywhere, the system states can be considered as linear combinations of the output basis as $\textbf{Q}^{1/2}=\hat{\textbf{U}}\hspace{0mm}^H\textbf{E}$. A truncated output matrix $\hat{\textbf{U}}_k$ can thus be used to reduce the number of outputs for optimal controller design since a significant proportion of the energy of the flow's response is usually concentrated in the first few output modes.

\subsection{Low-rank input-output bases}\label{sec:method_intro_pod}
\begin{figure}
    \centerline{
    \includegraphics[width=0.95\textwidth]{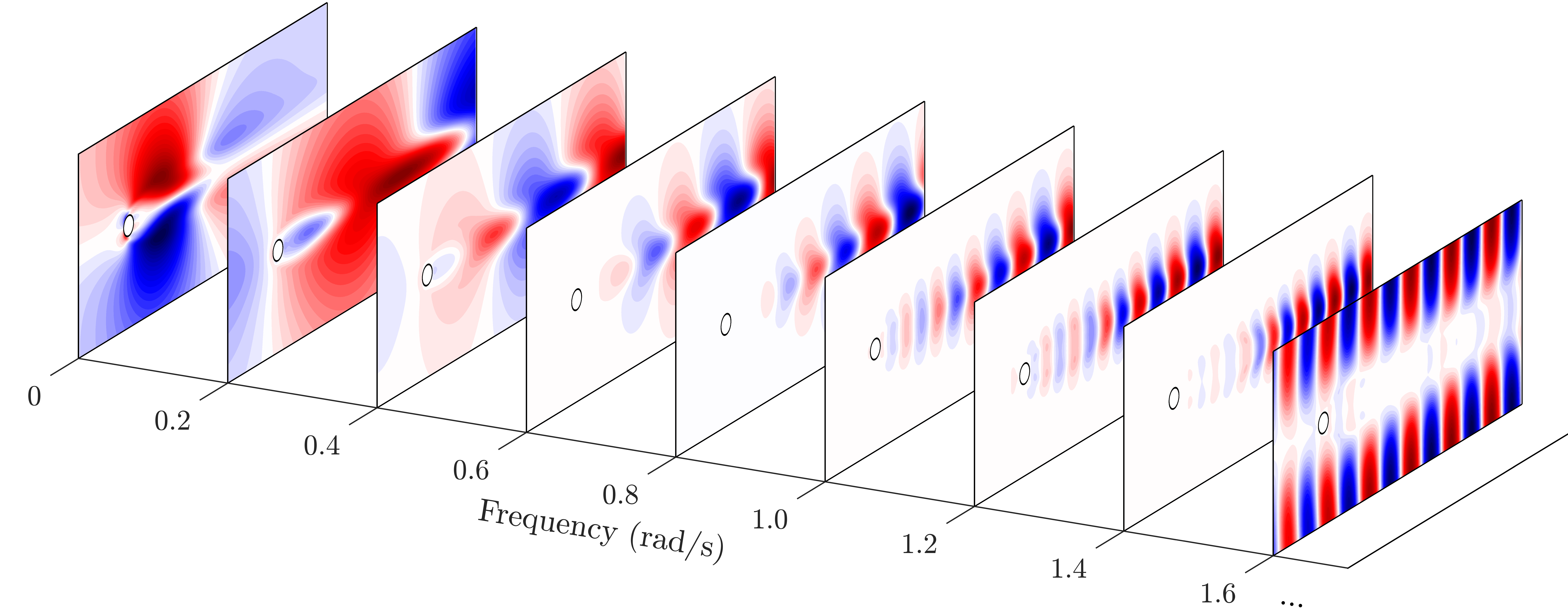}
    \llap{\parbox[b]{4.6in}{(a) Resolvent modes $\hat{\textbf{\textit{u}}}_1(\omega_i)$\\\rule{0ex}{1.75in}}}
    }
    \centerline{
    \includegraphics[width=0.95\textwidth]{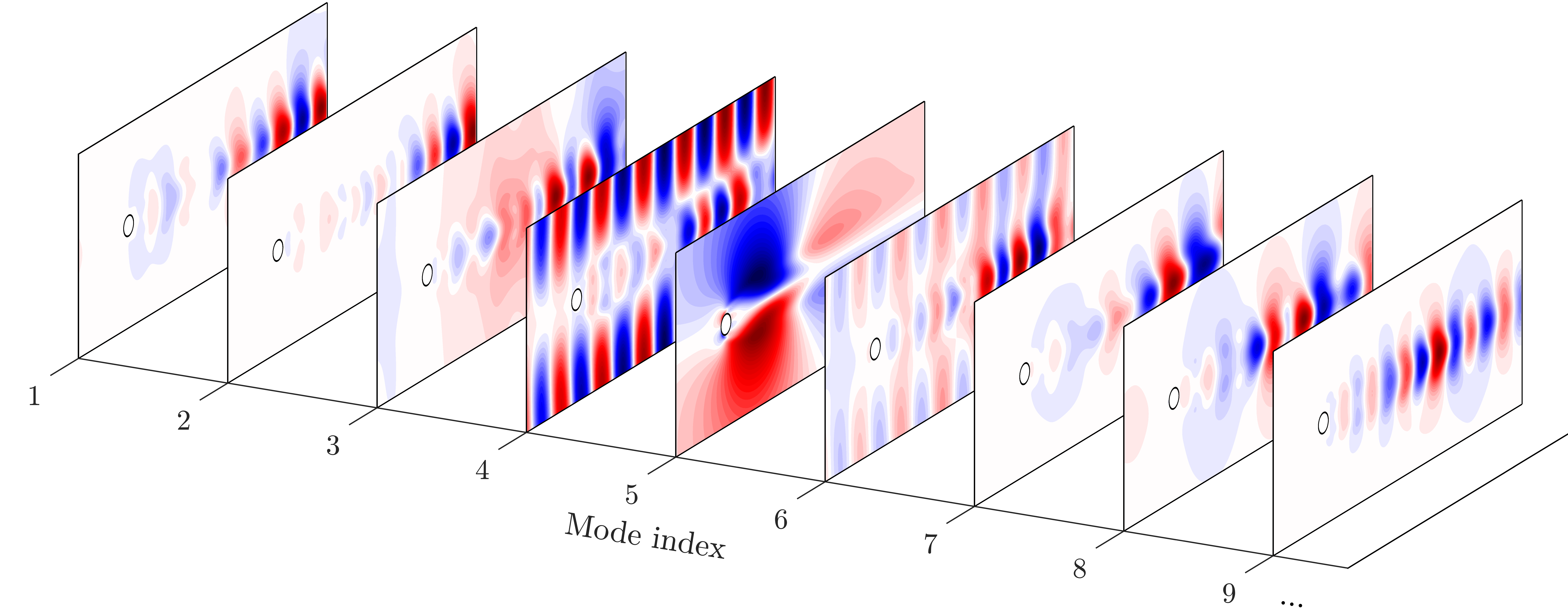}
    \llap{\parbox[b]{4.6in}{(b) POD modes $\Tilde{\textbf{\textit{u}}}_i$\\\rule{0ex}{1.75in}}}
    }
    \caption[Comparison of leading resolvent output modes and POD modes.]{(a) Leading resolvent output modes (transverse component, real part) at the sampled frequencies when $\Rey=90$. (b) The corresponding POD modes computed from the leading resolvent output modes. Note that these POD modes form an orthonormal basis that can reproduce resolvent modes via linear combinations.}
    \label{fig:linearsubspace_pod}
\end{figure}

In the previous section, we used truncated resolvent input and output matrices at one specified frequency (e.g.~the unstable frequency) to approximate the covariance and weight matrices appearing in the OE and FIC problems. This is because i) resolvent input and output modes are optimally ranked in terms of energy gain so that the first few modes capture the most energetic structures in fluid flows; ii) resolvent input-output modes are orthonormal with respect to the mass matrix and the truncated matrix forms a low-rank orthonormal basis that can be directly substituted into the Riccati equations. Note that it is important to use orthonormal bases because the original OE and FIC problems consider uncorrelated disturbances and independent measurements over the whole domain.

However, it is important to note that we aim to minimise the energy gain integrated over a wide frequency range, as defined by the cost function \eqref{equ:h2norm_resovent1}. The truncated input and output matrices $\hat{\textbf{F}}_k$ and $\hat{\textbf{U}}_k$ capture the most important features of the system at one specific frequency and may not be appropriate at other frequencies. This is most clearly shown in figure \ref{fig:linearsubspace_pod}(a) in which the leading resolvent output modes at different frequencies show significant differences. Therefore, we would like to construct a low-rank orthonormal basis that is appropriate across a range of frequencies. One way to achieve this is to perform proper orthogonal decomposition (POD) of the resolvent modes across a wide frequency range, which will result in an orthonormal basis that is ranked by the energy norm. For example, figure \ref{fig:linearsubspace_pod}(b) shows the POD modes computed from the leading resolvent output modes sampled at several frequencies, as shown in figure \ref{fig:linearsubspace_pod}(a). The leading resolvent output mode at each sampled frequency can be represented by a linear combination of these POD modes.

To construct the orthonormal basis, we first consider a series of truncated input matrices $\hat{\textbf{F}}\hspace{0mm}^{(\omega)}_k$ or output matrices $\hat{\textbf{U}}\hspace{0mm}^{(\omega)}_k$ at sampled frequencies $\omega\in \{\omega_1,\ \cdots,\ \omega_n\}$ as well as their complex conjugates. We choose the first $k$ resolvent modes from each input or output matrix so that the most important dynamics at each frequency are included. Therefore, we can build complex-valued matrices containing all information across the sampled frequencies:
\begin{subequations}\label{equ:complexmatpod}
    \begin{align}
        \hat{\textbf{H}}_{\text{F}}&=
        \begin{bmatrix}
            \hat{\textbf{F}}\hspace{0mm}^{(\omega_1)}_k&\hat{\textbf{F}}\hspace{0mm}^{(-\omega_1)}_k&\cdots&\hat{\textbf{F}}\hspace{0mm}^{(\omega_n)}_k&\hat{\textbf{F}}\hspace{0mm}^{(-\omega_n)}_k
        \end{bmatrix}
        \in \mathbb{C}^{N_u\times 2nk}\ ,\\
        \hat{\textbf{H}}_{\text{U}}&=
        \begin{bmatrix}
            \hat{\textbf{U}}\hspace{0mm}^{(\omega_1)}_k&\hat{\textbf{U}}\hspace{0mm}^{(-\omega_1)}_k&\cdots&\hat{\textbf{U}}\hspace{0mm}^{(\omega_n)}_k&\hat{\textbf{U}}\hspace{0mm}^{(-\omega_n)}_k
        \end{bmatrix}
        \in \mathbb{C}^{N\times 2nk}\ ,\\
        \hat{\textbf{H}}_{\rmSigma}&=\textbf{diag}
        \begin{bmatrix}
            \rmSigma\hspace{0mm}^{(\omega_1)}_k&\rmSigma\hspace{0mm}^{(-\omega_1)}_k&\cdots&\rmSigma\hspace{0mm}^{(\omega_n)}_k&\rmSigma\hspace{0mm}^{(-\omega_n)}_k
        \end{bmatrix}
        \in \mathbb{C}^{2nk\times 2nk}\ ,
    \end{align}
\end{subequations}
where $\textbf{diag}[\cdot]$ indicates a block-diagonal matrix. Now we want to find the smallest linear subspace of each matrix and then construct an orthonormal basis for this subspace. Notice that a complex vector $\bm{\psi}$ and its complex conjugate $\boldsymbol{\psi}^H$ are both linear combinations of the real part and the imaginary part with complex coefficients. That is, the linear span of these vectors has $span\{\boldsymbol{\psi},\boldsymbol{\psi}^H\}=span\{\Real(\boldsymbol{\psi}),\Imag(\boldsymbol{\psi})\}$. Thus, by replacing each complex conjugate pair of truncated matrices with a real-valued pair, one can build real-valued subspaces without losing any information:
\begin{subequations}\label{equ:accmuresolmodes}
    \begin{align}
    \textbf{H}_{\text{F}}&=
    \begin{bmatrix}
        \Real(\hat{\textbf{F}}\hspace{0mm}^{(\omega_1)}_k)&\Imag(\hat{\textbf{F}}\hspace{0mm}^{(\omega_1)}_k)&\cdots&\Real(\hat{\textbf{F}}\hspace{0mm}^{(\omega_n)}_k)&\Imag(\hat{\textbf{F}}\hspace{0mm}^{(\omega_n)}_k)
    \end{bmatrix}\hat{\textbf{H}}_{\rmSigma}^{1/2}
    \in \mathbb{R}^{N_u\times 2nk}\ ,\\
    \textbf{H}_{\text{U}}&=
    \begin{bmatrix}
        \Real(\hat{\textbf{U}}\hspace{0mm}^{(\omega_1)}_k)&\Imag(\hat{\textbf{U}}\hspace{0mm}^{(\omega_1)}_k)&\cdots&\Real(\hat{\textbf{U}}\hspace{0mm}^{(\omega_n)}_k)&\Imag(\hat{\textbf{U}}\hspace{0mm}^{(\omega_n)}_k)
    \end{bmatrix}\hat{\textbf{H}}_{\rmSigma}^{1/2}
    \in \mathbb{R}^{N\times 2nk}\ .
    \end{align}
\end{subequations}
Here, each subspace is weighted by the square root of the diagonal matrix $\hat{\textbf{H}}_{\rmSigma}$ and forms an eigenvalue problem of the symmetric, positive-semidefinite matrix (i.e.~$\textbf{H}_{\text{F}}^T\textbf{M}\textbf{H}_{\text{F}}$ or $\textbf{H}_{\text{U}}^T\textbf{Q}\textbf{H}_{\text{U}}$) with a smaller dimension $2nk$:
\begin{subequations}
    \begin{align}
        \textbf{H}_{\text{F}}^T\textbf{M}\textbf{H}_{\text{F}}\rmPhi_{\text{F}}&=\rmLambda^2_{\text{F}}\rmPhi_{\text{F}}\ ,\\
        \textbf{H}_{\text{U}}^T\textbf{Q}\textbf{H}_{\text{U}}\rmPhi_{\text{U}}&=\rmLambda^2_{\text{U}}\rmPhi_{\text{U}}\ ,
    \end{align}
\end{subequations}
where $\rmLambda_{\text{F}}$ and $\rmLambda_{\text{U}}$ are diagonal matrices with real and non-negative energy norms $\lambda_i$ in descending order, and the columns of $\rmPhi_{\text{F}}$ or $\rmPhi_{\text{U}}$ are real, orthonormal eigenvectors. The orthonormal projection bases (either with respect to $\textbf{M}$ or $\textbf{Q}$) of \eqref{equ:complexmatpod}, known as the POD modes, are given by
\begin{subequations}
    \begin{align}
        \Tilde{\textbf{F}}&=\textbf{H}_{\text{F}}\rmPhi_{\text{F}}\rmLambda_{\text{F}}^{-1}\in\mathbb{R}^{N_u\times 2nk}\ ,\\
        \Tilde{\textbf{U}}&=\textbf{H}_{\text{U}}\rmPhi_{\text{U}}\rmLambda_{\text{U}}^{-1}\in\mathbb{R}^{N\times 2nk}\ .
    \end{align}
\end{subequations}
These POD modes are sorted in descending order based on the magnitude of each $\lambda_i$, which indicate the importance of each mode.
A further truncation can be made so that we only need to consider the first $m$ POD modes for each case, as denoted by $\Tilde{\textbf{F}}_m$ and $\Tilde{\textbf{U}}_m$. In this study, we allow a relative mismatch of less than $10^{-6}$ across the specified frequency range and $m$ is chosen such that $\sum_{i=1}^{m}\lambda^2_i\approx (1-10^{-6})\ \rm{tr}(\rmLambda^2)$ to capture a sufficient fraction of the total energy norm. 

\subsection{Performance optimisation}\label{sec:method_intro_algorithm}
In general, using a low-rank basis results in a local optimum for the optimal design problem with full-state inputs and full-state outputs. This is because the low-rank basis that is good enough for the open-loop system may not appropriately approximate the most important input-output directions of the closed-loop system. Therefore, we need to perform POD not only for the open-loop system but also for the closed-loop system and iterate until the estimation or control performance converges.

We now present detailed steps for designing the optimal estimator, as summarised in algorithm \ref{code:algorithmic1}. The optimal full-state information controller can be designed in an analogous manner. Instead of full-velocity-state random disturbances of size $N_u$, we now consider reduced random disturbances $\textbf{\textit{d}}_m\in\mathbb{R}^m$ which disturb the system through random combinations of the input basis $\Tilde{\textbf{F}}_m$, i.e.~$\textbf{\textit{d}}=\Tilde{\textbf{F}}_m\textbf{\textit{d}}_m$. Thus, the covariance matrix in the Riccati equation is set to be $\textbf{W}^{1/2}=\textbf{P}\textbf{M}\Tilde{\textbf{F}}_m$. As discussed in \S\ref{sec:method_intro_pod}, the orthonormal input basis $\Tilde{\textbf{F}}_m$ is constructed from the matrix $\textbf{H}_{\text{F}}$ in which resolvent input modes with large energy amplification gains are stacked together (see equation (\ref{equ:accmuresolmodes}\textit{a})). The optimal estimator $\textbf{K}_{f}$ is designed to minimise the estimation error under the stimulation of resolvent input modes that give rise to the most energetic responses across the specified frequency range.

\newcommand{\algorithmicoe}{\textbf{OE problem:}}
\newcommand{\algorithmicfic}{\textbf{FIC problem:}}
\newcommand{\algorithmicioc}{\textbf{IOC problem:}}
\begin{algorithm}[ht]
\caption{Full-order Optimal Estimator Design}
\begin{algorithmic}[1]
\STATE Initialise. Set stopping criterion $\theta_c=10^{-4}$. Choose $n$ sampled frequencies.
\vspace{1mm}
\STATE Determine the resolvent input modes $\hat{\textbf{F}}_k$ of the open-loop system \eqref{equ:sec4.1_svd} at each sampled frequency (see \S\ref{sec:method_intro_resol_analysis}). Assemble the matrix $\textbf{H}_{\text{F}}$ defined by equation (\ref{equ:accmuresolmodes}\textit{a}).
\vspace{1mm}
\WHILE{$|\theta_t|>\theta_c$}
\vspace{1mm}
\STATE Determine the low-rank orthonormal basis $\Tilde{\textbf{F}}_m$ from the POD modes of $\textbf{H}_{\text{F}}$ (see $\S$\ref{sec:method_intro_pod}). The truncation criterion is chosen to be $10^{-6}$.
\vspace{1mm}
\STATE Set $\textbf{W}^{1/2}=\textbf{P}\textbf{M}\Tilde{\textbf{F}}_m$, solve the Riccati equation (see \S\ref{sec:intro_matrices}). Update the closed-loop error system with the estimator $\textbf{K}_{f}$:
\begin{equation}\nonumber
    \textbf{Z}_1(s)=\textbf{M}^{1/2}\textbf{P}^T(s\textbf{E}-\textbf{A}+\textbf{K}_{f}\textbf{C})^{-1}\textbf{P}(\textbf{M}^{1/2})^T\ .
\end{equation}
\STATE  Determine the resolvent input modes $\hat{\textbf{F}}_k$ of the closed-loop system at each sampled frequency (see \S\ref{sec:method_intro_resol_analysis}).\\
\vspace{1mm}
\STATE Form the weighted subspace $\textbf{H}_{\text{F}}^{\text{oe}}$ from the closed-loop resolvent results and update the matrix $\textbf{H}_{\text{F}}$ for computing POD modes:
\begin{equation}\nonumber
\begin{aligned}
    \textbf{H}_{\text{F}}&=[\sqrt{\textrm{Trace}\{\hat{\textbf{H}}\hspace{0mm}^{\text{oe}}_{\rmSigma}\}/\textit{m}}\ \Tilde{\textbf{F}}\hspace{0mm}_m\ \  \textbf{H}_{\text{F}}^{\text{oe}}]\ ,\\
    \textbf{H}_{\text{F}}^{\text{oe}}&=
    \begin{bmatrix}
        \cdots&\Real(\hat{\textbf{F}}_k)\ &\ \Imag(\hat{\textbf{F}}_k)&\cdots
    \end{bmatrix}(\hat{\textbf{H}}\hspace{0mm}^{\text{oe}}_{\rmSigma})^{1/2}\ .
\end{aligned}
\end{equation}
\STATE Evaluate the estimation performance (see \S\ref{sec:method_intro_resol_2norm}):
\begin{equation}\label{equ:performance_iteration}\nonumber
    \gamma^2_t(k,\omega_n)=\dfrac{1}{\pi}\int_{0}^{\omega_n}\sum_{i=1}^{k}\sigma^2_i(j\omega)\ d\omega+\norm{\textbf{Z}_2}_2^2\ .
\end{equation}
\STATE Compute the relative change of the performance at the current iteration
\begin{equation}\nonumber
    \theta_t=|\gamma^2_{t}/\gamma^2_{t-1}-1|\ .
\end{equation}
\ENDWHILE
\end{algorithmic}
\label{code:algorithmic1}
\end{algorithm}

However, disturbance reduction is performed based on the open-loop system, and a number of possible inputs discarded during the construction of the low-rank basis might become critical for the closed-loop estimation error system. Therefore, we update the matrix $\textbf{H}_{\text{F}}$ with the resolvent input modes evaluated for the closed-loop error system:
\begin{equation}
    \textbf{H}_{\text{F}}=[\sqrt{\textrm{Trace}\{\hat{\textbf{H}}\hspace{0mm}^{\text{oe}}_{\rmSigma}\}/\textit{m}}\ \Tilde{\textbf{F}}\hspace{0mm}_m\ \  \textbf{H}_{\text{F}}^{\text{oe}}]\ ,
\end{equation}
where matrices $\textbf{H}_{\text{F}}^{\text{oe}}$ and $\hat{\textbf{H}}\hspace{0mm}^{\text{oe}}_{\rmSigma}$ consist of the resolvent input modes and the corresponding singular values from the closed-loop error system. Here, the superscript $(\cdot)^{\text{oe}}$ denotes quantities from the closed-loop estimation error system.  Note that the low-rank basis vectors $\Tilde{\textbf{F}}\hspace{0mm}_m$ from the open-loop system are also stacked together in the matrix $\textbf{H}_{\text{F}}$ and scaled by the coefficient $\sqrt{\textrm{Trace}\{\hat{\textbf{H}}\hspace{0mm}^{oe}_{\rmSigma}\}/\textit{m}}$ so that it has the same total energy norm as $\textbf{H}_{\text{F}}^{\text{oe}}$. The optimal estimator designed for the updated low-rank basis is able to minimise the estimation error not only for the original open-loop system but also for the closed-loop system.

The steps of the iterative procedure are summarised in algorithm \ref{code:algorithmic1}. We aim to iteratively optimise the estimation performance by minimising the integrated energy gains of the closed-loop error system (see \S\ref{sec:method_intro_resol_2norm}). We monitor the relative change of the estimation performance $\gamma^2(k,\omega_n)$ at each iteration until it reaches the stopping criterion $\theta_c$. By choosing a sufficient number of resolvent modes over a sufficiently large frequency range, it should result in an optimal estimator that converges to the global optimal performance when disturbances are applied everywhere. 

\section{Numerical Results}\label{sec:results}
In order to test our numerical method for the optimal estimation and control of fluid flows, we consider linear feedback control of the cylinder wake at $\Rey=90$ as an example. Specifically, a single sensor is placed downstream at $\textbf{\textit{x}}_s=(7.70,\ 0.73)$ measuring the perturbation velocity $\textbf{\textit{u}}'(\textbf{\textit{x}}_s, t)$ and a single actuator is placed near the cylinder surface (at $\textbf{\textit{x}}_a=(2.56,\ 1.18)$) stimulating the flow in the streamwise and transverse directions independently. This section is organised as follows. The validation of low-rank bases and a convergence analysis of algorithm \ref{code:algorithmic1} are presented in \S\ref{sec:SubspaceAndConvergence} using different parameter choices. The implementation of the feedback control design is then presented, which consists of three parts: i) optimal estimator design when only a single sensor is available for measurement; ii) optimal full-state information controller design when only a single actuator is available for control ; iii) feedback control with a single sensor for measurement and a single actuator for control. 

\subsection{Validation and convergence analysis}\label{sec:SubspaceAndConvergence}
\begin{figure}
    \centerline{
    \hspace{0mm}
    \includegraphics[width=0.47\textwidth]{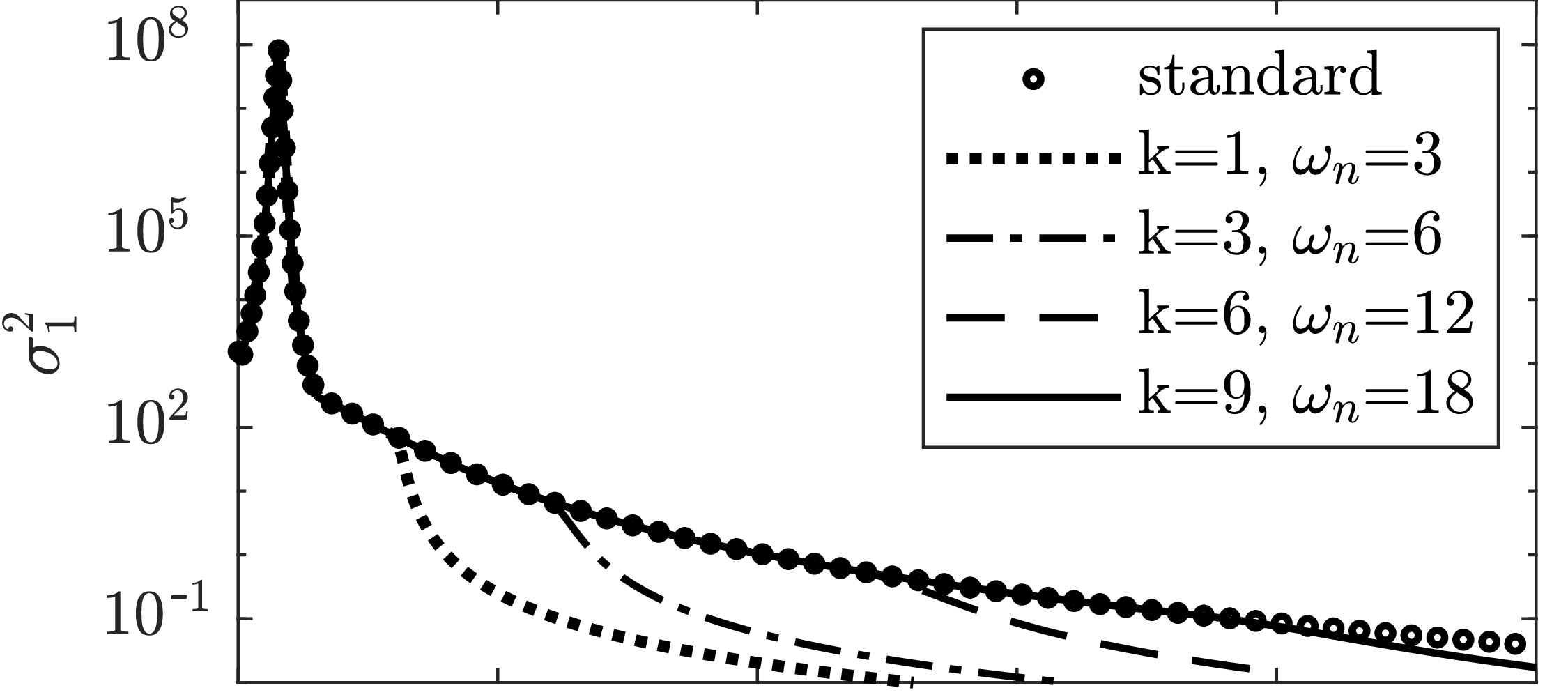}
    \llap{\parbox[b]{2.325in}{(a)\\\rule{0ex}{0.9in}}}
    \hspace{0mm}
    \includegraphics[width=0.47\textwidth]{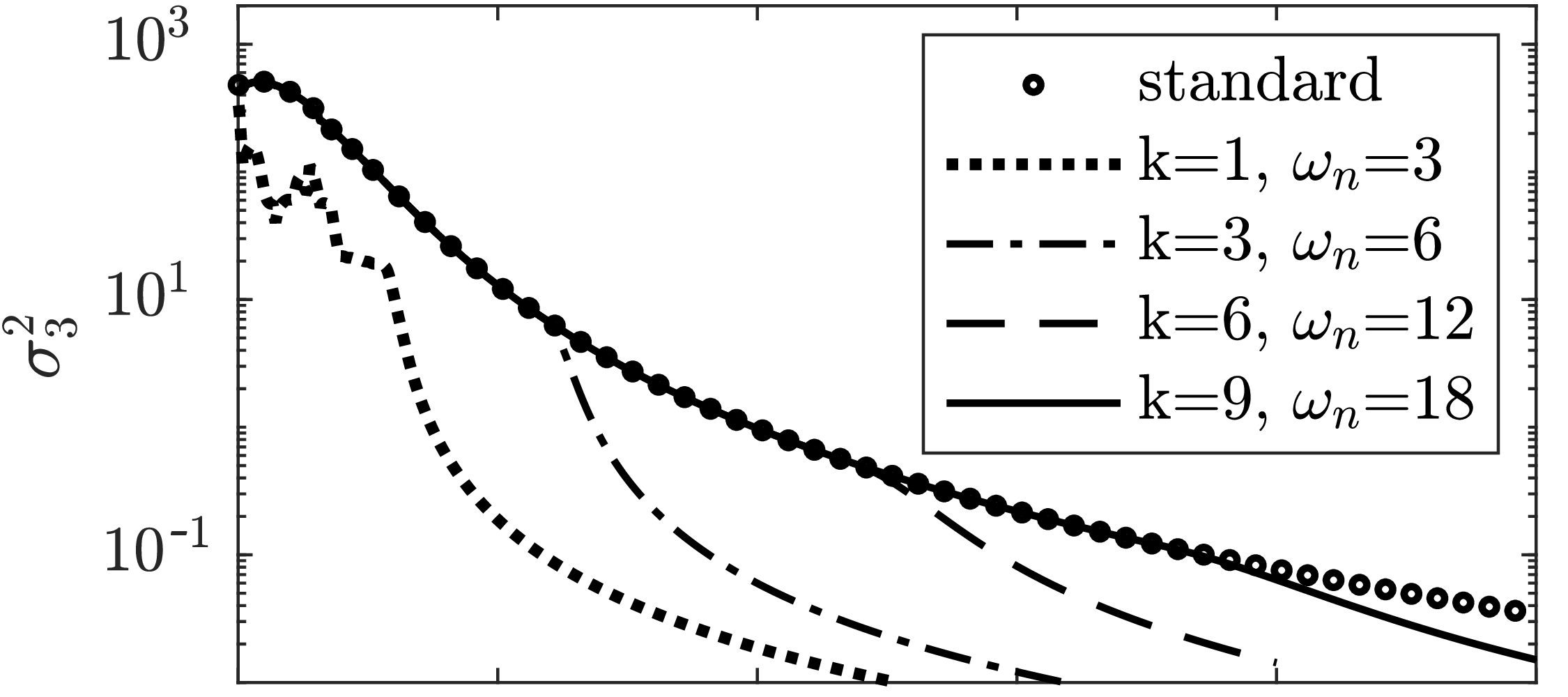}
    \llap{\parbox[b]{2.325in}{(b)\\\rule{0ex}{0.9in}}}
    }
    \vspace{0mm}
    \centerline{
    \hspace{0mm}
    \includegraphics[width=0.47\textwidth]{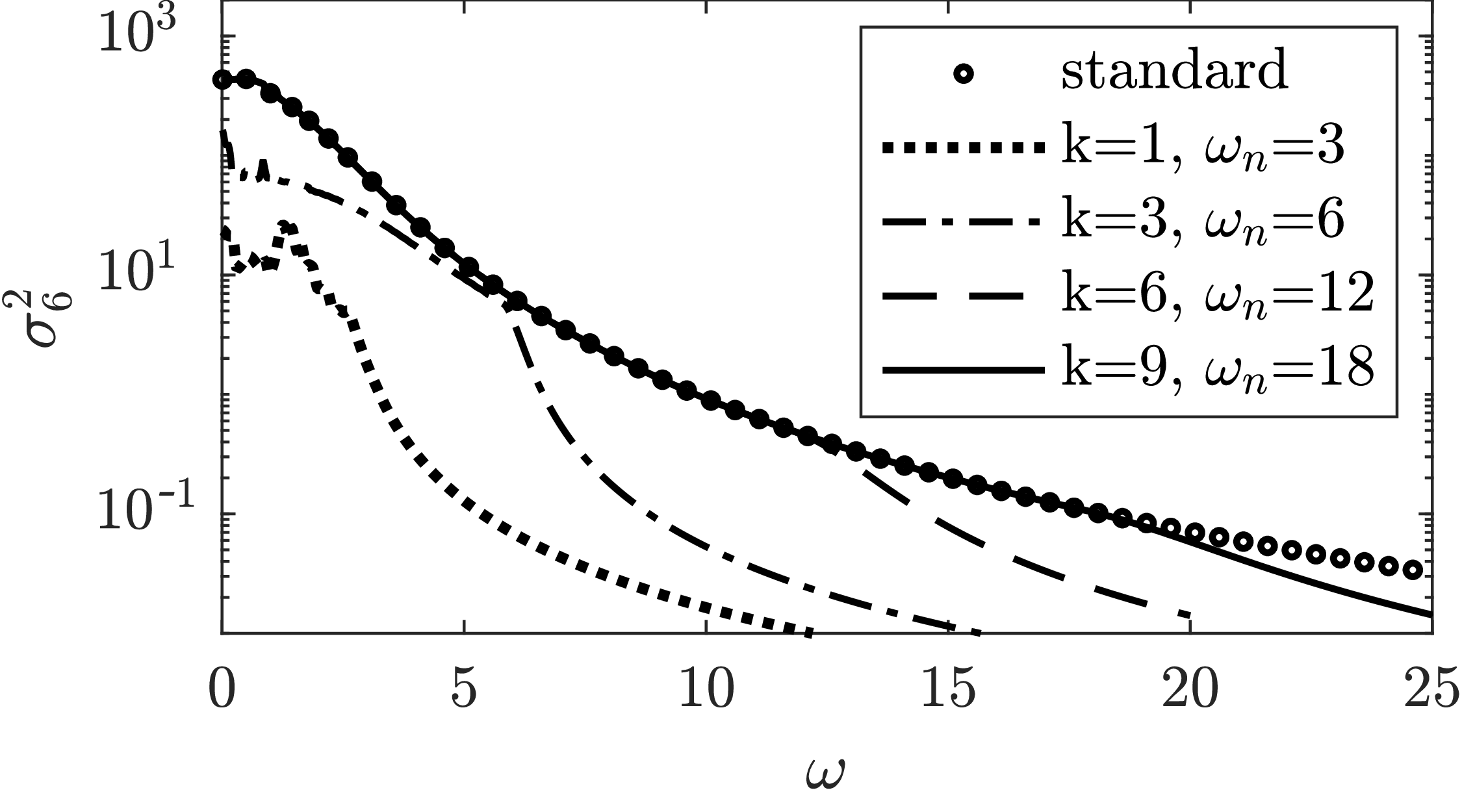}
    \llap{\parbox[b]{2.325in}{(c)\\\rule{0ex}{1.175in}}}
    \hspace{0mm}
    \includegraphics[width=0.47\textwidth]{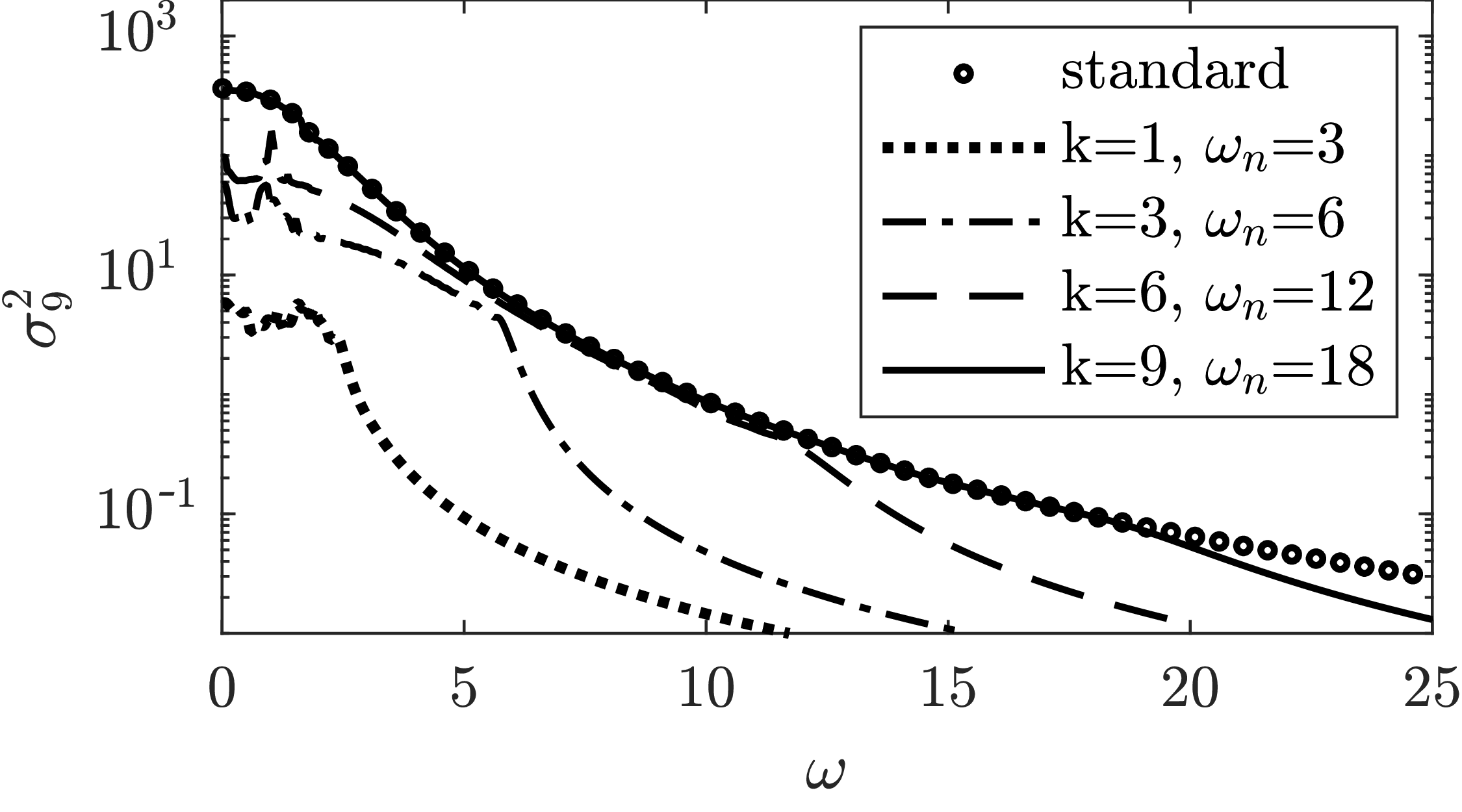}
    \llap{\parbox[b]{2.325in}{(d)\\\rule{0ex}{1.175in}}}
    }
    \caption[Comparison of resolvent spectra between the standard perturbation system and the terminal-reduced systems]{Comparison of (a) the first singular value $\sigma_1$, (b) the third singular value $\sigma_3$, (c) the sixth singular value $\sigma_6$ and (d) the ninth singular value $\sigma_9$ as a function of frequency. They are computed either from the standard perturbation system ($\circ$) or from the terminal-reduced systems (lines). The terminal reduction is performed using different choices of $k$ and $\omega_n$ which correspond to including the first $k$ resolvent input modes within the frequency range $[-\omega_n,\ \omega_n]$.}
    \label{fig:lqe_valid_sigs}
\end{figure}
The algorithm for solving the OE and FIC problems is based on an iterative strategy. Each step consists of the projection of the original system onto a low-rank orthonormal basis and gives a terminal-reduced system that retains only the most energetic mechanisms. For the sake of validation, we consider the terminal-reduced system's ability to capture the dynamics of the original system by comparing their resolvent spectra. For brevity, we consider input reduction of the open-loop system \eqref{equ:sec4.1_svd} where the size of the random disturbances is reduced to $m$ after setting the matrix $\textbf{W}^{1/2}=\textbf{P}\textbf{M}\Tilde{\textbf{F}}_m$. The construction of the low-rank orthonormal basis $\Tilde{\textbf{F}}_m$ relates to two parameters: $k$ and $\omega_n$ which correspond to including the first $k$ resolvent input modes over the frequency range $[-\omega_n, \omega_n]$. With a good choice of parameters, any largely amplified linear inputs, as well as their energy gains, can be well modelled by the terminal-reduced system.

Figure \ref{fig:lqe_valid_sigs} compares resolvent spectra from the standard perturbation system ($\circ$) to those from the terminal-reduced systems (lines) using four different parameter choices. As can be seen from the figure, the terminal-reduced system is better able to capture the dynamics of the original system for a wider frequency range if larger values of $k$ and $\omega_n$ are chosen. In particular, the comparisons from figure \ref{fig:lqe_valid_sigs} indicate two important features of the terminal-reduced system. First of all, the first $k$ resolvent spectra from the terminal-reduced system well match those from the original system within the frequency range $\omega \in [0,\ \omega_n]$ (i.e.~relative error less than $10^{-6}$). Second, the remaining resolvent spectra, which correspond to resolvent input modes that are excluded in the terminal-reduction procedure, poorly fit those from the original system. For example, the dash-dotted lines, which correspond to the parameter choice $k=3$ and $\omega_n=6\ rad/s$, match the circle markers reasonably well for the first and third resolvent spectra within the frequency range $[0,6]\ rad/s$, as shown in figure \ref{fig:lqe_valid_sigs}($\textit{a}$, $\textit{b}$). Nevertheless, in figure \ref{fig:lqe_valid_sigs}($\textit{c}$, $\textit{d}$), which shows a comparison of the sixth spectra $\sigma^2_6$ and the ninth spectra $\sigma^2_9$, the dash-dotted lines in no case are following the circle markers. 
\begin{figure}
    \centerline{
    \hspace{0mm}
    \includegraphics[width=0.47\textwidth]{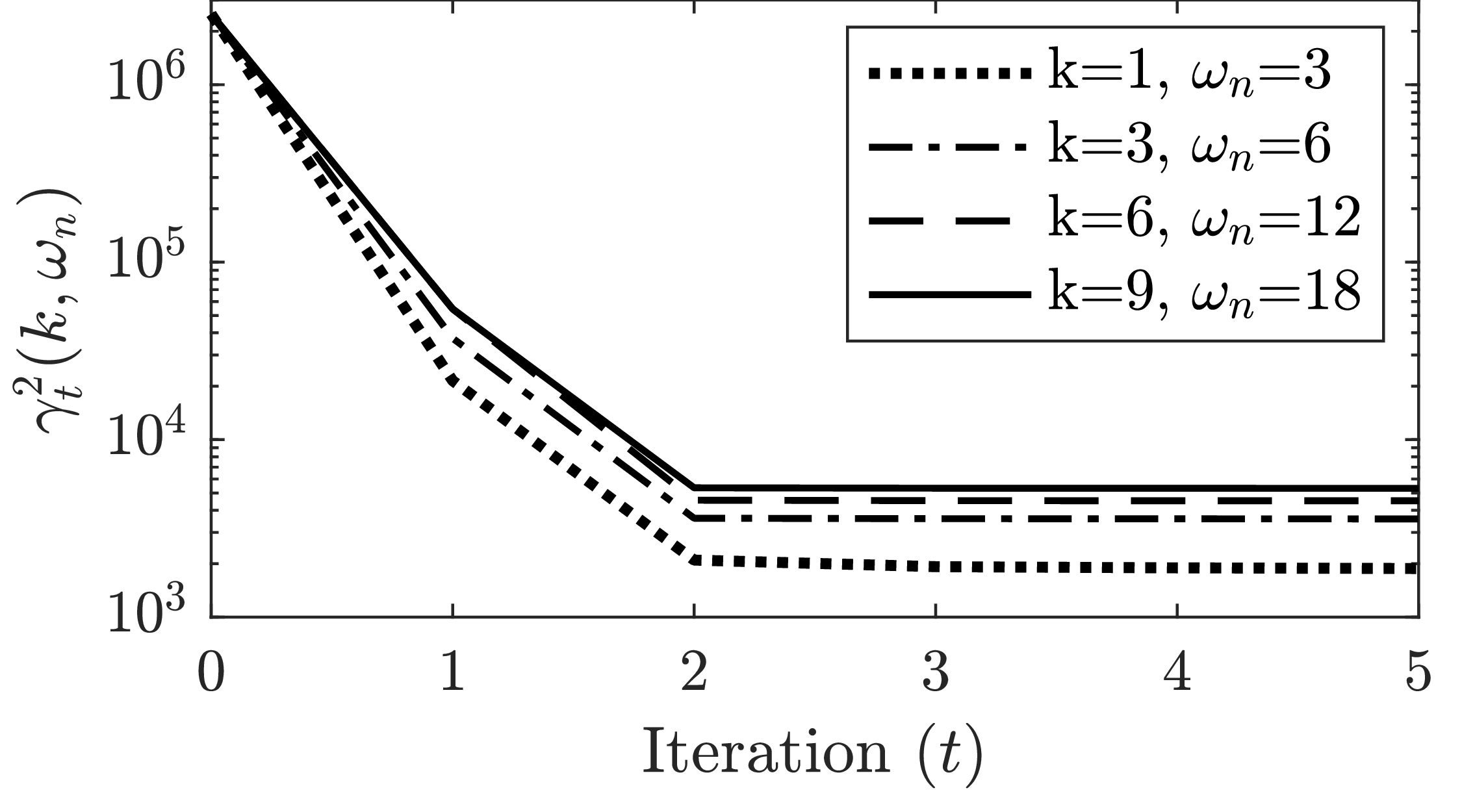}
    \llap{\parbox[b]{2.325in}{(a)\\\rule{0ex}{1.175in}}}
    \hspace{0mm}
    \includegraphics[width=0.47\textwidth]{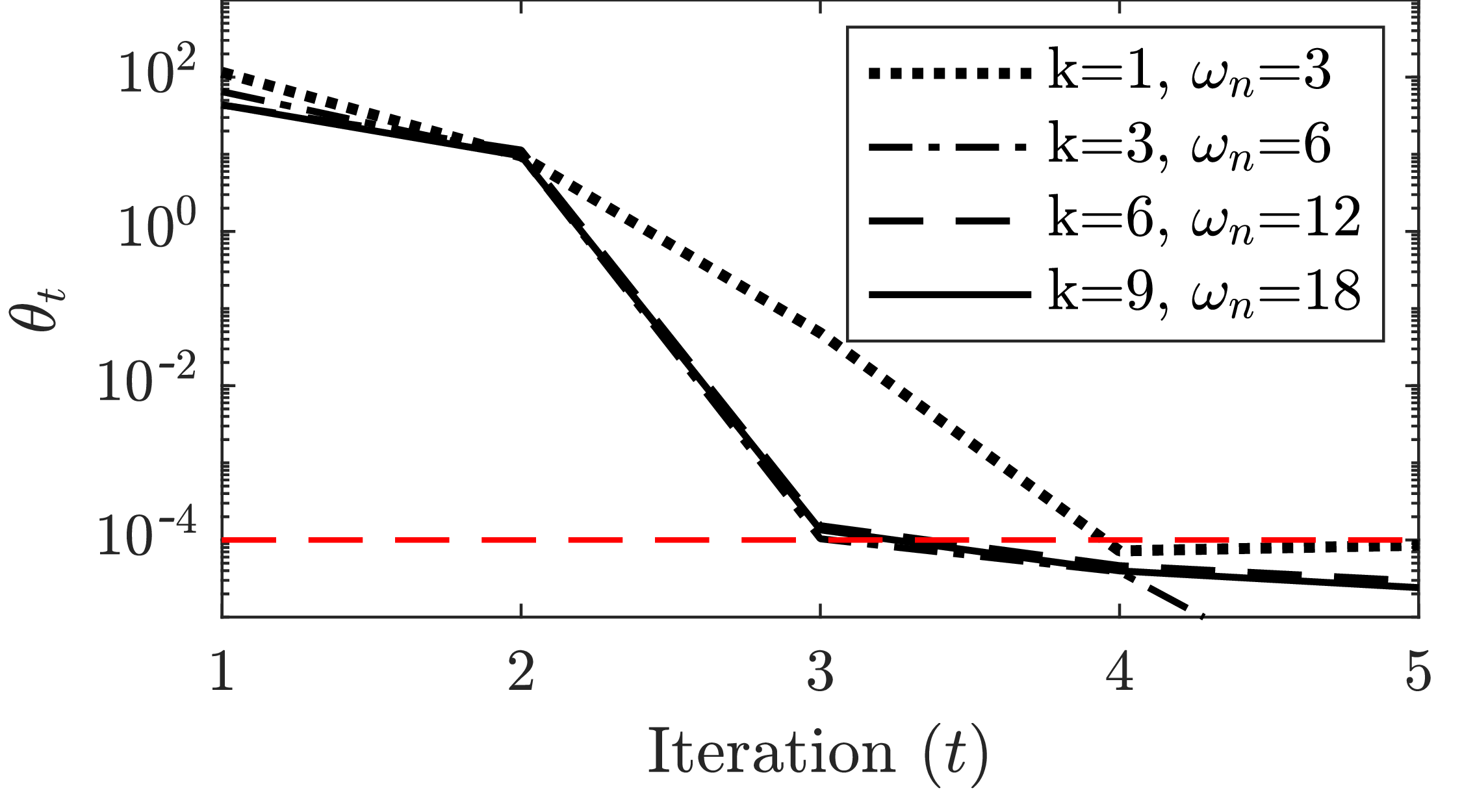}
    \llap{\parbox[b]{2.325in}{(b)\\\rule{0ex}{1.175in}}}
    }
    \caption[Convergence of the iterative algorithm using four different parameter choices]{Convergence of the iterative algorithm \ref{code:algorithmic1} with four different parameter choices. (a) The convergence of the cost function $\gamma^2_t(k,\omega_n)$. (b) The convergence of the relative change of the cost function $\theta_t$. The red dashed line indicates the chosen stop criterion $\theta_c$.}
    \label{fig:lqe_valid_conv}
\end{figure}

We further consider the implementation of algorithm \ref{code:algorithmic1} using the four above-mentioned parameter choices for convergence analysis. In figure \ref{fig:lqe_valid_conv}(\textit{a}, \textit{b}), the cost function $\gamma^2_t(k,\omega_n)$ and the corresponding relative change $\theta_t$ are shown as a function of the iteration number $t$. We notice that different choices of $k$ and $\omega_n$ show similar convergence features and each cost function converges to a constant after only two iterations. Their relative errors drop below the stop criterion $\theta_c=10^{-4}$ after four iterations, which is deemed sufficient to guarantee the converged estimator gain $\textbf{K}_f$. As already mentioned, each iteration includes a resolvent-based terminal reduction procedure (see \S\ref{sec:method_intro_pod}) as well as iteratively solving a nonlinear Riccati equation \cite{SaaKB19-mmess-2.0}. Choosing either a higher-value truncation criterion (e.g.~capture more than 99.5$\%$ energy norm while truncating POD matrices) or a lower Riccati solver tolerance will lead to further decreasing of the relative error $\theta_t$.

The results of the output reductions, as well as the convergence features of the FIC problem, are similar to those presented above, in which the convergence of the controller is always guaranteed after four iterations. The choice of parameters $k$ and $\omega_n$ is critical to the design of optimal estimators and controllers. We will show further details about their effects on estimation and control performance in the following sections. 

\subsection{Optimal estimator design}
\subsubsection{Parameters and estimation performance}
\begin{table}
  \begin{center}
\def~{\hphantom{0}}
  \begin{tabular}{p{1cm}p{1cm}p{1cm}p{1cm}p{1.5cm}p{1.5cm}}
  \toprule
      Case& $k$& $\omega_n$ & $m$  &$\gamma^2{(k,\omega_n)}$& $\textbf{\textit{J}}_{\text{OE}}$ \\[3pt]
    \midrule
       (a)& 1 & 3 & 29 & 1878.28 & 21170.07 \\
       (b)& 3 & 6 & 155 & 3574.74 & 21101.26 \\
       (c)& 6 & 12 & 446 & 4522.91 & 21094.52 \\ 
       (d)& 9 & 18 & 768 & 5316.83 & 21094.01 \\ 
    \bottomrule
  \end{tabular}
  \caption[A summary of results from the OE problem with different parameter choices]{The rank of the reduced disturbance $m$, the cost function for the optimal estimator design $\gamma^2{(k,\omega_n)}$ and the mean energy of the total estimation error $\textbf{\textit{J}}_{\text{OE}}$ (random disturbances applied everywhere) are listed for different parameter choices ($k$, $\omega_n$).  The optimal estimators are designed at $\Rey=90$ with a sensor placed at $\textbf{\textit{x}}_s=(7.70,0.73)$ and sensor noise of magnitude $\alpha=10^{-4}$.}
  \label{tab:lqe_valid_optjoe}
  \end{center}
\end{table}
We now consider optimal estimator design at $\Rey=90$ by implementing algorithm \ref{code:algorithmic1} with four different choices of the parameters $k$ and $\omega_n$. The sensor noise has a negligible magnitude of $\alpha=10^{-4}$. With this small choice of $\alpha$, the performance of the estimator remains insensitive to the sensor noise and thus we are able to analyse the effects of disturbances on the estimation error. Table \ref{tab:lqe_valid_optjoe} lists the rank of the reduced disturbances $m$, the cost function for the optimal estimator design $\gamma^2{(k,\omega_n)}$ and the resulting estimation performance $\textbf{\textit{J}}_{\text{OE}}$ for different choices of $k$ and $\omega_n$. The cost function $\gamma^2{(k,\omega_n)}$ is a norm squared when considering a limited number of system inputs and outputs (i.e.~low-rank inputs and outputs) whereas $\textbf{\textit{J}}_{\text{OE}}$ represents the estimation error when random disturbances are applied everywhere in the domain (i.e.~full-rank inputs and outputs).

We immediately see that the rank of the reduced disturbances $m$ increases with increasing values of $k$ and $\omega_n$. The physical meaning of the cost function $\gamma^2{(k,\omega_n)}$ is the mean energy of the estimation error while the system is excited by the first $k$ resolvent input modes across the frequency range $[-\omega_n,\ \omega_n]$. As expected, a larger rank of the reduced disturbance $m$, i.e.~considering more input directions, gives rise to a larger cost function $\gamma^2{(k,\omega_n)}$. By performing numerical simulations, we are able to evaluate the estimation performance $\textbf{\textit{J}}_{\text{OE}}$ based on the definition of the original OE problem: the mean energy of the estimation error when random disturbances are applied everywhere in the domain. With estimators designed with larger values of $k$ and $\omega_n$, we observe that the total estimation error $\textbf{\textit{J}}_{\text{OE}}$ decreases and eventually converges to a constant (relative change converges around $10^{-4}$). 

\begin{figure}
    \centerline{
    \hspace{2mm}
    \includegraphics[width=0.46\textwidth]{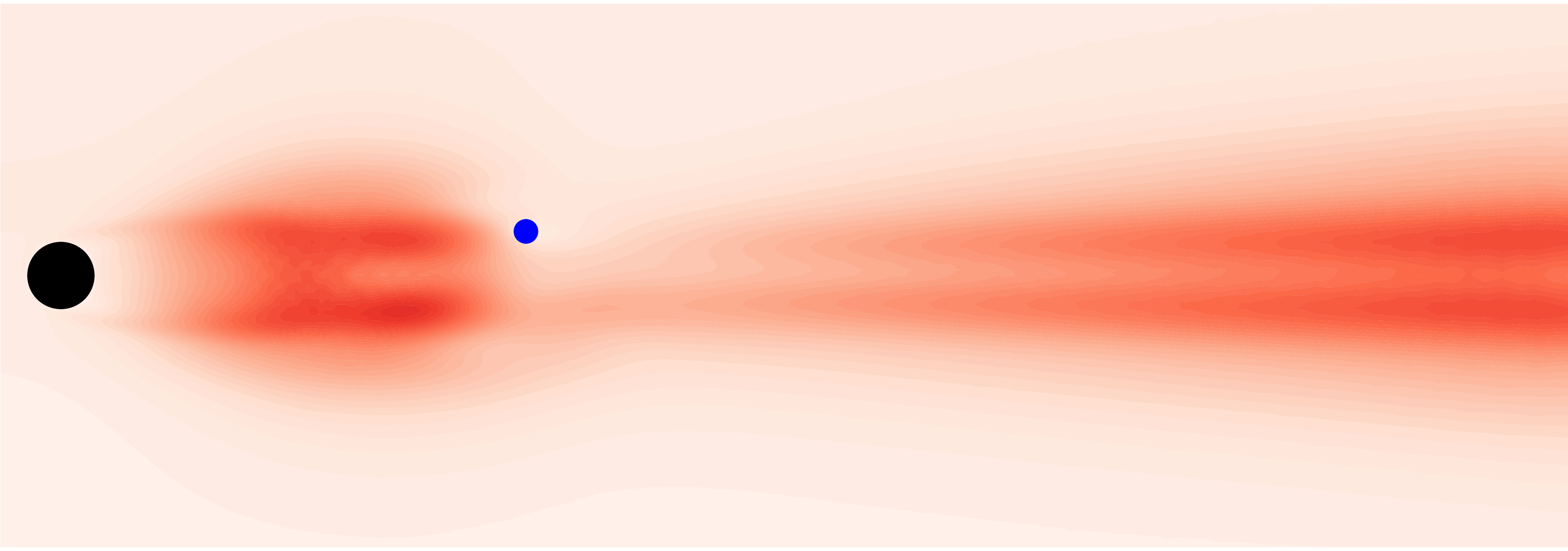}
    \llap{\parbox[b]{2.35in}{(a)\\\rule{0ex}{0.7in}}}
    \hspace{2mm}
    \includegraphics[width=0.46\textwidth]{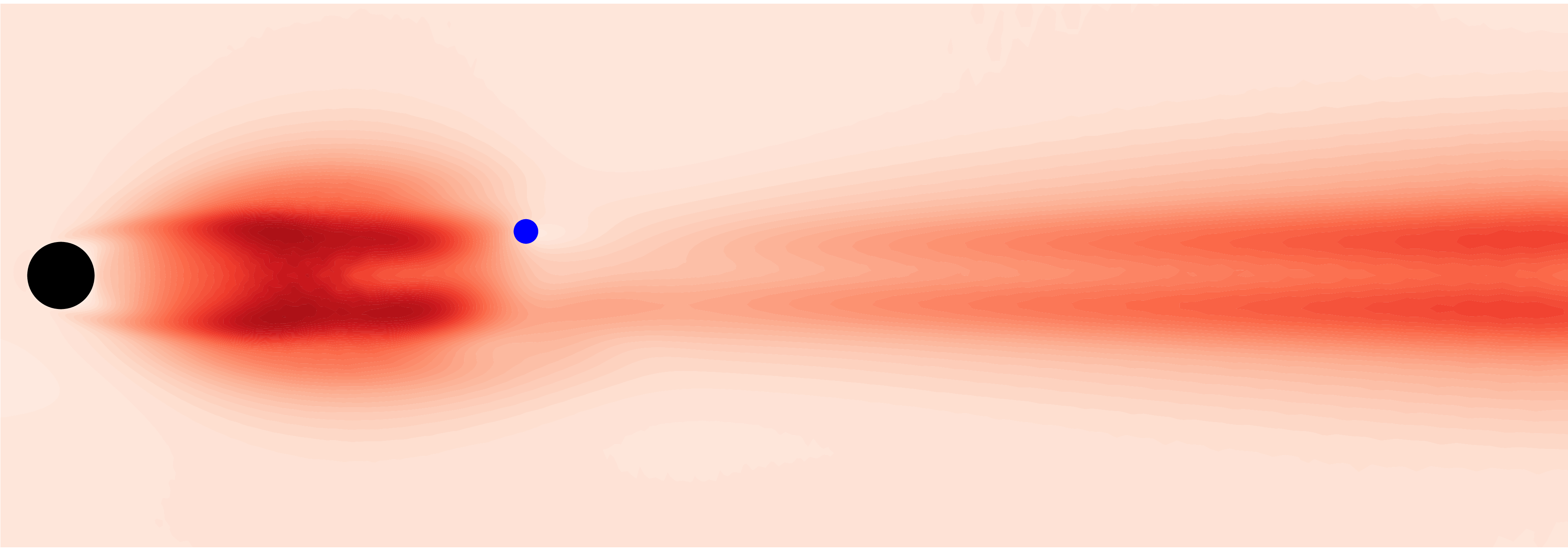}
    \llap{\parbox[b]{2.375in}{(b)\\\rule{0ex}{0.7in}}}
    }
    \centerline{
    \hspace{2mm}
    \includegraphics[width=0.46\textwidth]{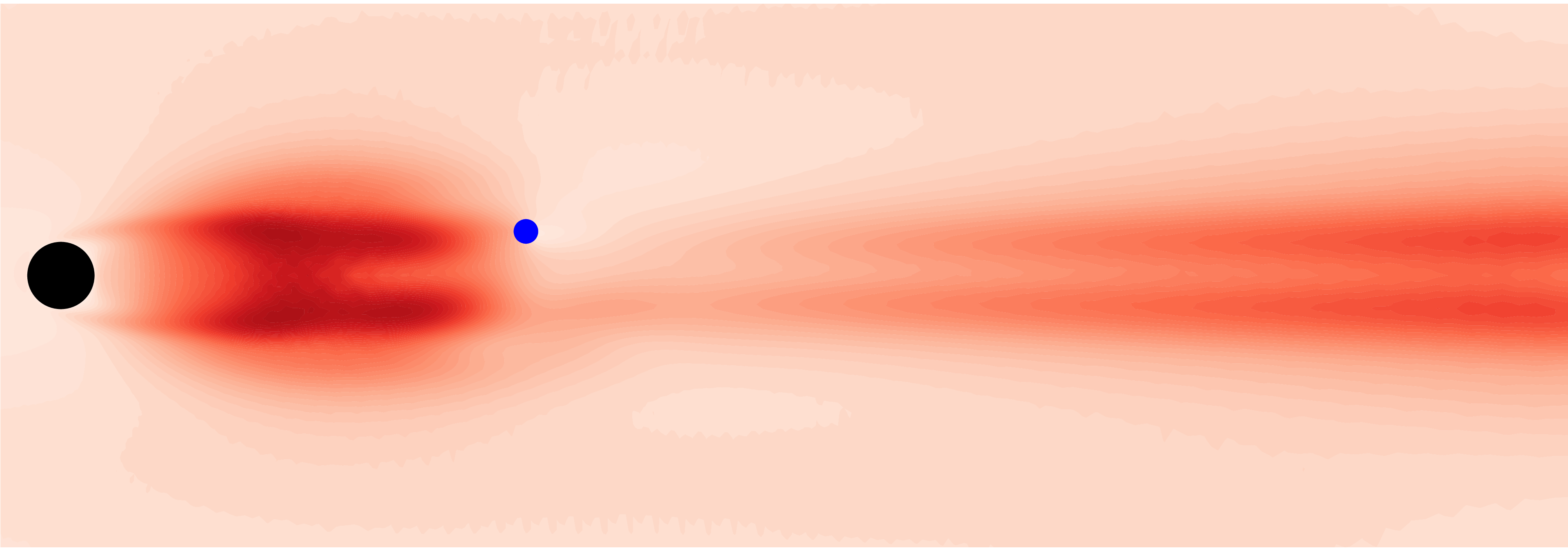}
    \llap{\parbox[b]{2.35in}{(c)\\\rule{0ex}{0.7in}}}
    \hspace{2mm}
    \includegraphics[width=0.46\textwidth]{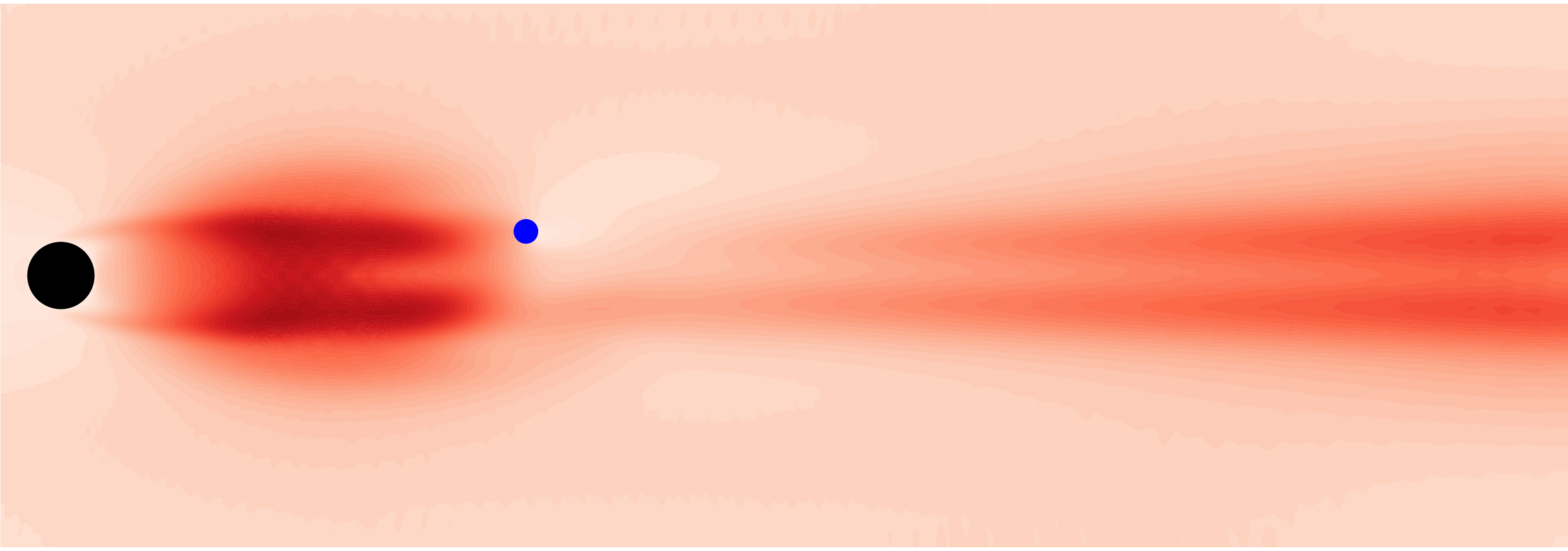}
    \llap{\parbox[b]{2.375in}{(d)\\\rule{0ex}{0.7in}}}
    }
    \centerline{
    \hspace{2mm}
    \includegraphics[width=0.46\textwidth]{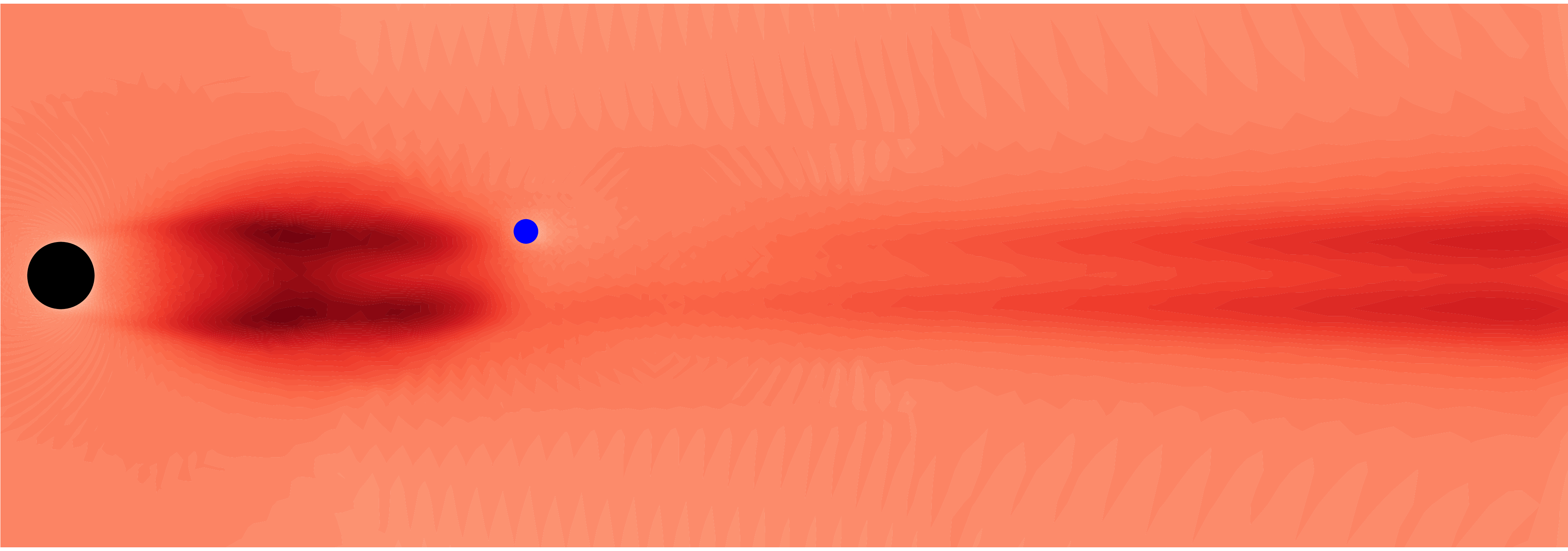}
    \llap{\parbox[b]{2.35in}{(e)\\\rule{0ex}{0.7in}}}
    \hspace{2mm}
    \includegraphics[width=0.46\textwidth]{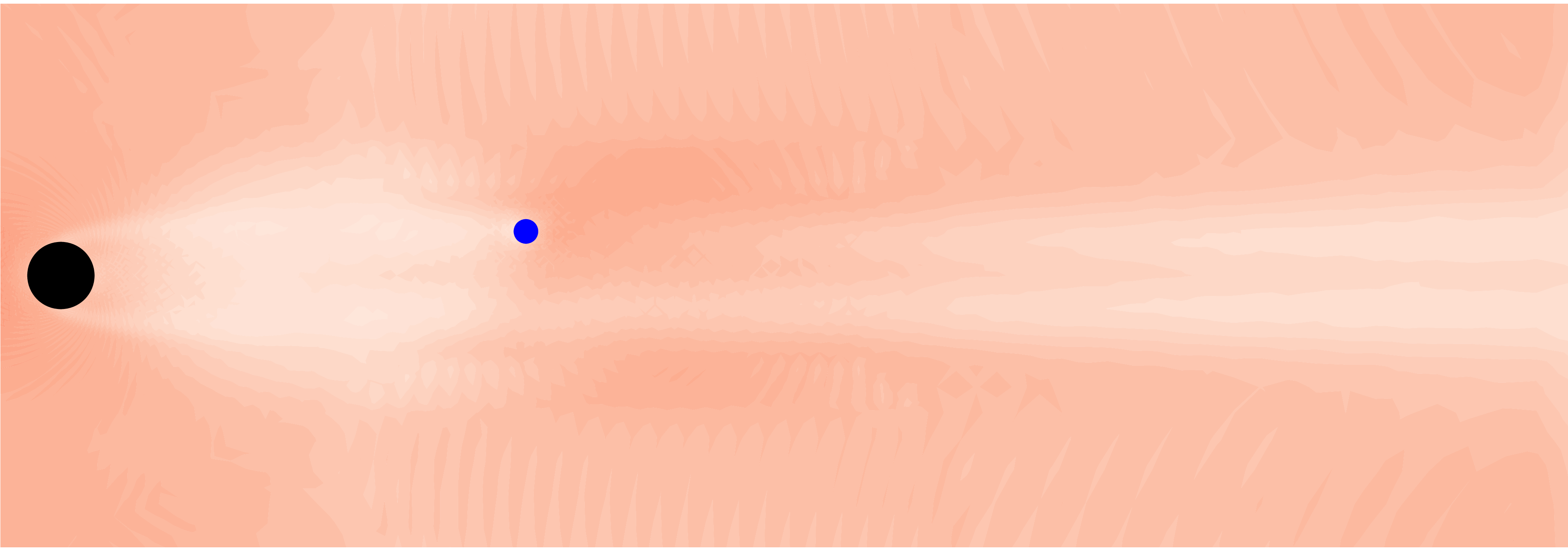}
    \llap{\parbox[b]{2.375in}{(f)\\\rule{0ex}{0.7in}}}
    }
    \caption[Spatial distribution of the estimation error $\epsilon_{\text{OE}}(x,y)$ for optimal estimators designed using different parameters]{Spatial distribution of the estimation error $\epsilon_{\text{OE}}(x,y)$ with a sensor placed at $(7.70,\ 0.73)$, which is marked by a blue dot ($\protect\bluedot$). The optimal estimators are designed using parameters (a) $k=1$,  $\omega_n=3$; (b) $k=3$, $\omega_n=6$; (c) $k=9$,  $\omega_n=12$ and (d) $k=9$, $\omega_n=18$. (e) The estimation error $\epsilon_{\text{OE}}$ when disturbances are applied everywhere over the domain. (\textit{f}) The differences of the estimation error $\epsilon_{\text{OE}}$ between (d) and (e) (i.e.~between $\gamma^2{(k,\omega_n)}$ and $\textbf{\textit{J}}_{\text{OE}}$). All plots share the same linear colour scale $[0,\ 10.5]$.}
    \label{fig:lqe_valid_errdist}
\end{figure}

It is also insightful to see the energy distribution of the estimation error throughout the domain and the effect of the disturbances applied. Thus, we define a root mean square value $\epsilon_{\text{OE}}(x,y)$ such that
\begin{equation}\label{equ:lqe_epsilon}
    \gamma^2{(k,\omega_n)}=\int_{\Omega}\epsilon^2_{\text{OE}}(x,y)\ d\Omega\ ,
\end{equation}
where $\epsilon_{\text{OE}}(x,y)$ denotes the contribution of the estimation error throughout the domain (see Appendix \ref{sec:app.a}). The contour plots of $\epsilon_{\text{OE}}(x,y)$ for the four cases listed in table \ref{tab:lqe_valid_optjoe} are shown in figure \ref{fig:lqe_valid_errdist}(a, b, c, d), respectively. We notice that these contour plots show quite similar structures. In particular, the most significant contributions to the estimation error are concentrated in two horizontal streaks which are approximately symmetric and located around $y\approx 0.7$. As expected, the smallest $\epsilon_{\text{OE}}$ occurs at the sensor location ($\protect\bluedot$) which divides the two streaks into two regions: a near-wake area (between the cylinder and the sensor) and a far-wake area (downstream of the sensor). The only difference among these plots is the slightly larger magnitude of $\epsilon_{\text{OE}}$ for cases with larger $k$ and $\omega_n$. This can be most clearly seen by comparing figure \ref{fig:lqe_valid_errdist}(a) to the remaining three plots. We further consider the spatial distribution $\epsilon_{\text{OE}}$ from numerical simulations when disturbances are applied everywhere, i.e.~$\epsilon_{\text{OE}}$ defined using $\textbf{\textit{J}}_{\text{OE}}$, which is shown in figure \ref{fig:lqe_valid_errdist}(e). Here, we use the same estimator as that used in case (d) but we apply random disturbances everywhere in the domain instead of the reduced disturbances of rank $m$. By comparing figure \ref{fig:lqe_valid_errdist}(d, e) we see that the reduced disturbances produce lower estimation error in the freestream but give almost the same characteristic structures (streaks described above) as the case with disturbances applied everywhere. This can be more clearly seen in figure \ref{fig:lqe_valid_errdist}(f), which plots the difference between figure \ref{fig:lqe_valid_errdist}(d) and (e): the region in which the streaks lie is white in space whereas the `background' or `freestream' area is of approximately uniform colour. 

\begin{figure}
    \centerline{
    \hspace{0mm}
    \includegraphics[width=0.47\textwidth]{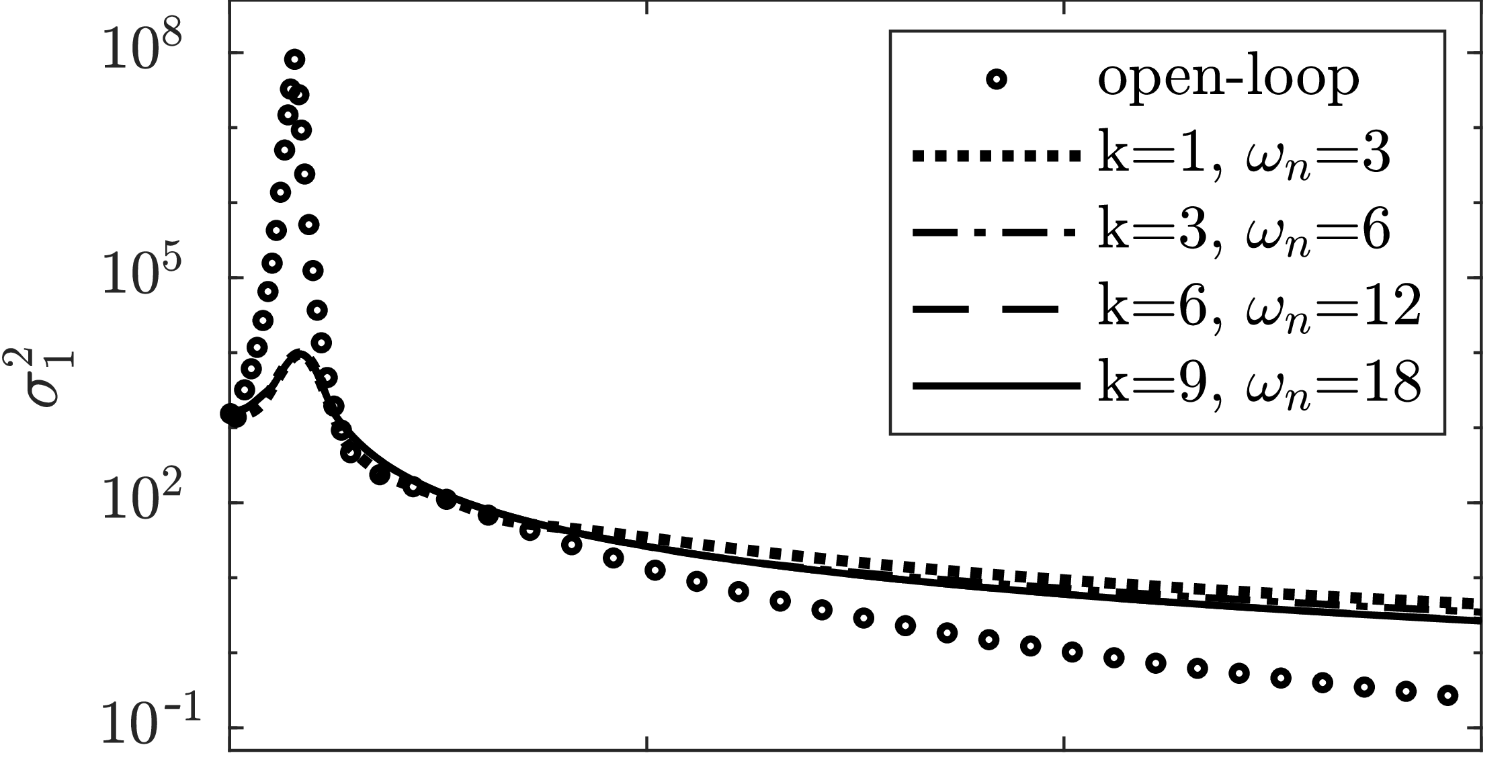}
    \llap{\parbox[b]{2.275in}{(a)\\\rule{0ex}{1.05in}}}
    \hspace{0mm}
    \includegraphics[width=0.47\textwidth]{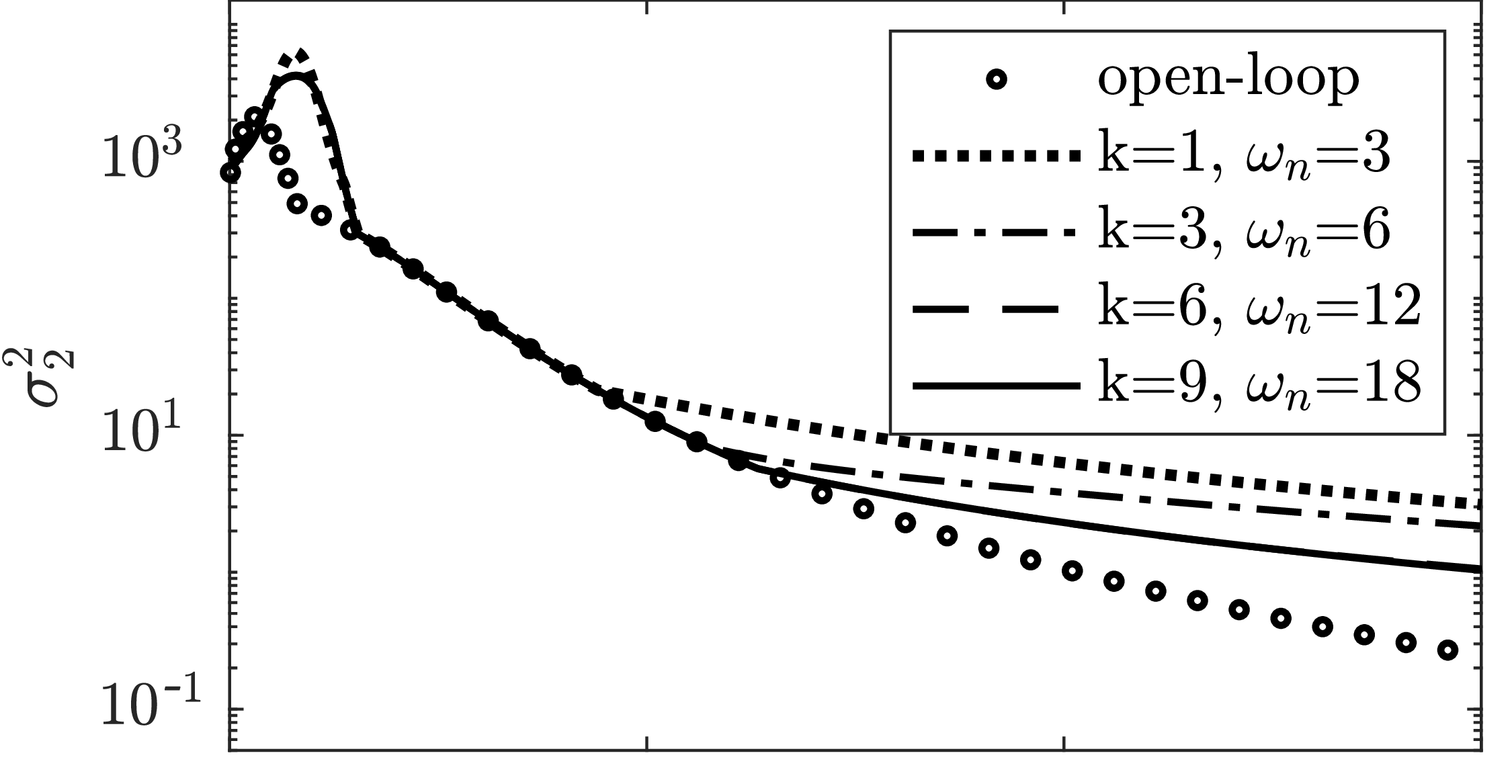}
    \llap{\parbox[b]{2.325in}{(b)\\\rule{0ex}{1.05in}}}
    }
    \vspace{0mm}
    \centerline{
    \hspace{0mm}
    \includegraphics[width=0.47\textwidth]{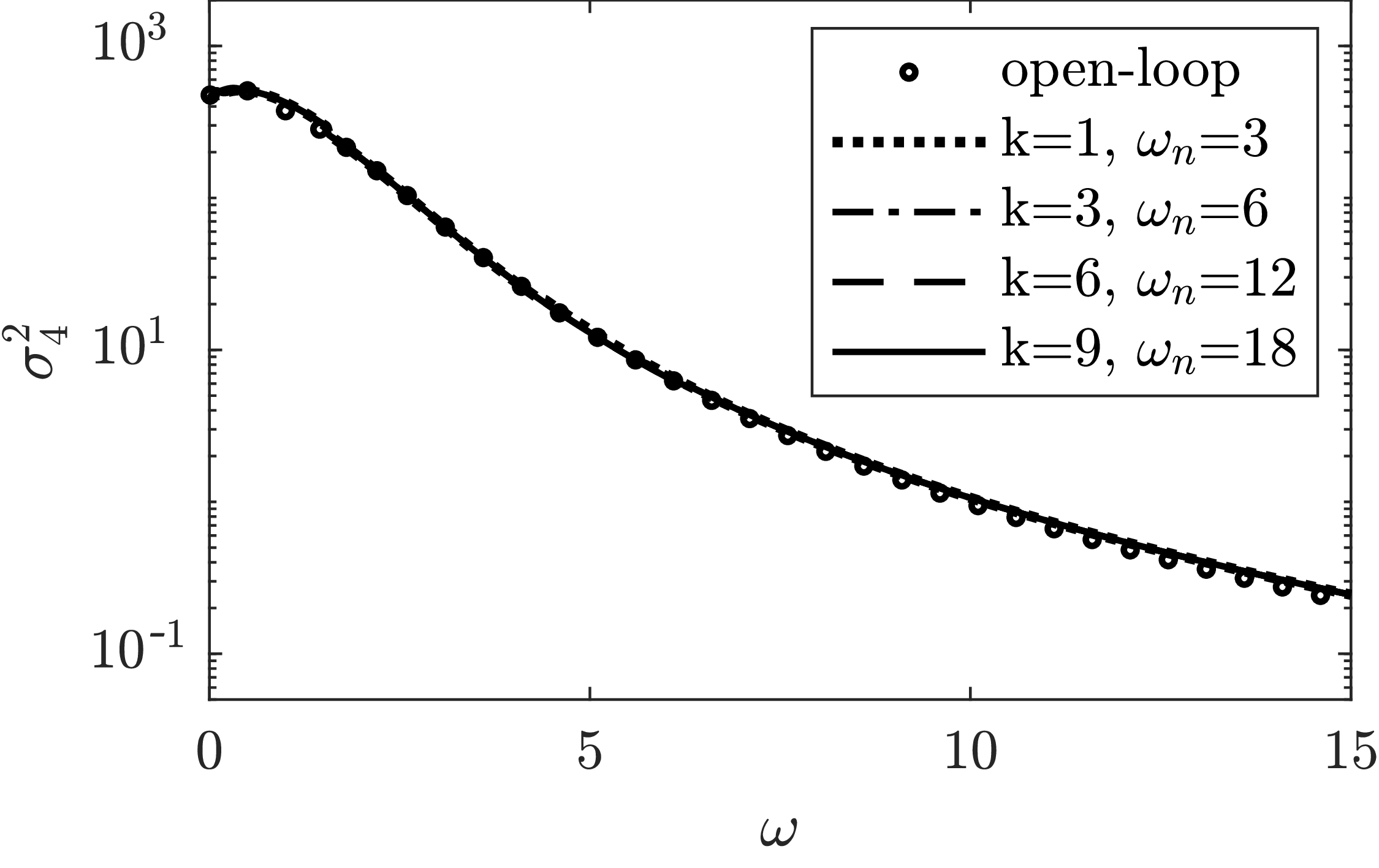}
    \llap{\parbox[b]{2.275in}{(c)\\\rule{0ex}{1.325in}}}
    \hspace{0mm}
    \includegraphics[width=0.47\textwidth]{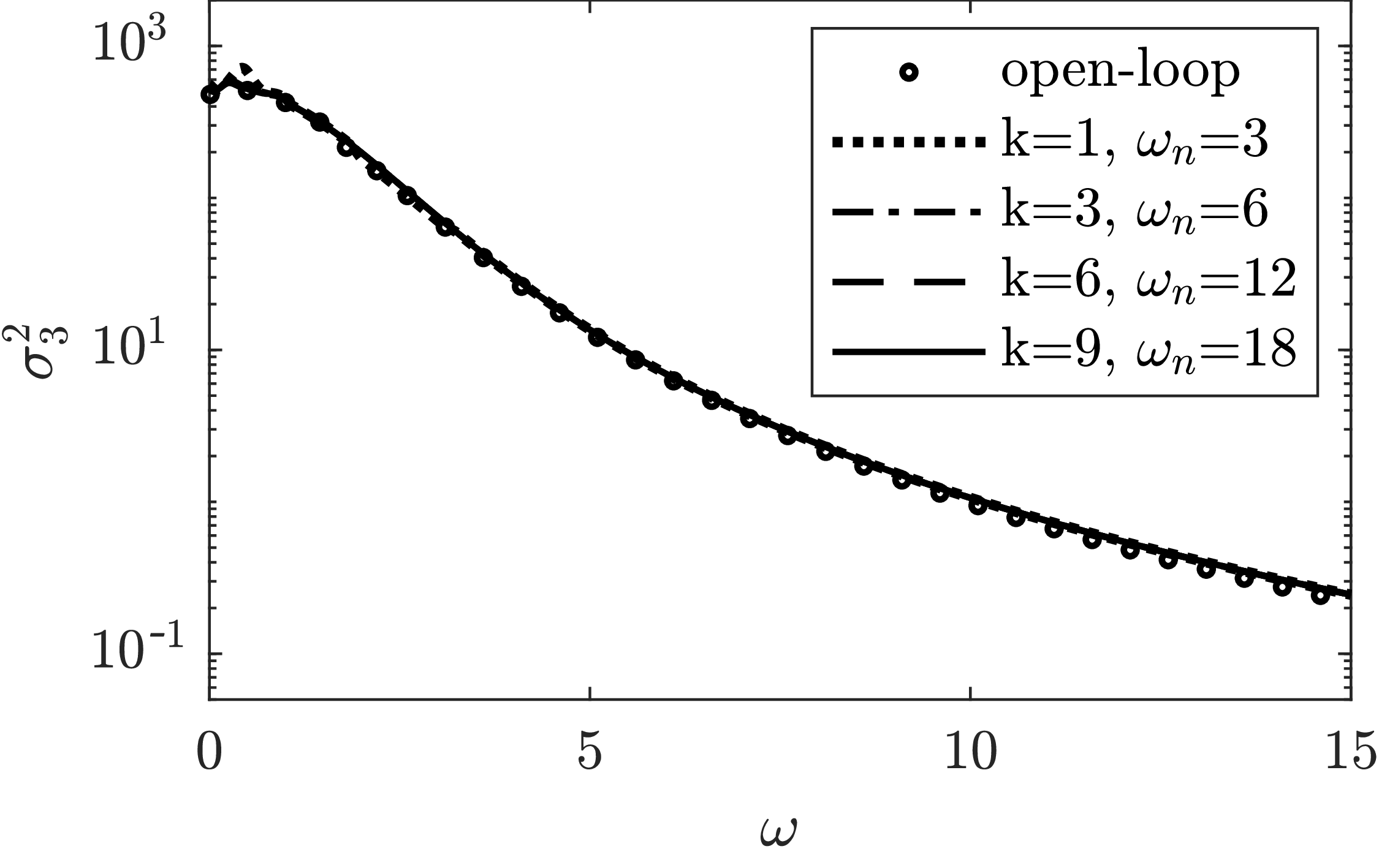}
    \llap{\parbox[b]{2.325in}{(d)\\\rule{0ex}{1.325in}}}
    }
    \caption[Comparison of the first four resolvent spectra between the open-loop system and the closed-loop error systems]{Comparison of (a) the first singular value $\sigma_1$, (b) the second singular value $\sigma_2$, (c) the third singular value $\sigma_3$ and (d) the fourth singular value $\sigma_4$ from resolvent analysis of the open-loop system ($\circ$) and the closed-loop systems $\textbf{Z}(s)$ (lines) at $\Rey=90$. Estimators are designed for a sensor placed at $(7.70,\ 0.73)$ using four different parameter choices.}
    \label{fig:lqe_valid_closedsigs}
\end{figure}

To quantitatively characterise the estimation performance, figure \ref{fig:lqe_valid_closedsigs} compares the first four resolvent spectra computed from the open-loop system (without estimator) to those from the closed-loop error system $\textbf{Z}(s)$ for the four cases.
A comparison of the four plots reveals two important facts: first, with increasing values of $k$ and $\omega_n$, the resolvent spectra of the closed-loop error system eventually converge, as shown by the matched dashed lines ($k=6,\ \omega_n=12$) and solid lines ($k=9,\ \omega_n=18$); second, the estimator significantly changes the first two resolvent spectra but barely modifies either the third or the fourth singular spectra. 
These two observations further validate the iterative method for solving the OE problem: it is unnecessary to optimise all energy gains over the whole frequency range and a reasonable choice of $k$ and $\omega_n$ should give convergence to the global optimal estimator. 
Note that $\gamma^2{(k,\omega_n)}$ is defined as the integral over the first $k$ resolvent spectra within the frequency range $[-\omega_n,\omega_n]$ whereas $\textbf{\textit{J}}_{\text{OE}}$ is over all resolvent spectra across all frequencies. 
The gap between the cost function $\gamma^2{(k,\omega_n)}^2$ and $\textit{\textbf{J}}_{OE}$ represents the effect of `background' or `free-stream' modes which do not affect the estimator design. In other words, the reduced disturbances used for the optimal estimator design procedure allow us to ignore these `background' modes (which are not important for estimation) and thus reduce the complexity of computation.

\subsubsection{Effect of sensor noise}
\begin{figure}
    \centerline{
    \includegraphics[width=0.47\textwidth]{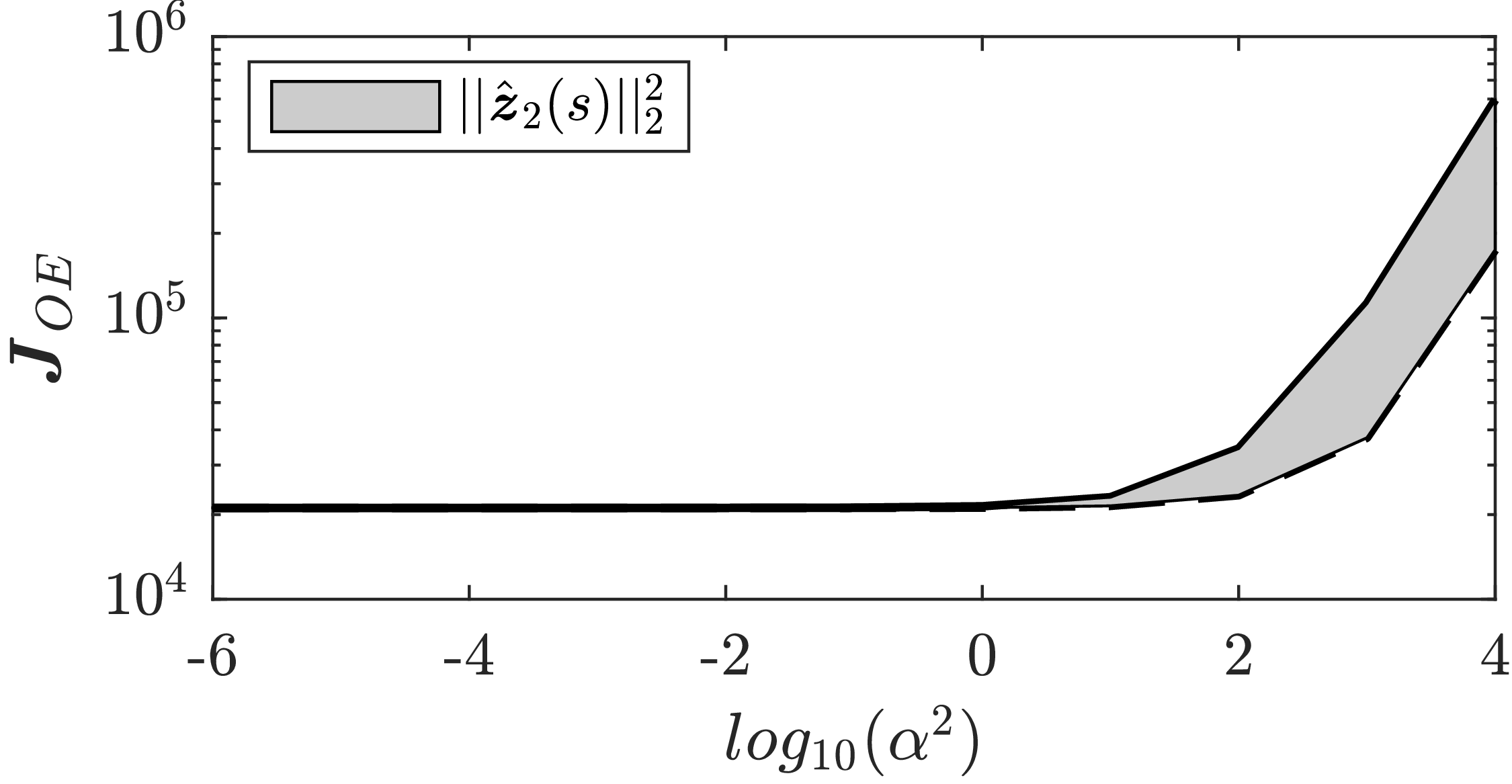}
    \llap{\parbox[b]{2.3in}{(a)\\\rule{0ex}{1.075in}}}
    \includegraphics[width=0.47\textwidth]{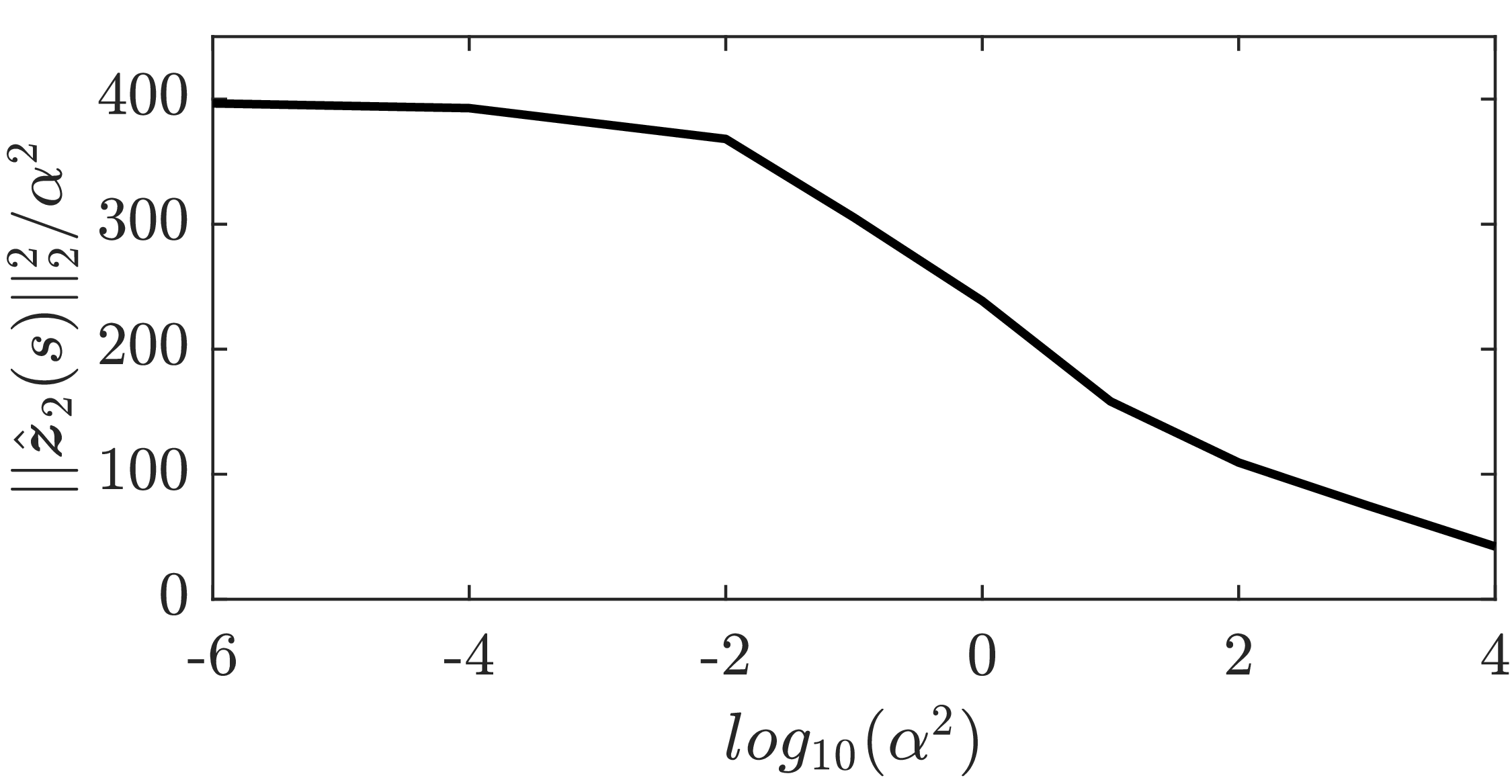}
    \llap{\parbox[b]{2.3in}{(b)\\\rule{0ex}{1.075in}}}
    }
    \caption[The performance of the optimal estimator designed for different magnitudes of the sensor noise]{The performance of the optimal estimator designed for different magnitudes of the sensor noise. (a) The total cost function $\textbf{\textit{J}}_{\text{OE}}$ as a function of the noise magnitude $\alpha$. The grey area indicates the contribution of the sensor noise ($||\textbf{Z}_2(s)||_2^2$) to the total cost function. (b) The mean energy of the estimation error due to unitary sensor noise.}
    \label{fig:lqe_noise}
\end{figure}

The choice of sensor noise also affects the optimal estimation performance. The magnitude of the sensor noise is determined by the matrix $\textbf{V}^{1/2}=\alpha\textbf{I}$ where the positive scalar $\alpha$ represents the relative size of the sensor noise $\textbf{\textit{n}}$ to the disturbances $\textbf{\textit{d}}$. To quantify their contribution to the total estimation error, we compute the $H_2$ norm squared of the corresponding transfer functions individually from numerical simulations. Here, let $\textbf{Z}_1(s)$ be the transfer function from the disturbances $\textbf{\textit{d}}$ to the estimation error, and let $\textbf{Z}_2(s)$ be the transfer function from the sensor noise $\textbf{\textit{n}}$ to the estimation error. It can be shown that $||\textbf{Z}_1(s)||_2^2+||\textbf{Z}_2(s)||_2^2=\textbf{\textit{J}}_{\text{OE}}$. The performance of the optimal estimator for different choices of $\alpha$ and the corresponding contributions from $\textbf{\textit{d}}$ and $\textbf{\textit{n}}$ are summarised in figure \ref{fig:lqe_noise}.

In figure \ref{fig:lqe_noise}(a), we can see that for $\alpha<1$, the estimation performance is insensitive to the precise choice of $\alpha$. With further increasing values of $\alpha$, the estimation error energy  $\textbf{\textit{J}}_{\text{OE}}$ increases logarithmically. The deterioration of the estimation performance is owing not only to the contribution of the sensor noise, as indicated by the grey area, but also to the contribution of disturbances (lower bound of the grey area). The value of $\alpha$ controls the balance between measurement uncertainty and model uncertainty in the estimation problem. The estimator tends to reduce the effect of the measurement uncertainty if it is designed using a larger value of $\alpha$, as shown by figure \ref{fig:lqe_noise}(b). The normalised energy norm $||\textbf{Z}_2(s)||_2^2/\alpha^2$, which reflects the effect of unitary measurement uncertainty on the estimator performance, decreases with increasing values of $\alpha$. As can be seen from the figure, under the same amount of measurement uncertainty, the optimal estimator designed using a larger value of $\alpha$ is more effective at reducing the estimation error due to sensor noise. In practice, one should choose $\alpha$ based on a trade-off between the uncertainty of the measurements and the uncertainty of the model.

\subsection{Optimal full-state information controller design}
Analogous to the OE problem, the FIC problem can be solved by following algorithm \ref{code:algorithmic1} and replacing $\textbf{Q}^{1/2}=\Tilde{\textbf{U}}\hspace{0mm}^{T}_m\textbf{E}$ in the Riccati equation. The number of outputs in the performance measure $\textbf{Z}_1(s)$ is thus reduced to the rank of the POD basis $m$. In this section, we consider the optimal full-state information controller design for the actuator placement $\textbf{\textit{x}}_a=(2.56,1.18)$ at $\Rey=90$. The control penalty is chosen to be $\beta=10^{-4}$ to allow the most aggressive controller that minimises the energy gains of the closed-loop system for the first $k$ resolvent output modes across the frequency range $[-\omega_n, \omega_n]$. The rank of the reduced performance measuring $m$, the cost function for the optimal control design $\gamma^2{(k,\omega_n)}$ and the resulting control performance $\textbf{\textit{J}}_{\text{FIC}}$ for different choices of parameters $k$ and $\omega_n$ are summarised in table \ref{tab:lqr_valid_optfic}.

\begin{table}
  \begin{center}
\def~{\hphantom{0}}
  \begin{tabular}{p{1cm}p{1cm}p{1cm}p{1cm}p{1.5cm}p{1.5cm}}
  \toprule
      Case& $k$& $\omega_n$ & $m$ & $\gamma^2{(k,\omega_n)}$& $\textbf{\textit{J}}_{\text{FIC}}$ \\[3pt]
    \midrule
       (a)& 1 & 3 & 25 &  1916.46 & 20950.07 \\
       (b)& 3 & 6 & 156 &  3416.41 & 20856.78 \\
       (c)& 6 & 12 & 467 &  4350.22 & 20840.40 \\ 
       (d)& 9 & 18 & 789 &  5132.85 & 20837.68 \\ 
    \bottomrule
  \end{tabular}
  \caption[A summary of results from the FIC problem with different parameter choices]{The rank of the reduced performance measuring $m$, the cost function $\gamma^2{(k,\omega_n)}$ for the optimal control design and the resulting control performance $\textbf{\textit{J}}_{\text{FIC}}$ are listed with different parameter choices ($k$, $\omega_n$). The optimal FIC controllers are designed at $\Rey=90$ with an actuator placed at $\textbf{\textit{x}}_a=(2.56,1.18)$ and the control penalty is chosen to be $\beta=10^{-4}$.}
  \label{tab:lqr_valid_optfic}
  \end{center}
\end{table}

As expected, the rank of the reduced performance measuring $m$ as well as the cost function $\gamma^2{(k,\omega_n)}$ increase with increasing values of $k$ and $\omega_n$ since more resolvent output modes and their energy gains are considered. We further perform numerical simulations of the closed-loop system to evaluate the original control performance $\textbf{\textit{J}}_{\text{FIC}}$ defined by \eqref{equ:lqr_jfic}. We observe that the control performance $\textbf{\textit{J}}_{\text{FIC}}$ slightly decreases with increasing values of $k$ and $\omega_n$ and eventually converges to a constant with a relative change around $10^{-4}$. For the cost function in the FIC problem, a root mean square value $\epsilon_{\text{FIC}}(x,y)$ can also be defined as in \eqref{equ:lqe_epsilon} (see Appendix \ref{sec:app.a}). In this case, $\epsilon_{\text{FIC}}(x,y)$ indicates the contribution of each single disturbance to the energy norm of the performance measure $\textbf{Z}(s)$. In other words, it shows the spatial distribution of the receptivity to disturbances applied in the closed-loop system and the region with larger values of $\epsilon_{\text{FIC}}$ has higher receptivity and is more sensitive to disturbances (worst disturbance rejection ability) under the control of the actuator.
\begin{figure}
    \centerline{
    \hspace{2mm}
    \includegraphics[width=0.46\textwidth,clip=true,trim=0mm 5mm 50mm 5mm]{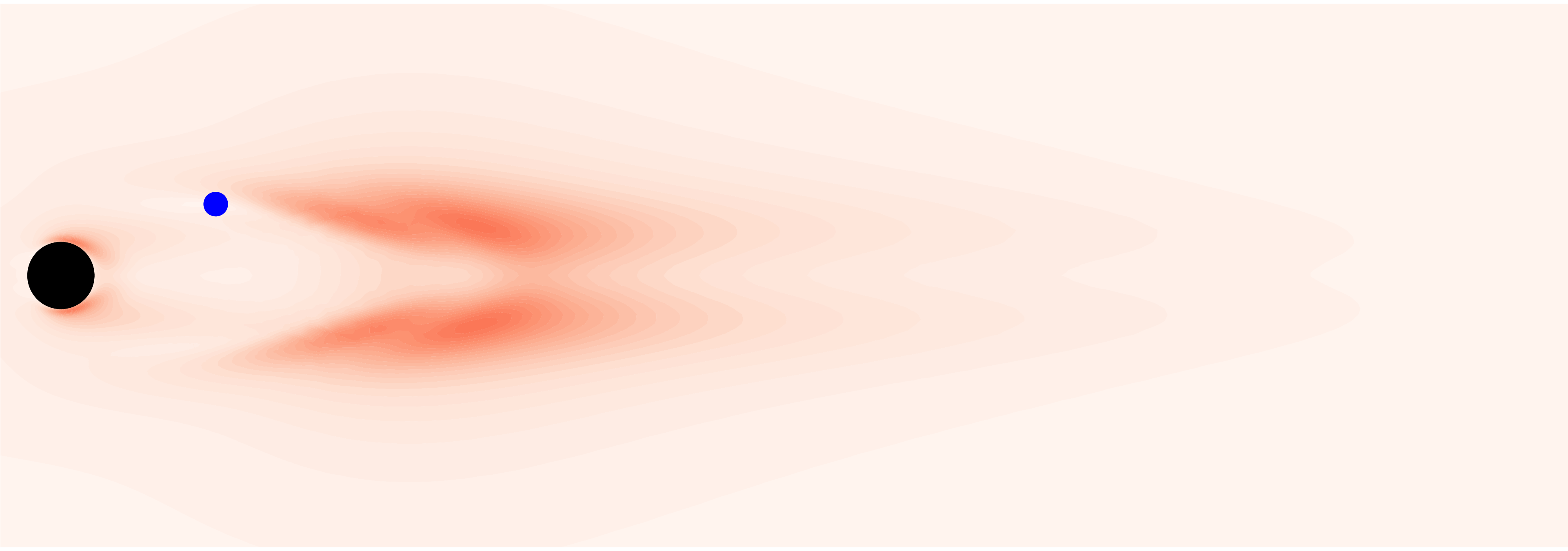}
    \llap{\parbox[b]{2.35in}{(a)\\\rule{0ex}{0.8in}}}
    \hspace{2mm}
    \includegraphics[width=0.46\textwidth,clip=true,trim=0mm 5mm 50mm 5mm]{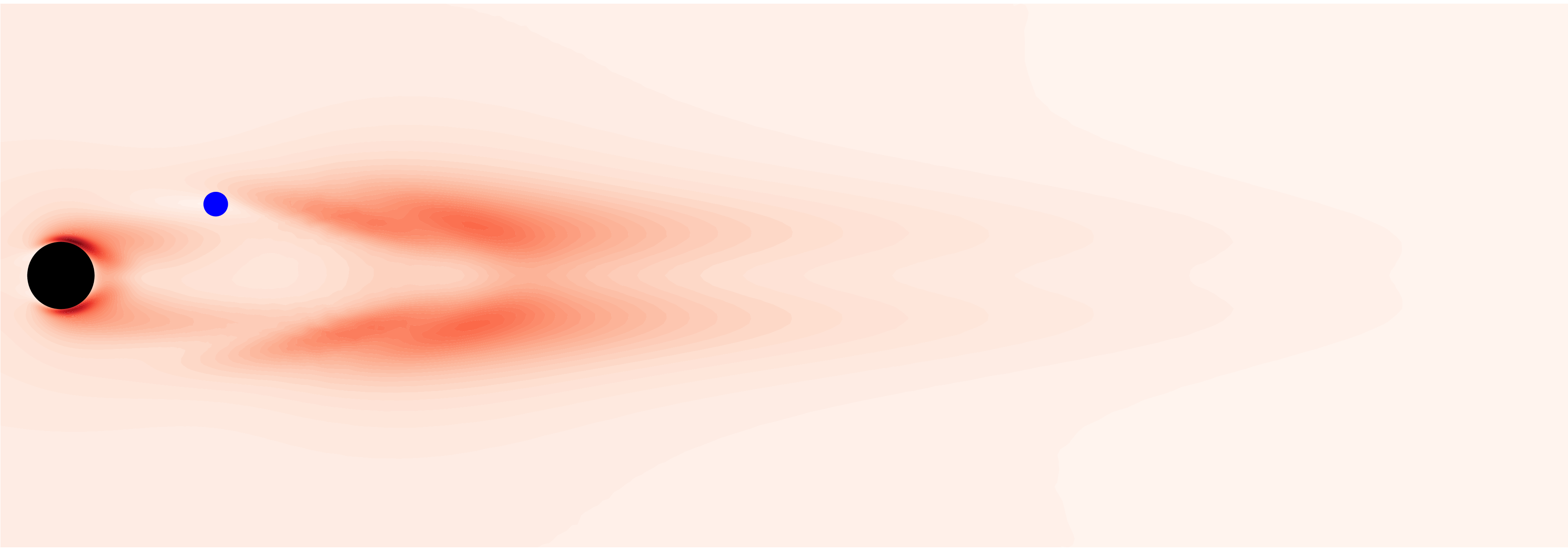}
    \llap{\parbox[b]{2.35in}{(b)\\\rule{0ex}{0.8in}}}
    }
    \centerline{
    \hspace{2mm}
    \includegraphics[width=0.46\textwidth,clip=true,trim=0mm 5mm 50mm 5mm]{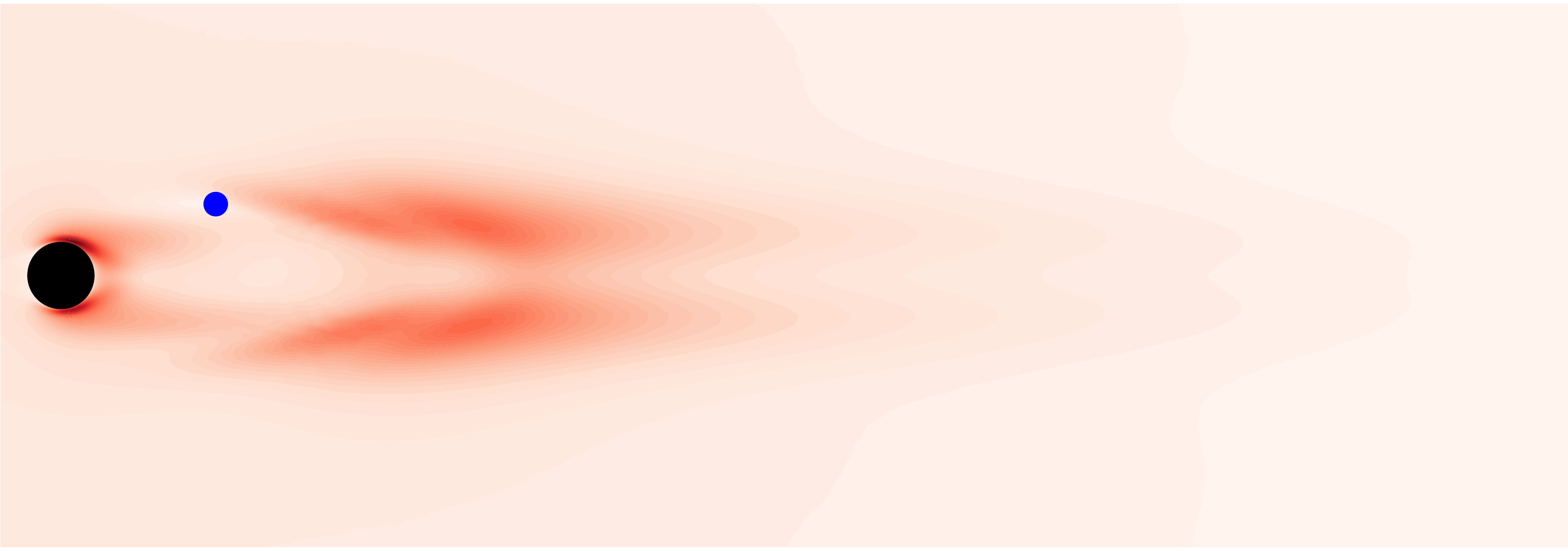}
    \llap{\parbox[b]{2.35in}{(c)\\\rule{0ex}{0.8in}}}
    \hspace{2mm}
    \includegraphics[width=0.46\textwidth,clip=true,trim=0mm 5mm 50mm 5mm]{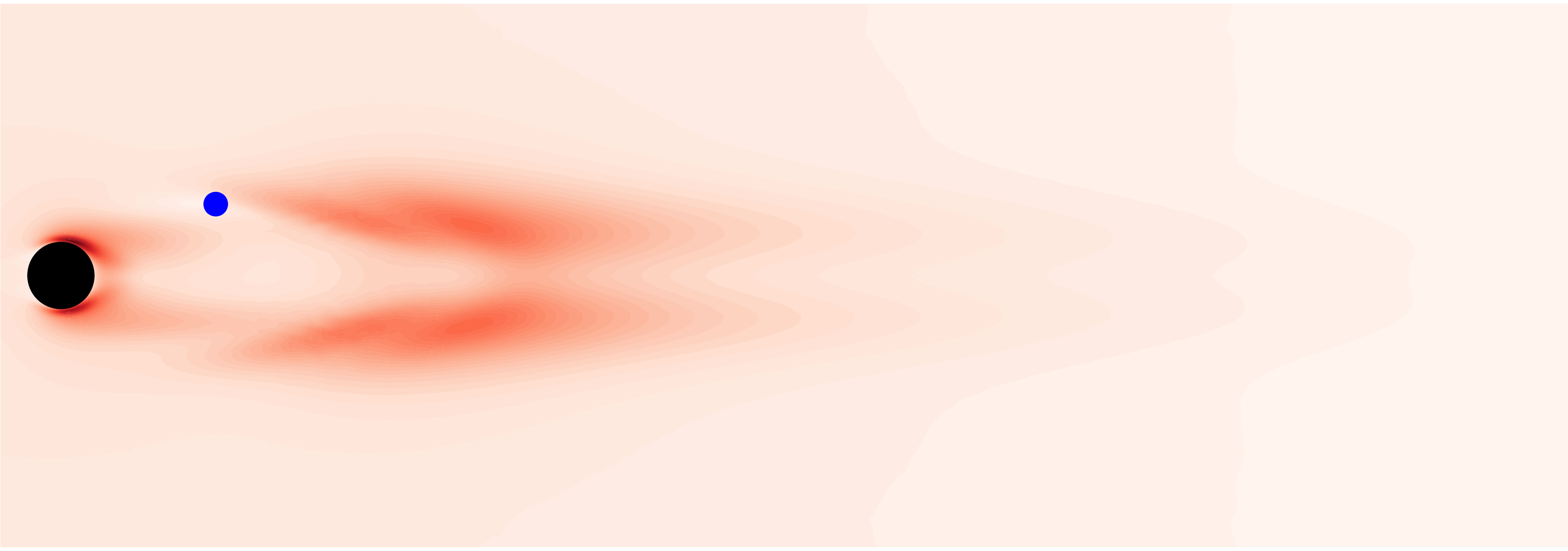}
    \llap{\parbox[b]{2.35in}{(d)\\\rule{0ex}{0.8in}}}
    }
    \centerline{
    \hspace{2mm}
    \includegraphics[width=0.46\textwidth,clip=true,trim=0mm 5mm 50mm 5mm]{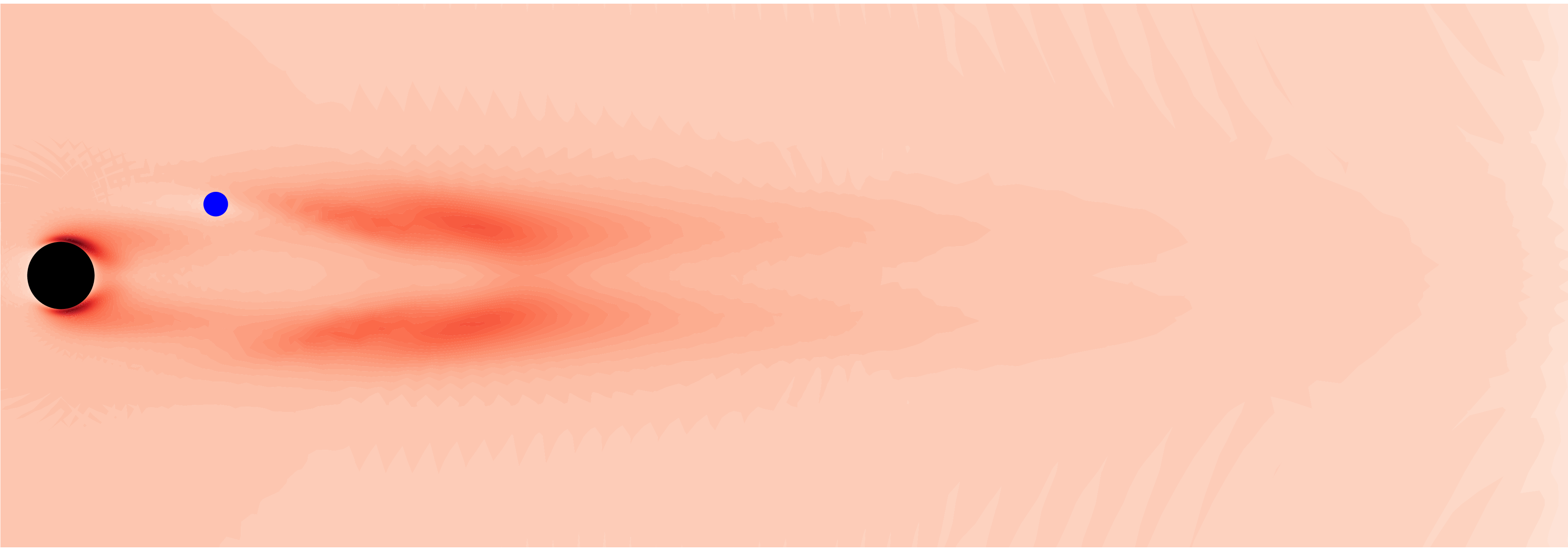}
    \llap{\parbox[b]{2.35in}{(e)\\\rule{0ex}{0.8in}}}
    \hspace{2mm}
    \includegraphics[width=0.46\textwidth,clip=true,trim=0mm 5mm 50mm 5mm]{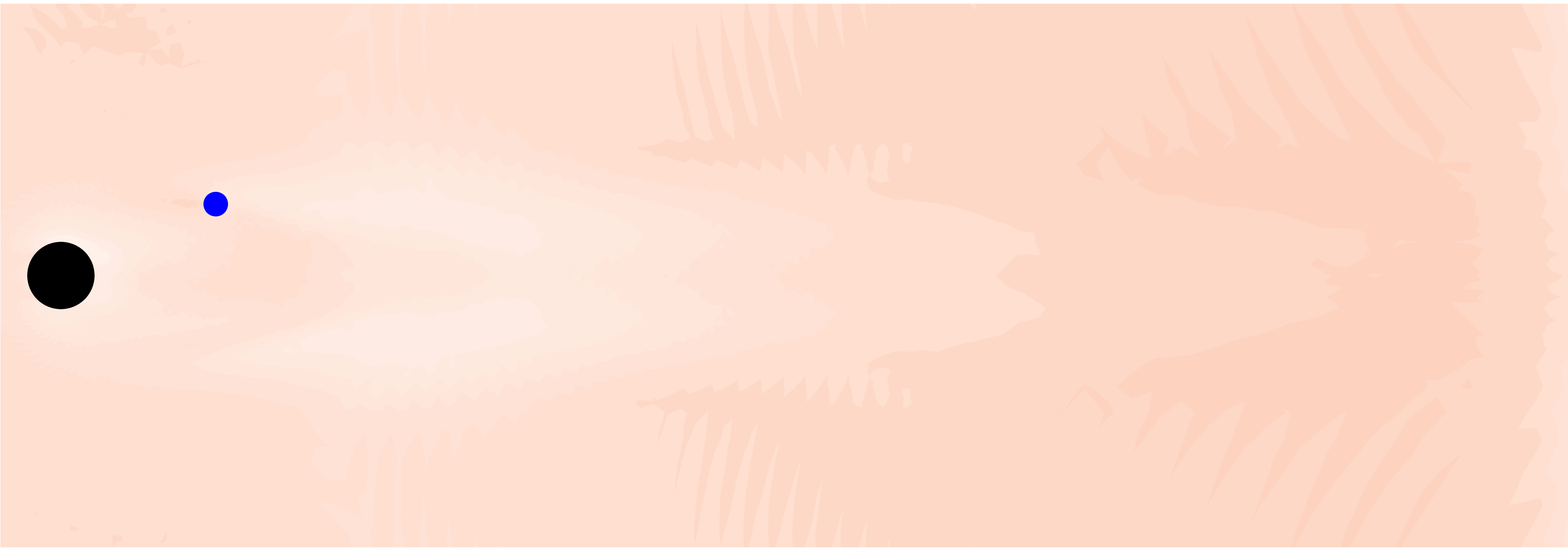}
    \llap{\parbox[b]{2.35in}{(f)\\\rule{0ex}{0.8in}}}
    }
    \caption[Spatial distribution of the receptivity to disturbances $\epsilon_{\text{FIC}}(x,y)$ for the closed-loop system with optimal FIC controllers designed using different parameters]{Spatial distribution of the receptivity to disturbances $\epsilon_{\text{FIC}}(x,y)$ under the control of the actuator placed at $(2.56,\ 1.18)$, which is marked by a blue dot ($\protect\bluedot$). The optimal full-state information controllers are designed using parameters (a) $k=1$,  $\omega_n=3$; (b) $k=3$, $\omega_n=6$; (c) $k=9$,  $\omega_n=12$ and (d) $k=9$, $\omega_n=18$. (e) The receptivity $\epsilon_{\text{FIC}}$ with full-rank performance measuring. (\textit{f}) The differences of the receptivity $\epsilon_{\text{FIC}}$ between (d) and (e) (i.e.~between $\gamma^2{(k,\omega_n)}$ and $\textbf{\textit{J}}_{\text{FIC}}$). All plots share the same linear colour scale $[0,\ 30]$.}
    \label{fig:lqr_valid_errdist}
\end{figure}

In figure \ref{fig:lqr_valid_errdist}(a, b, c, d), we plot the receptivity $\epsilon_{\text{FIC}}$ defined from the cost function $\gamma^2{(k,\omega_n)}$ for the four cases listed in table \ref{tab:lqr_valid_optfic}. We first notice that the receptivity of all four cases has almost the same spatial distribution regardless of the values of $k$ and $\omega_n$. In particular, the minimum of $\epsilon_{\text{FIC}}$ occurs near the actuator position where linear perturbations are best controlled. Another significant observation from these contour plots is that the spatial distribution of $\epsilon_{\text{FIC}}$ can also be divided into two regions: the first is near the cylinder and the second is in the far wake (downstream of the actuator). Analogous to the OE problem, the only difference among these plots is the slightly larger magnitude of $\epsilon_{\text{FIC}}$ for those cases with larger $k$ and $\omega_n$. This can be most clearly seen by comparing figure \ref{fig:lqr_valid_errdist}(a, d). In figure \ref{fig:lqr_valid_errdist}(e), we further consider the receptivity of the perturbation system controlled by the same controller designed for case (d) but with full-state performance measuring. That is, $\epsilon_{\text{FIC}}$ is defined from $\textbf{\textit{J}}_{\text{FIC}}$, which can be computed from numerical simulations of the corresponding adjoint system (see Appendix \ref{sec:app.a}). From a comparison of figure \ref{fig:lqr_valid_errdist}(d, e), their differences only occur in the `freestream' area and the receptivity defined from $\gamma^2{(k,\omega_n)}$ (with output reduction) exactly reproduces the characteristic structures found using $\textbf{\textit{J}}_{\text{FIC}}$ (without any reductions). This can be most clearly seen in figure \ref{fig:lqr_valid_errdist}(f) where regions with high values of $\epsilon_{\text{FIC}}$ are white (negligible difference between (d) and (e)) whereas the `freestream' area is of approximately uniform colour.

\begin{figure}
    \vspace{2mm}
    \centerline{
    \hspace{0mm}
    \includegraphics[width=0.47\textwidth]{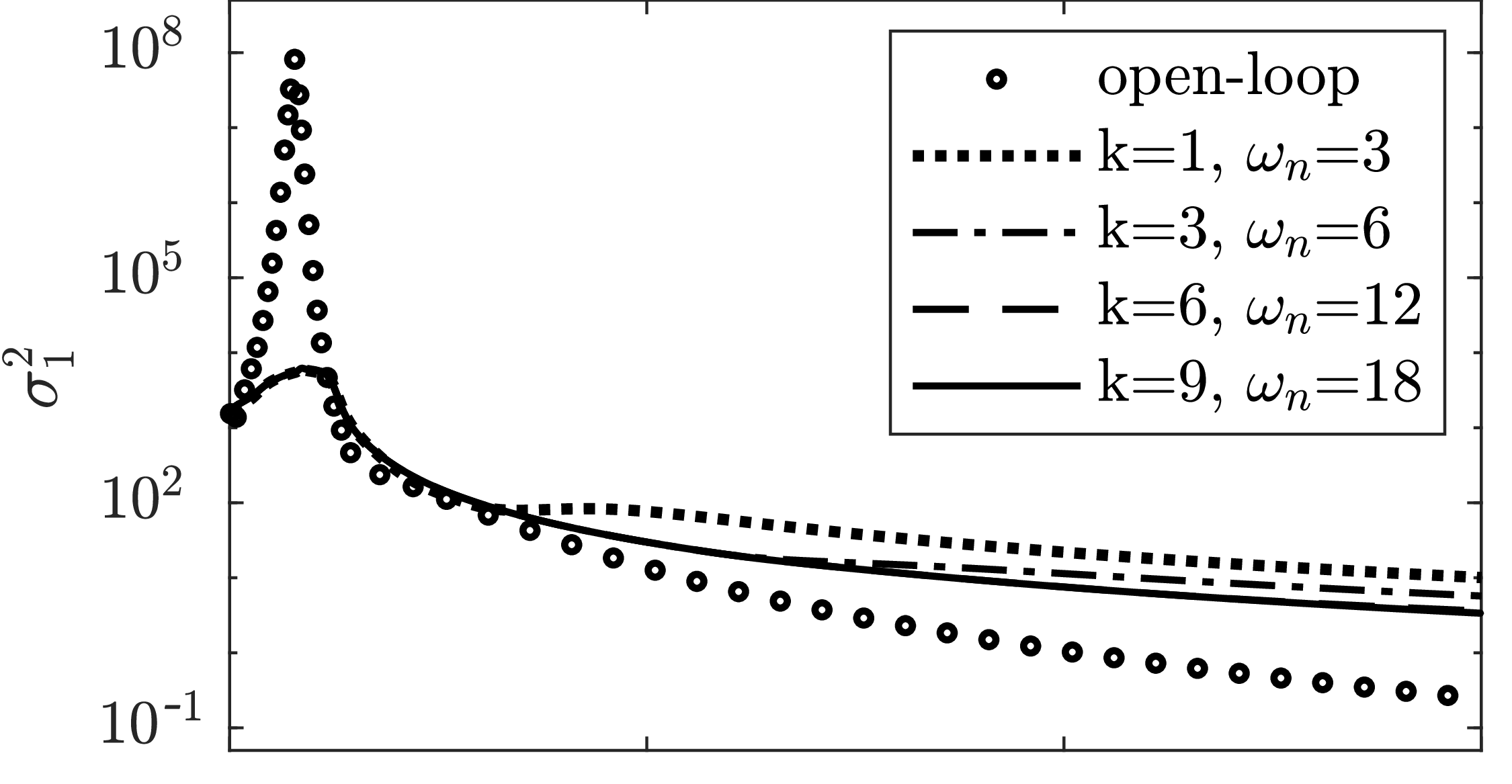}
    \llap{\parbox[b]{2.3in}{(a)\\\rule{0ex}{1.075in}}}
    \hspace{0mm}
    \includegraphics[width=0.47\textwidth]{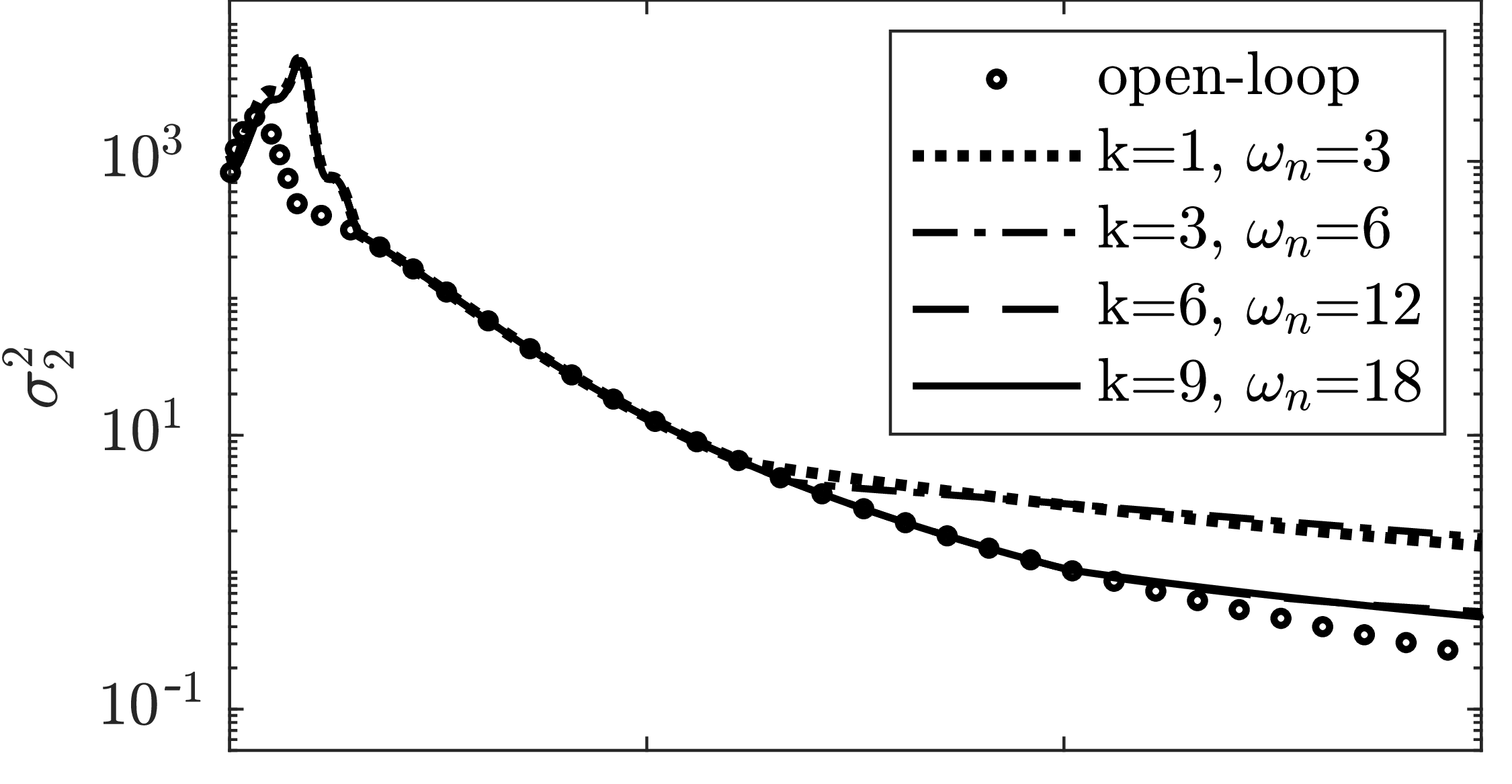}
    \llap{\parbox[b]{2.325in}{(b)\\\rule{0ex}{1.075in}}}
    }
    \vspace{0mm}
    \centerline{
    \hspace{0mm}
    \includegraphics[width=0.47\textwidth]{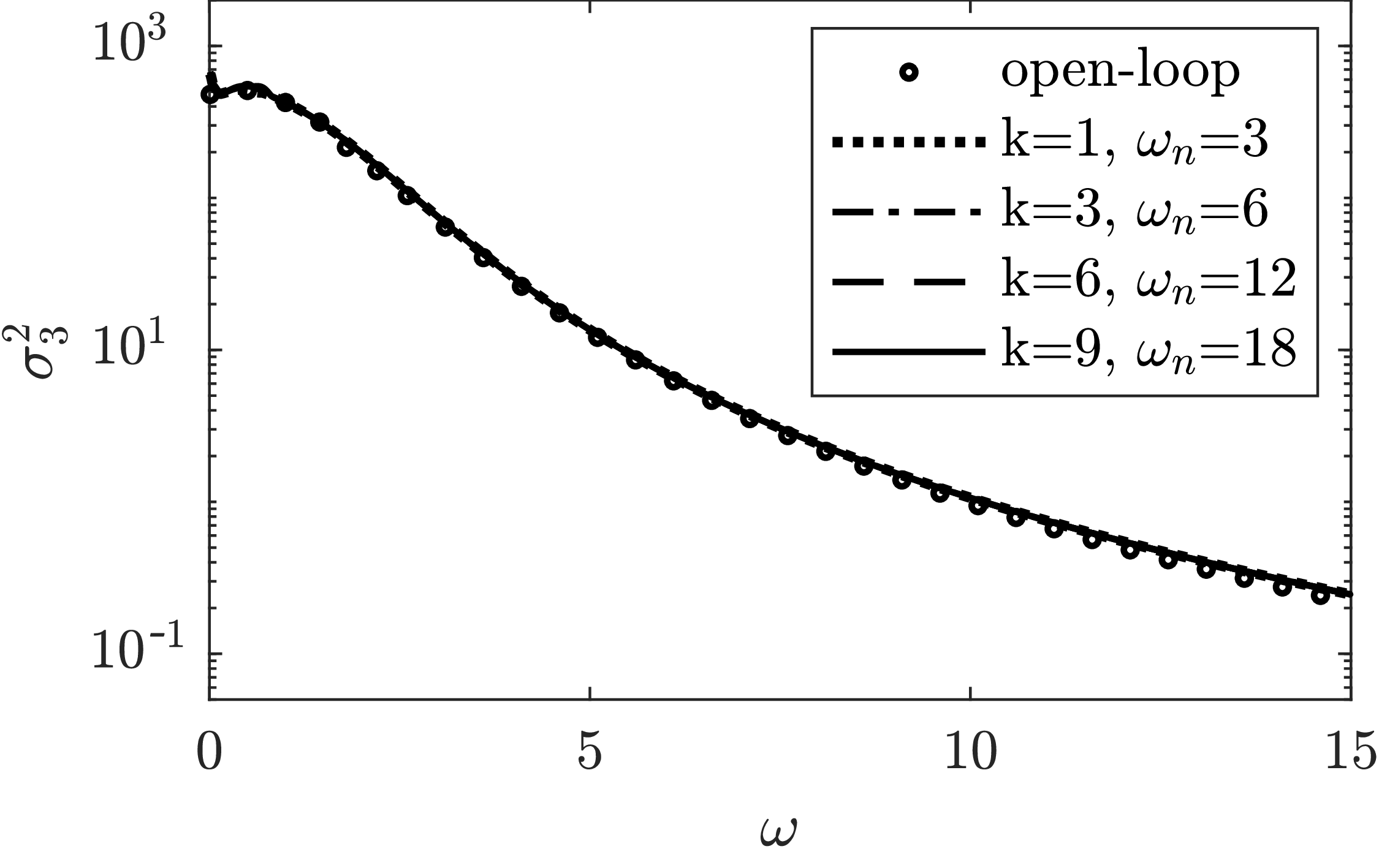}
    \llap{\parbox[b]{2.3in}{(c)\\\rule{0ex}{1.325in}}}
    \hspace{0mm}
    \includegraphics[width=0.47\textwidth]{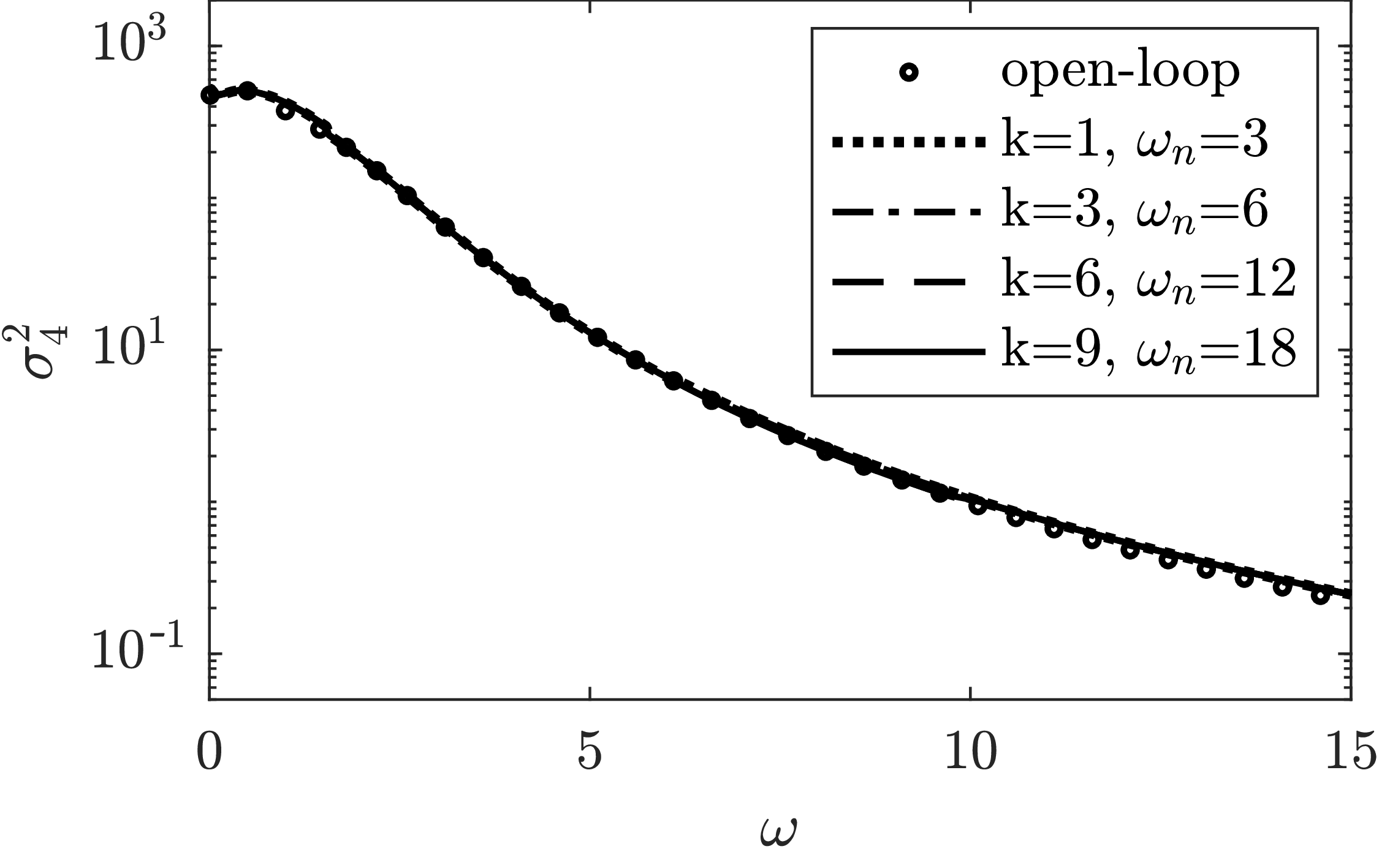}
    \llap{\parbox[b]{2.325in}{(d)\\\rule{0ex}{1.325in}}}
    }
    \caption[Comparison of the first four resolvent spectra between the open-loop system and the closed-loop systems]{Comparison of (a) the first singular value $\sigma_1$, (b) the second singular value $\sigma_2$, (c) the third singular value $\sigma_3$ and (d) the fourth singular value $\sigma_4$ from resolvent analysis of the open-loop system ($\circ$) and the closed-loop systems $\textbf{Z}(s)$ (lines) at $\Rey=90$. Controllers are designed for an actuator placed at $(2.56,\ 1.18)$ using four sets of parameters.}
    \label{fig:lqr_valid_closedsigs}
\end{figure}

Figure \ref{fig:lqr_valid_closedsigs} shows comparisons of the first four resolvent spectra between the open-loop system and the closed-loop system $\textbf{Z}(s)$ for the four cases listed in table \ref{tab:lqr_valid_optfic}. We observe that the resolvent spectra of the closed-loop system converge for larger values of $k$ and $\omega_n$, as shown by the matched dashed lines ($k=6,\ \omega_n=12$) and solid lines ($k=9,\ \omega_n=18$). In figure \ref{fig:lqr_valid_closedsigs}(c, d), the third and the fourth singular spectra of the closed-loop system are almost identical to those from the open-loop system regardless of the values of $k$ and $\omega_n$. 
Analogous to the OE problem, it is unnecessary to optimise all energy gains over the whole frequency range to solve the FIC problem, and a reasonable choice of $k$ and $\omega_n$ is sufficient to obtain a converged controller that achieves a globally optimal performance. 
The cost function $\gamma^2{(k,\omega_n)}$, which excludes the effect of `background' or `freestream' modes (that do not need to be controlled for good performance), reduces the complexity of the original FIC problem and makes feasible the optimal full-state information controller design procedure.

\subsubsection{Effect of control penalty}
\begin{figure}
    \centerline{
    \includegraphics[width=0.47\textwidth]{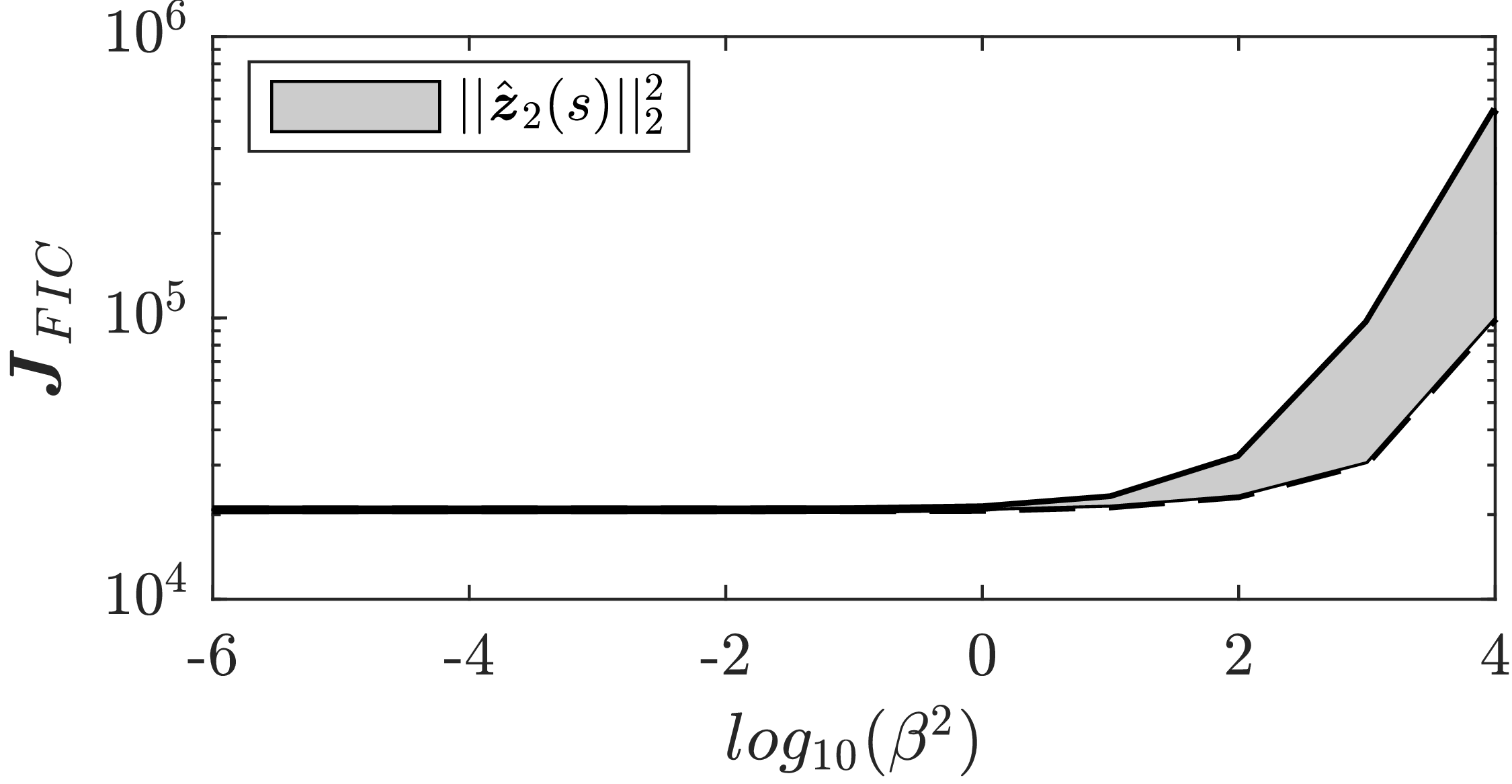}
    \llap{\parbox[b]{2.325in}{(a)\\\rule{0ex}{1.05in}}}
    \includegraphics[width=0.47\textwidth]{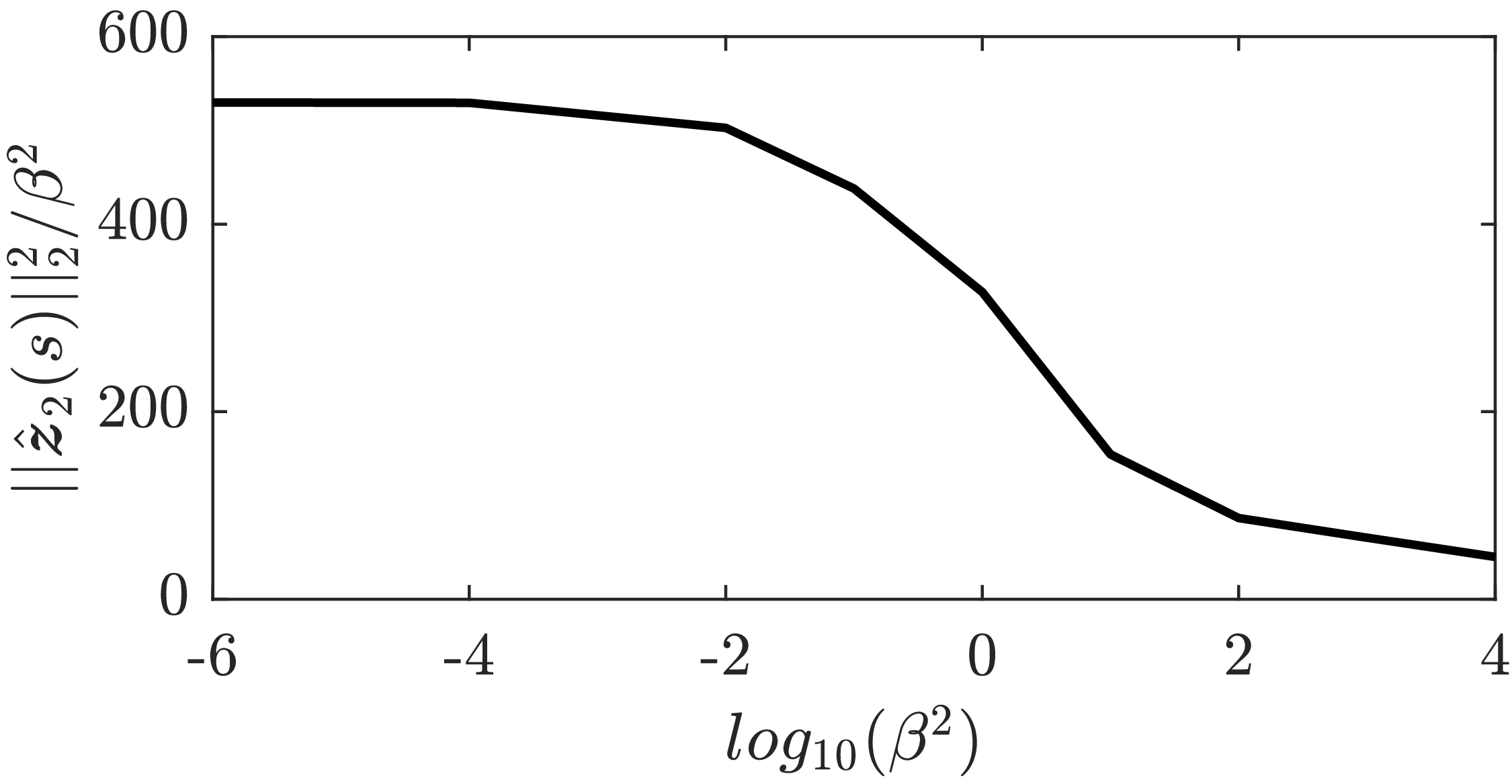}
    \llap{\parbox[b]{2.325in}{(b)\\\rule{0ex}{1.05in}}}
    }
    \caption[The performance of the optimal controller designed using different control penalty]{The performance of the optimal controller designed using different control penalties. (a) The total cost function $\textbf{\textit{J}}_{\text{FIC}}$ as a function of the control penalty $\beta$. The grey area indicates the contribution of the control cost ($||\textbf{Z}_2(s)||_2^2$) to the total cost function. (b) The normalised control cost as a function of the control penalty $\beta$.}
    \label{fig:lqr_controlcost}
\end{figure}

We now consider the effect of the control penalty on the optimal controller design. The control performance is composed of two parts: i) a contribution from the flow perturbations $\textbf{\textit{z}}_1=\textbf{Q}^{1/2}\textbf{\textit{w}}$ and ii) the control signal $\textbf{\textit{z}}_2=\textbf{R}^{1/2}\textbf{\textit{q}}$, where $\textbf{Q}^{1/2}=\textbf{M}^{1/2}\textbf{P}^T$ and $\textbf{R}^{1/2}=\beta\textbf{I}$. The balance between minimising perturbations $\textbf{\textit{w}}$ and minimising the control signal $\textbf{\textit{q}}$ is thus controlled by the positive scalar $\beta$. To quantify their contributions individually, the $H_2$ norm squared of the corresponding transfer functions $\textbf{Z}_1(s)$ and $\textbf{Z}_2(s)$ are computed from numerical simulations of the closed-loop system. It can be shown that $||\textbf{Z}_1(s)||_2^2+||\textbf{Z}_2(s)||_2^2=\textbf{\textit{J}}_{\text{FIC}}$. For controllers designed using different values of $\beta$, the contributions from $\textbf{\textit{w}}$ and $\textbf{\textit{q}}$ are shown in figure \ref{fig:lqr_controlcost} as well as the total control performance $\textbf{\textit{J}}_{\text{FIC}}$.

From figure \ref{fig:lqr_controlcost}(a), we observe that with a small value of $\beta<1$, the control performance is insensitive to the precise choice of $\beta$. This is consistent with previous findings for a spatially developing one-dimensional flow \cite{chen2011h}. A larger value of $\beta$ is not only less effective at reducing the perturbation magnitude (lower bound of the grey area) but also gives a control cost that increases logarithmically with $\beta$ (grey area). Both of these contributions result in the logarithmic deterioration of the control performance $\textbf{\textit{J}}_{\text{FIC}}$. A measure of the control signal is given by the normalised $H_2$ norm squared $||\textbf{Z}_2(s)||_2^2/\beta^2$ and a controller tends to stabilise the perturbation system with a smaller control signal if it is designed using a larger value of $\beta$, i.e.~less aggressive control. This can be most clearly seen in figure \ref{fig:lqr_controlcost}(b) where the norm squared of the control signal decreases with increasing values of $\beta$. In practice, one should choose $\beta$ based on a trade-off between minimising perturbations while maintaining reasonably-sized control inputs.

\subsection{Optimal feedback controller design}
We now turn our attention to the IOC problem where the optimal estimator and the optimal full-state information controller are combined to form a feedback controller as shown in figure \ref{fig:controlsetup}(c). In particular, we consider feedback control of the cylinder flow at $\Rey=90$ where a single sensor is placed at $\textbf{\textit{x}}_s=(7.7,0.73)$ with sensor noise of magnitude $\alpha=10^{-4}$ and a single actuator is placed at $\textbf{\textit{x}}_a=(2.56,1.18)$ with control penalty $\beta=10^{-4}$. 
In this case, we aim to minimise the energy gain of the performance measure $\textbf{Z}(s)$ defined in table \ref{tab:sysmats_summary} for the IOC problem. 
The cost function $\gamma^2{(k,\omega_n)}$ for the optimal feedback control design and the resulting control performance $\textbf{\textit{J}}_{\text{IOC}}$ for different choices of $k$ and $\omega_n$ are summarised in table \ref{tab:lqg_valid_optjioc}. 

\begin{table}
  \begin{center}
\def~{\hphantom{0}}
  \begin{tabular}{p{1cm}p{1cm}p{1cm}p{2.5cm}p{2.5cm}p{2.5cm}}
  \toprule
      Case& $k$& $\omega_n$ &$\gamma^2{(k,\omega_n)}(\times 10^5)$& $\textbf{\textit{J}}_{\text{IOC}}(\times 10^5)$\\[3pt]
      \midrule
       (a)& 1 & 3 & 1.370014 & 2.498071  \\
       (b)& 3 & 6 & 1.663140 & 2.027560  \\
       (c)& 6 & 12 & 1.759455 & 1.965255 \\ 
       (d)& 9 & 18 & 1.793839 & 1.958116  \\ 
       \bottomrule
  \end{tabular}
  \caption[A summary of results from the IOC problem with different parameter choices]{The cost function $\gamma^2{(k,\omega_n)}$ for the optimal feedback control design and the resulting feedback control performance $\textbf{\textit{J}}_{\text{IOC}}$ are listed for different parameter choices ($k$, $\omega_n$). The feedback controllers are found from solutions of the OE and FIC problems at $\Rey=90$ with a sensor placed at $\textbf{\textit{x}}_s=(7.70,0.73)$ and an actuator placed at $\textbf{\textit{x}}_a=(2.56,1.18)$. The sensor noise and control penalty are chosen to be $\alpha=\beta=10^{-4}$.}
  \label{tab:lqg_valid_optjioc}
  \end{center}
\end{table}

Similar to the OE and FIC problems, the cost function $\gamma^2{(k,\omega_n)}$ increases with increasing values of $k$ and $\omega_n$ for the IOC problem since more resolvent modes and their energy gains are considered. The feedback controllers designed for the four cases are implemented in numerical simulations with random disturbances applied everywhere and the control performance $\textbf{\textit{J}}_{\text{IOC}}$ is evaluated according to \eqref{equ:lqr_jfic}, which is summarised in table \ref{tab:lqg_valid_optjioc}. We first notice the convergence of $\textbf{\textit{J}}_{\text{IOC}}$ with increasing values of $k$ and $\omega_n$, and the final relative change decreases to around $10^{-3}$, which is larger than those observed in the OE and FIC problems. This is mainly caused by the iterative algorithm which leads to converged results but with small relative errors when compared to the true optimum. Moreover, the control performance rapidly deteriorates when the control setup is switched from the optimal actuation or sensing ($\textbf{\textit{J}}_{OE}$ or $\textbf{\textit{J}}_{FIC}$) to the feedback control problem ($\textbf{\textit{J}}_{\text{IOC}}$). In this case, the deterioration of control performance is due not only to the reduced ability of actuation and sensing, but also to the severe time delay between the actuator and the sensor \cite{oehler2018sensor,jin2019feedback}. 
\begin{figure}
    \vspace{2mm}
    \centerline{
    \hspace{0mm}
    \includegraphics[width=0.47\textwidth]{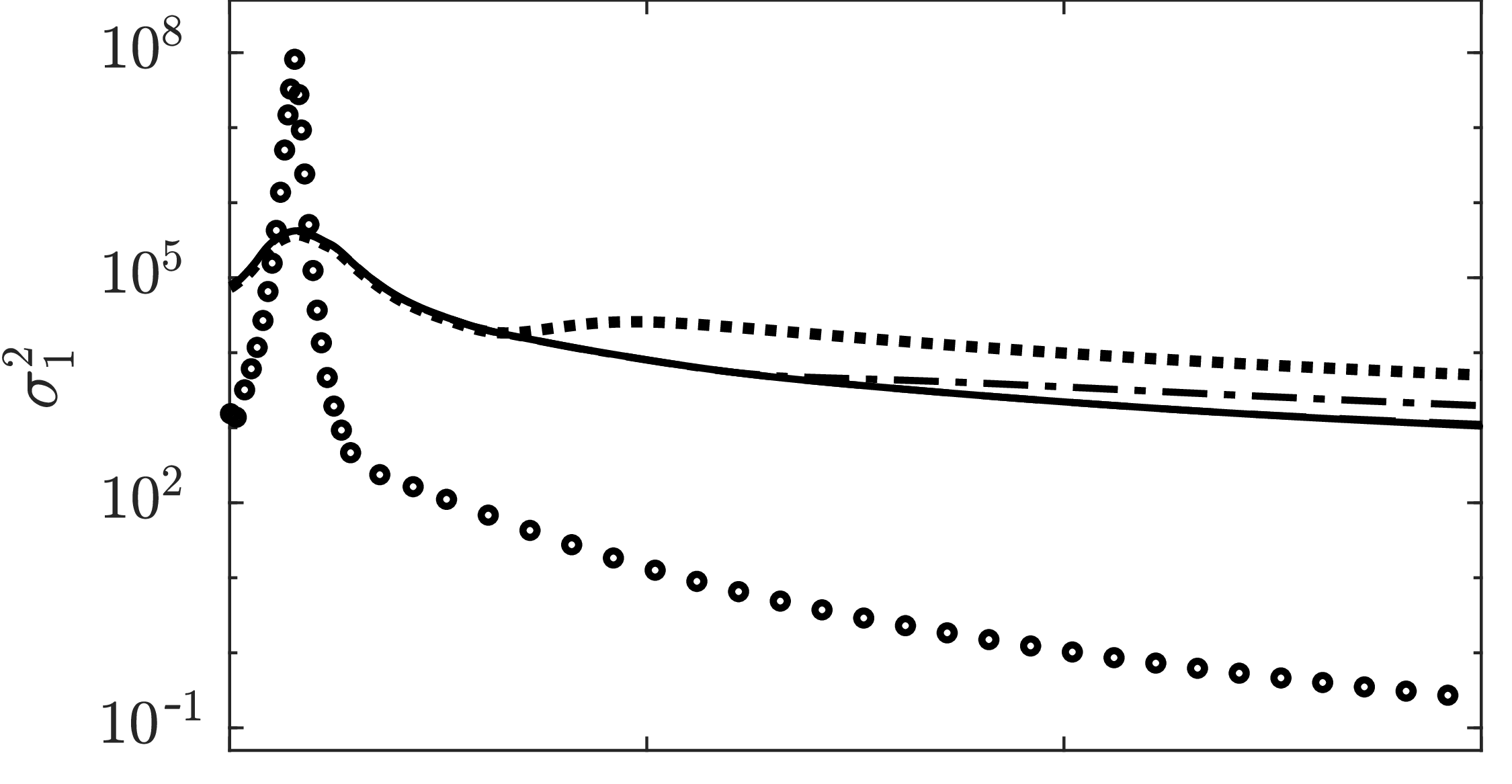}
    \llap{\parbox[b]{2.3in}{(a)\\\rule{0ex}{1.075in}}}
    \hspace{0mm}
    \includegraphics[width=0.47\textwidth]{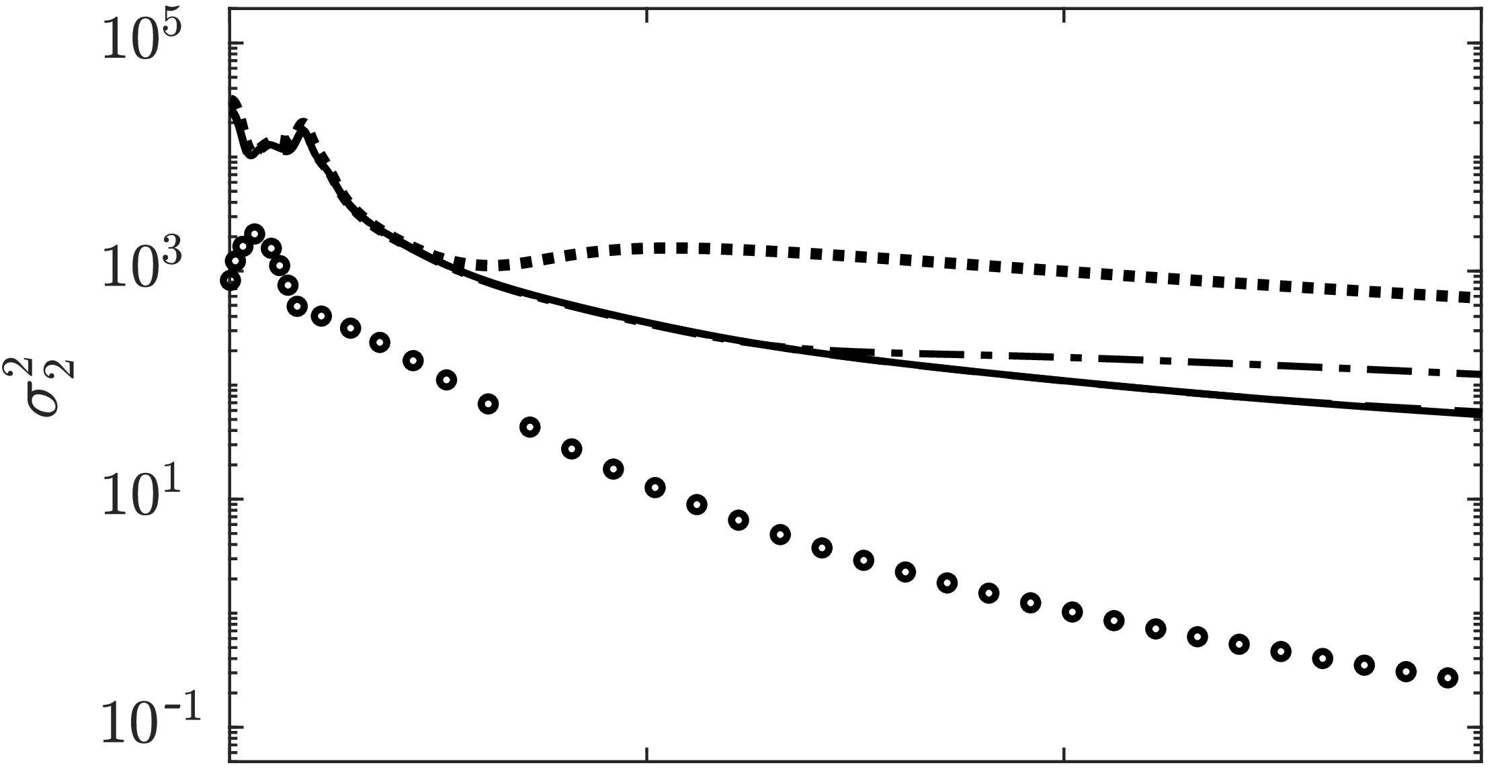}
    \llap{\parbox[b]{2.325in}{(b)\\\rule{0ex}{1.075in}}}
    }
    \vspace{0mm}
    \centerline{
    \hspace{0mm}
    \includegraphics[width=0.47\textwidth]{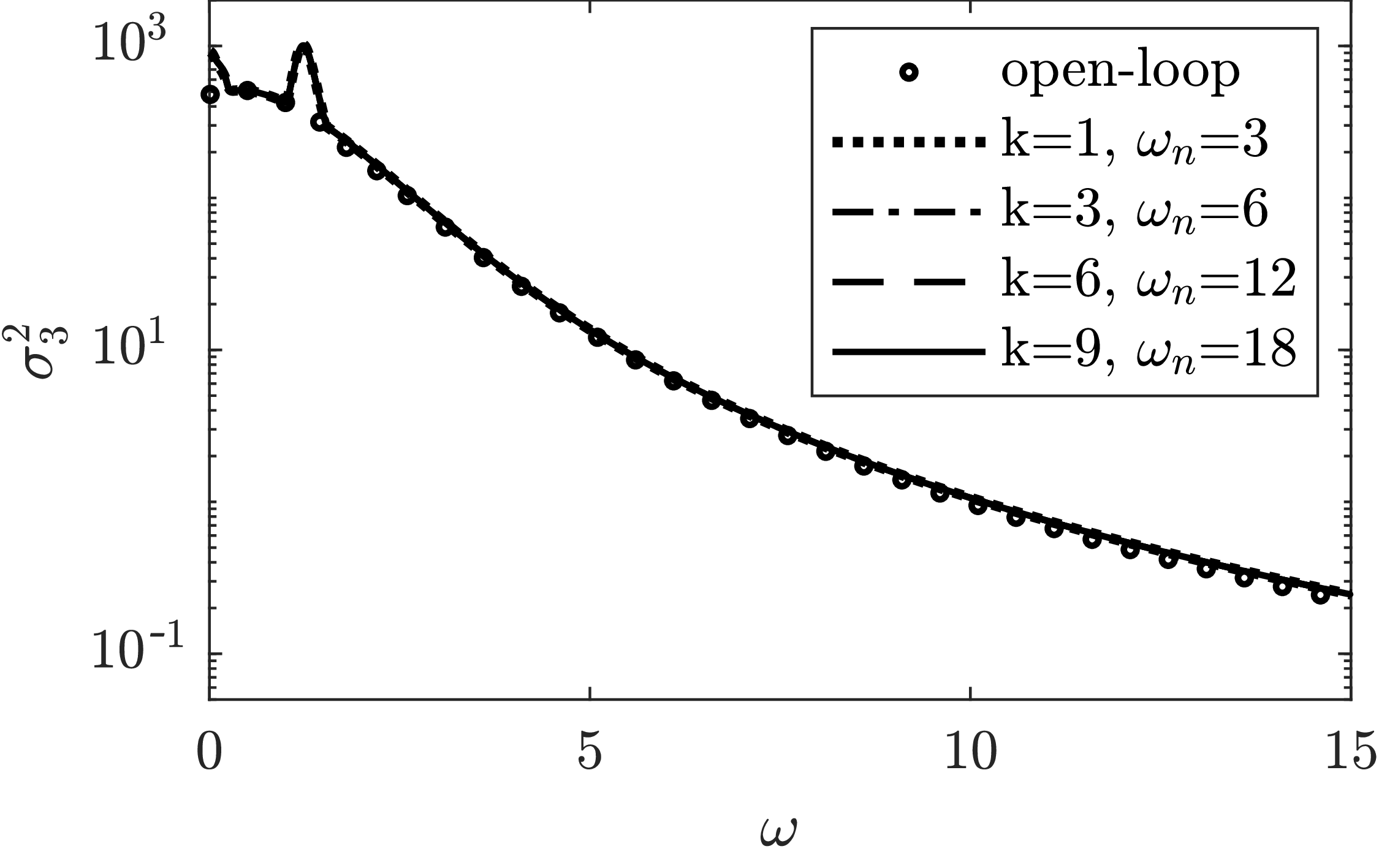}
    \llap{\parbox[b]{2.3in}{(c)\\\rule{0ex}{1.325in}}}
    \hspace{0mm}
    \includegraphics[width=0.47\textwidth]{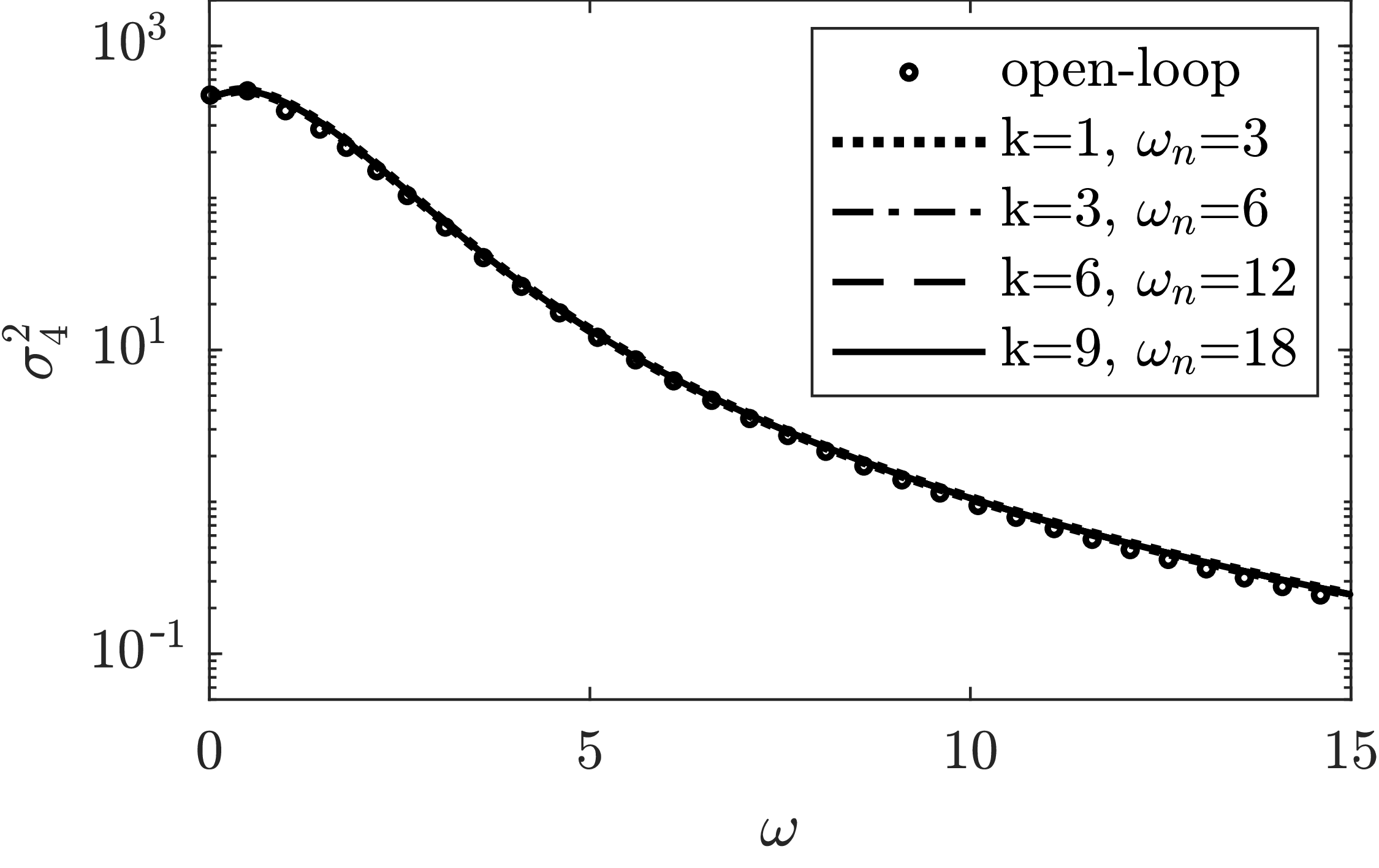}
    \llap{\parbox[b]{2.325in}{(d)\\\rule{0ex}{1.325in}}}
    }
    \caption[Comparison of the first four resolvent spectra between the open-loop system and the closed-loop systems]{Comparison of (a) the first singular value $\sigma_1$, (b) the second singular value $\sigma_2$, (c) the third singular value $\sigma_3$ and (d) the fourth singular value $\sigma_4$ from resolvent analysis of the open-loop system ($\circ$) and closed-loop systems $\textbf{Z}(s)$ (lines) at $\Rey=90$. The feedback controller is assembled from the optimal estimator and the optimal full-state information controller designed at $\Rey=90$ with a sensor placed at $\textbf{\textit{x}}_s=(7.70,0.73)$ and an actuator placed at $\textbf{\textit{x}}_a=(2.56,1.18)$.}
    \label{fig:lqg_valid_closedsigs}
\end{figure}

In this study, we characterise the control performance $\textbf{\textit{J}}_{\text{IOC}}$ using the energy gains between orthonormal inputs and outputs (i.e.~resolvent analysis) over the whole frequency space. Therefore, figure \ref{fig:lqg_valid_closedsigs} shows comparisons of the first four resolvent spectra between the open-loop system and the closed-loop system $\textbf{Z}(s)$ for the four cases listed in table \ref{tab:lqg_valid_optjioc}. Analogous to the OE and FIC problems, we first observe convergence of the resolvent spectra for larger values of $k$ and $\omega_n$, as shown by the matched dashed lines ($k=6,\ \omega_n=12$) and solid lines ($k=9,\ \omega_n=18$). Furthermore, in the first two resolvent spectra, i.e.~figure \ref{fig:lqg_valid_closedsigs}(a, b), we observe large separations between the open-loop system ($\circ$) and the closed-loop systems (lines). Although slight deviations are observed in the third resolvent spectra for $\omega<1.6$, figure \ref{fig:lqg_valid_closedsigs}(c, d) shows almost identical resolvent spectra for the open-loop system and the closed-loop systems regardless of the values of $k$ and $\omega_n$. Together with similar observations for the OE and FIC problems, these findings are consistent with those of \cite{dergham2013stochastic,jin2020energy} where two types of resolvent modes were identified: i) those corresponding to important physical mechanisms in the flow system that are required to be estimated and controlled and ii) those representing the advection and diffusion of perturbations in the freestream that can be ignored in the optimal control design procedure.

\section{Conclusions}\label{sec:conclusion}
We proposed a novel method to design $H_2$-optimal estimators and controllers for high-dimensional fluid flows. This involves solving high-dimensional Riccati equations with the full-rank covariance matrices (i.e.~$\textbf{W}$ and $\textbf{Q}$) which is challenging for existing numerical solvers. This challenge has been overcome by exploiting low-rank orthonormal bases that are constructed by performing a proper orthogonal decomposition of resolvent modes across a wide frequency range. We have thus transformed the problem from one of minimising energy gains between inputs and outputs over all frequencies and all possible directions to one of minimising the first few resolvent spectra over a specified frequency range. An iterative algorithm was implemented to update the low-rank orthonormal bases at each iteration using resolvent modes from the closed-loop system. Our numerical results indicate that by choosing sufficiently many resolvent modes over a sufficiently large frequency range, the performance converges to the global optimum in a few iterations. 

We demonstrated the effectiveness of the algorithm by considering the linear feedback stabilisation problem arising from the two-dimensional cylinder flow using a single sensor and a single actuator. The optimal feedback controller is solved by considering two basic problem: i) an optimal estimation (OE) problem and ii) a full-state information control (FIC) problem. In the OE problem, the optimal estimation performance is robust when the sensor noise is smaller than stochastic disturbances (i.e.~$\alpha<1$). Similarly, the optimal FI control performance is insensitive to the control penalty if we take minimising the perturbation magnitude as the priority (i.e~$\beta<1$).

The algorithm has been implemented for different choices of parameters, i.e.~the number of resolvent spectra $k$ and the frequency range $\omega_n$. The convergence of the corresponding results was shown by comparing the resolvent spectra for the original open-loop system and the closed-loop systems resulting from four different choices of $k$ and $\omega_n$. We observed that the optimal estimators and controllers affect only the first few resolvent spectra within a specified frequency range regardless of the parameter choice. This links directly to the findings of \cite{dergham2013stochastic,jin2020energy} where two types of resolvent modes were identified: i) those accounting for important physical mechanisms in the flow system that are required to be estimated and controlled and ii) those representing the advection and diffusion of perturbations in the freestream that can be ignored in the design procedure. We therefore see that, although the discretised Navier-Stokes equations are high-dimensional, only a limited number of physical mechanisms are important for effective estimation and control. 



\appendix
\section{Transfer functions for estimation and control}\label{sec:app.a0}
\paragraph{\textbf{Optimal estimation (OE).}} The transfer function governing the estimation performance $\textbf{Z}(s)$ is subject to disturbances and sensor noise: 
\begin{gather}\label{equ:app.OE_measure_zs}
    \begin{aligned}
        \textbf{Z}(s)=&
        \begin{bmatrix}
        \textbf{Z}_1(s)&\textbf{Z}_2(s)
        \end{bmatrix}\\
        =&\textbf{Q}^{1/2}[s\textbf{E}-\left(\textbf{A}-\textbf{K}_f\textbf{C} \right)]^{-1}
        \begin{bmatrix}
        \textbf{W}^{1/2}&-\textbf{K}_f\textbf{V}^{1/2}
        \end{bmatrix}\ ,
    \end{aligned}
\end{gather}
which consists of two subsystems: i) from disturbances $\textbf{\textit{d}}$ to the weighted estimation error $\textbf{Q}^{1/2}\textbf{\textit{e}}$; and ii) from the sensor noise $\textbf{\textit{n}}$ to the weighted estimation error $\textbf{Q}^{1/2}\textbf{\textit{e}}$, denoted by $\textbf{Z}_1(s)$ and $\textbf{Z}_2(s)$, respectively. 
\paragraph{\textbf{Full-state information control (FIC).}}
The transfer function measuring the control performance $\textbf{Z}(s)$ incorporates two subsystems: from disturbances $\textbf{\textit{d}}$ i) to the weighted perturbations $\textbf{Q}^{1/2}\textbf{\textit{w}}$ and ii) to the weighted control signal $\textbf{R}^{1/2}\textbf{\textit{q}}$, which are denoted by $\textbf{Z}_1(s)$ and $\textbf{Z}_2(s)$, respectively:
\begin{equation}\label{equ:app.FIC_measure_zs}
    \textbf{Z}(s)=
    \begin{bmatrix}
    \textbf{Z}_1(s)\\\textbf{Z}_2(s)
    \end{bmatrix}=
    \begin{bmatrix}
    \textbf{Q}^{1/2}\\
    -\textbf{R}^{1/2}\textbf{K}_r
    \end{bmatrix}    
    [s\textbf{E}-\left(\textbf{A}-\textbf{B}\textbf{K}_r \right)]^{-1}\textbf{W}^{1/2}\ .
\end{equation}
\section{The root mean square of the norm}\label{sec:app.a}
The root-mean-square value $\epsilon(x,y)$ is defined as:
\begin{equation}\label{equ:app.pngilon}
    \gamma^2(k,\omega)=\int_{\Omega}\epsilon^2(x,y)\ d\Omega\ ,
\end{equation}
where $\gamma^2(k,\omega)$ is the cost function of the associated design problem. In the OE problem, we can solve for $\epsilon_{\textrm{OE}}^2$ by using resolvent response modes or the time-series of the estimation error:
\begin{gather}
    \begin{aligned}
        \epsilon_{\textrm{OE}}^2=&\sum_{u,v}\Big\{\dfrac{1}{2\pi}\int_{-\omega_n}^{\omega_n}\sum_{i=1}^{k}\left(\sigma_i\hat{\textbf{u}}_i\right)^{\odot 2}\ d\omega\Big\}\\
        =&\lim_{t\rightarrow\infty}\dfrac{1}{T}\int_0^T\sum_{u,v}\Big\{\textbf{Q}^{-1}\left(\ \textbf{Q}^{1/2}\textbf{\textit{e}}\ \right)^{\odot 2}\Big\}\ dt\ ,
    \end{aligned}
\end{gather}
where the singular values $\sigma_i$ and resolvent response modes $\hat{\textbf{u}}_i$ are computed from resolvent analysis of the closed-loop system $\textbf{Z}(s)$. Here, $()^{\odot}$ is the Hadamard power and $\sum_{u,v}$ denotes summation over the streamwise and transverse components. In the FIC problem, a similar root-mean-square value $\epsilon_{\textrm{FIC}}$ can be defined from the resolvent forcing modes or results from numerical simulations of the adjoint system. 


%
\section*{Conflict of interest}
The authors declare that they have no conflict of interest.

\bibliographystyle{spmpsci}      
\bibliography{Bibliography}
\end{document}